\documentclass[notitlepage,10pt]{article}
\usepackage{amssymb,mathrsfs,amsmath,esint,enumitem,graphicx,color,comment}

\catcode`\@=11
\@addtoreset{equation}{section}

\catcode`\@=12

\newtheorem{Theorem}{Theorem}[section]
\newtheorem{Lemma}[Theorem]{Lemma}

\newtheorem{Remark}[Theorem]{Remark}

\def\sech{\mathrm{\hspace{1.5pt}sech\hspace{1.5pt}}}
\def\R{{\mathbb R}}
\def\S{{\mathbb S}}
\def\E{{\mathcal E}}
\def\F{{\mathcal F}}
\def\G{{\mathcal G}}
\def\L{{\mathcal L}}
\def\M{{\mathcal M}}
\def\N{{\mathcal N}}
\def\eps{\varepsilon}
\def\normal{N}
\def\Proof{\noindent\textit{Proof. }}
\def\qed{$~\square$\goodbreak \medskip}

\title{Embedded tori with prescribed mean curvature}
\author{Paolo Caldiroli\footnote{Dipartimento di Matematica, Universit\`a di Torino, via Carlo Alberto, 10 -- 10123 Torino, Italy. Email: \tt{paolo.caldiroli@unito.it}}, Monica Musso\footnote{Department of Mathematical Sciences, University of Bath, North Rd, Bath, BA2 7AY, 
United Kingdom. Email: \tt{m.musso@bath.ac.uk}}}

\date{}

\begin{document}
\maketitle

\begin{abstract}
\noindent
We construct a sequence of compact, oriented, embedded, two-dimensional surfaces of genus one into Euclidean 3-space with prescribed, almost constant, mean curvature of the form $H(X)=1+{A}{|X|^{-\gamma}}$ for $|X|$ large, when $A<0$ and $\gamma\in(0,2)$. Such surfaces are close to sections of unduloids with small neck-size, folded along circumferences centered at the origin and with larger and larger radii. The construction involves a deep study of the corresponding Jacobi operators, an application of the Lyapunov-Schmidt reduction method and some variational argument.
\smallskip

\smallskip

\end{abstract}

\section{Introduction and main result}

We are interested in the following problem: given a regular mapping $H\colon\R^{3}\to\R$ with suitable reasonable properties, look for
compact, oriented, embedded, two-dimensional surfaces $\Sigma$ of genus one into Euclidean 3-space, whose mean curvature at every point $X\in\Sigma$ equals $H(X)$. As prescribed mean curvature functions, we consider a class of radially symmetric mappings $H$ whose prototype is
\begin{equation}
\label{H-example}
H(X)=1+\frac{A}{|X|^{\gamma}}\quad\text{for $X\in\R^{3}$ with large $|X|$,}
\end{equation}
where $A\in\R$ and $\gamma>0$. The wanted surfaces are, roughly speaking, embedded tori like in figure 1, whose revolution axis is a circumference centered at the origin and with large radius, but with non-uniform, not necessarily circular, cross section.  Moreover, for (\ref{H-example}), such surfaces are asked to have almost constant mean curvature.

The interest toward this problem, stated also in a list of open problems raised by S.T. Yau in \cite{Yau82}, is motivated by the following remark: according to a famous result by A.D. Alexandrov \cite{Alexandrov1958}, the only embedded, oriented, compact, constant mean curvature (CMC) surfaces in $\R^{3}$ are round spheres. No CMC embedded tori in $\R^{3}$ exist, but only immersed tori, exhibited by H. Wente in his striking work \cite{Wente1986} of 1986. When $H$ is almost constant but non constant, it seems that the problem has never been considered up to now. Our goal is to show that for a class of prescribed mean curvature functions like (\ref{H-example}), which are suitable small perturbations of a constant, the issue of existence of embedded tori may radically change with respect to the case of constant mean curvature.
\bigskip

\centerline{\includegraphics[scale=0.2]{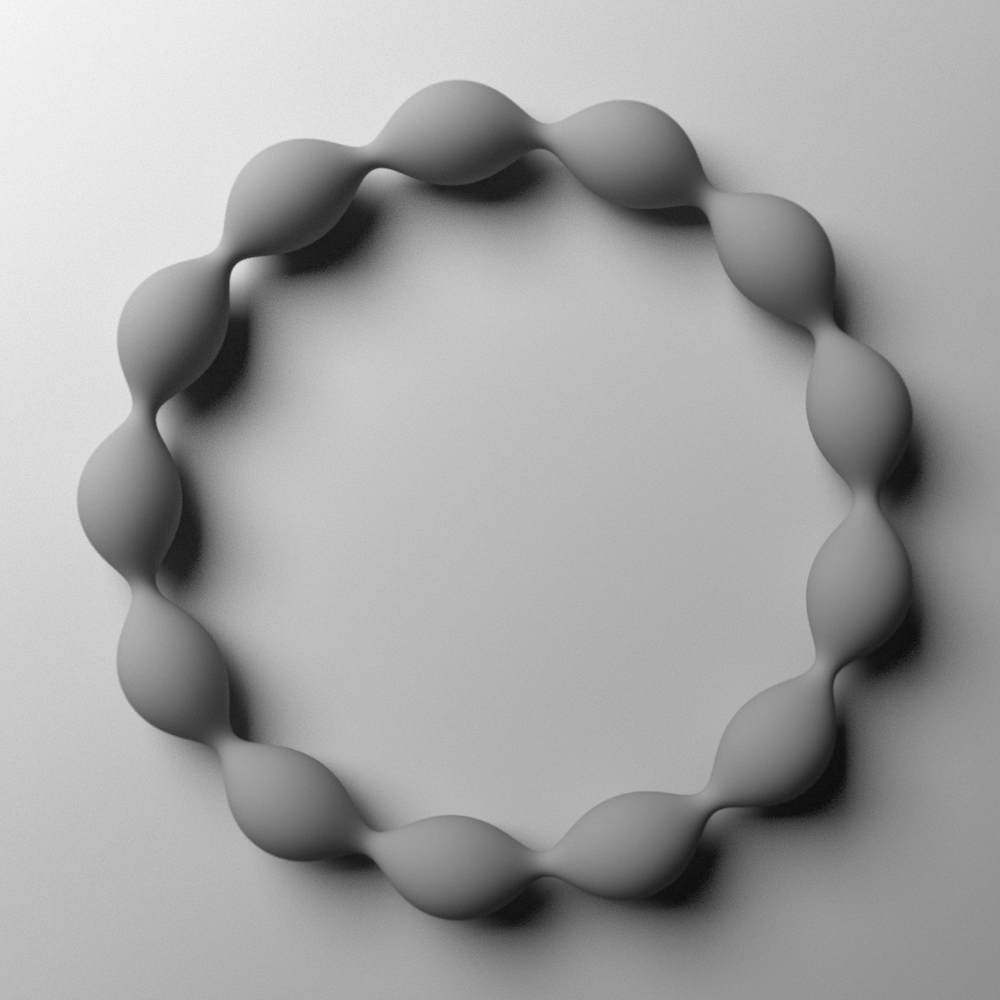}}
\centerline{\footnotesize{Fig.1. The reference surface (source: http://service.ifam.uni-hannover.de/$\sim$geometriewerkstatt/gallery/)}}
\bigskip

In order to state our result, let us start by constructing the reference surfaces. First, we introduce unduloids, i.e., CMC embedded surfaces of revolution in $\R^{3}$, generated by the trace of a focus of an ellipse which rolls without sliding on the axis of revolution. In particular, we denote $\Sigma_{a}$ the unduloid with mean curvature 1, axis of revolution on the $x_{3}$-axis, and neck-size $a\in\big(0,\frac{1}{2}\big]$ (the neck-size is defined as the distance of $\Sigma_{a}$ from the symmetry axis). It is convenient to consider conformal parametrizations of $\Sigma_{a}$, given by the mappings
\begin{equation}
\label{eq:Xa-vector}
X_{a}(t,\theta)=\left[\begin{array}{c}
x_{a}(t)\cos\theta\\
x_{a}(t)\sin\theta\\
z_{a}(t)\end{array}\right]
\end{equation}
with $x_{a}$ and $z_{a}$ solutions to some ode's (see Lemma \ref{L:c-a}).
%\begin{equation}\label{eq:xaza}
%x''_{a}=(1-2\gamma_{a})x_{a}-2x_{a}^{3}\quad\text{and}\quad z'_{a}=\gamma_{a}+x_{a}^{2}\text{~~where~~}\gamma_{a}:=a(1-a)~\!.
%\end{equation}
Functions $x_{a}$ and $z'_{a}$ turn out to be periodic of some period $2\tau_{a}>0$, with $\tau_{a}\sim-2\log a$ as $a\to 0$. Then for every $n\in\mathbb{N}$ sufficiently large, say $n\ge 4$, we introduce the mappings
\begin{equation}
\label{eq:toroidal-unduloid}
X_{n,a}(t,\theta):=\left[\begin{array}{c}x_{a}(t)\cos\theta\\ \left(\frac{nh_{a}}{\pi}+x_{a}(t)\sin\theta\right)\cos\frac{\pi z_{a}(t)}{nh_{a}}\\ \left(\frac{nh_{a}}{\pi}+x_{a}(t)\sin\theta\right)\sin\frac{\pi z_{a}(t)}{nh_{a}}\end{array}\right]
\end{equation}
where $2h_{a}$ equals the height of a complete period of the unduloid $\Sigma_{a}$.

The map ${X}_{n,a}$ turns out to be doubly-periodic with respect to a rectangle $[-n\tau_{a},n\tau_{a}]\times[-\pi,\pi]$ and provides a parametrization of an embedded toroidal surface $\Sigma_{n,a}$, obtained
by cutting a section of the unduloid $\Sigma_{a}$ made by $n$ periods, hence with length $2nh_{a}$, and folding it along a circumference  centered at the origin and with length $2nh_{a}$, as in figure 1.

We look for surfaces with some discrete symmetry and which are normal graphs over $\Sigma_{n,a}$.
%, namely, admitting a regular parametrization of the form $X_{n,a}+\varphi{\normal}_{n,a}$ with $\varphi\colon\R^{2}\to\R$ doubly-periodic with respect to the rectangle $[-\tau_{a},\tau_{a}]\times[-\pi,\pi]$.
This corresponds to solve the equation
\begin{equation}\label{main-equation}
\mathfrak{M}(X_{n,a}+\varphi{\normal}_{n,a})=H(X_{n,a}+\varphi{\normal}_{n,a})
\end{equation}
where $\mathfrak{M}$ is the mean curvature operator (defined in (\ref{eq:mean-curvature-def})), $N_{n,a}$ is the Gauss map corresponding to $X_{n,a}$ (defined as in (\ref{eq:normal-definition})), and the mapping $\varphi\colon\R/_{2\tau_{a}}\times\R/_{2\pi}\to\R$ is the unknown.
Our main result can be stated as follows:

\begin{Theorem}
\label{T:main}
Let $H\colon\R^{3}\to\R$ be a radially symmetric mapping of class $C^{2}$ satisfying:
\begin{itemize}
\item[$(H_{1})$]
$H(X)=1+A|X|^{-\gamma} +O(|X|^{-\gamma-\beta})$ as $|X|\to\infty$, with $A\in\R$, $0<\gamma<2$ and $\beta>0$;

\item[$(H_{2})$]
$|H''(|X|)|\le C|X|^{-\gamma-2}$ for some constant $C\in\R$ and $|X|$ large.
\end{itemize}
If $A<0$ then there exist $n_{A}\in\mathbb{N}$, a sequence $a_{n}\in(0,a_{0})$ and a sequence of regular, doubly-periodic mappings $\varphi_{n}\colon\R/_{2\tau_{a_{n}}}\times\R/_{2\pi}\to\R$ satisfying
\begin{gather}
\label{anbn}
a_{n}=\frac{b_{n}}{n^{\gamma}\log n}\text{~~with~~}\underbar{b}\le b_{n}\le\overline{b}~\!,\quad\text{for some $\underbar{b},\overline{b}>0$,}
\\
\label{varfin}
|\varphi_{n}(t,\theta)|\le\frac{C(\sech t)^{\widetilde{\mu}}}{n^{\min\{1,\gamma\}}}
\quad\forall(t,\theta)\in[-\tau_{a_{n}},\tau_{a_{n}}]\times[-\pi,\pi]\quad\text{for some $\widetilde{\mu}\in(0,1)$ and $C>0$,}
\end{gather}
and such that for every $n\ge n_{A}$ the map $X_{n,a_{n}}+\varphi_{n}N_{n,a_{n}}$ is an immersion of $\R^{2}$ into $\R^{3}$, doubly-periodic with respect to the rectangle $[-n\tau_{a_{n}},n\tau_{a_{n}}]\times[-\pi,\pi]$; the surface $\widetilde{\Sigma}_{n}$ parametrized by $X_{n,a_{n}}+\varphi_{n}N_{n,a_{n}}$ has mean curvature $H$, is symmetric with respect to the horizontal plane and invariant with respect to rotations of an angle $2\pi/n$ about the vertical axis. Moreover $n_{A}\to\infty$ as $A\to 0$ and for $n$ large enough $\widetilde{\Sigma}_{n}$ is embedded.
\end{Theorem}

Properties (\ref{anbn})--(\ref{varfin}) say that, when $n\to\infty$, the immersion $X_{n,a_{n}}+\varphi_{n}N_{n,a_{n}}$ parametrizes an embedded toroidal surface shaped on a collar made by $n$ almost tangent, almost round spheres joined by more and more narrow necks. This is consistent with a very recent result by Ciraolo and Maggi \cite{CiraoloMaggi2017} stating that almost constant mean curvature compact embedded surfaces are close, in a suitable sense, to a finite family of disjoint tangent round spheres with equal radii. The restriction on $\gamma$ means that the ``perturbation'' $A|X|^{-\gamma}$ has to be strong enough to make possible the presence of a sequence of larger and larger tori, and this gives clear evidence of a situation with a lack of compactness. The restriction on the sign of $A$ plays an essential role both for the existence of surfaces close to toroidal unduloids and for the embeddedness property.

The proof of Theorem \ref{T:main} is an application of the Lyapunov-Schmidt reduction method, combined with some variational argument. More precisely, we study the operator obtained by linearizing problem (\ref{main-equation}), for $a>0$ fixed, in the limit $n\to\infty$. This corresponds to study the Jacobi operator $\mathfrak{L}_{a}$ about right unduloids. Such operators have a non trivial kernel which, in view of the symmetries of our problem, is two-dimensional. For a suitable choice of functional spaces, after detaching the kernel, roughly speaking, we can invert the Jacobi operator and convert the equation (\ref{main-equation}) into a fixed point problem, where the contraction principle can be applied. Actually, because of the two-dimensional kernel of $\mathfrak{L}_{a}$, we arrive to solve (\ref{main-equation}) for every $a>0$ small enough and for every $n\in\mathbb{N}$ large, apart from a couple of Lagrange multipliers. In order to remove them, we exploit the variational nature of (\ref{main-equation}).

In fact, (\ref{main-equation}) corresponds ultimately to the Euler-Lagrange equation of a suitable energy functional in the space of parametrizations of toroidal surfaces. One of the Lagrange multipliers can be eliminated by taking variations with respect to the neck-size parameter. This leads to an equation for $a$ with respect to $n$, which is solved in correspondence of some $a_{n}$, as in the statement of the Theorem. The other Lagrange multiplier comes from the rotational invariance and it disappears for free as soon as the first one vanishes.

Concerning the embeddedness property, we need to consider a neighborhood of the reference surface $\Sigma_{n,a}$, characterized as a foliation whose leaves are surfaces parametrized by $X_{n,a}+r x_{a}N_{n,a}$. Here $x_{a}$ is the mapping appearing in (\ref{eq:Xa-vector}) and $r$ is a real parameter with $|r|$ sufficiently small. Observe that for $r=0$ the surface parametrized by $X_{n,a}$ has circular cross sections, with variable diameters; as soon as $r\ne 0$, the surface parametrized by $X_{n,a}+r x_{a}N_{n,a}$ has almost circular cross sections, with variable diameters, bounding planar non-convex, star-shaped domains. Then we have to show that the surface parametrized by $X_{n,a_{n}}+\varphi_{n}N_{n,a_{n}}$ is contained in the above-described neighborhood of $\Sigma_{n,a_{n}}$, and this fact rests on a very sharp estimate on $\varphi_{n}$. In turn, this depends on a good choice of the functional spaces where the problem is set. %They are H\"older spaces with some special weighted norms.

We point out that our study, in few initial parts, shares some features with the work \cite{MazzeoPacard2001} by Mazzeo and Pacard concerning the construction of complete, non compact, CMC surfaces, with a number of ends going to infinity made by semi-unduloids. In fact, this kind of research originates from the already mentioned breakthrough paper \cite{Wente1986} by Wente, and in a different direction, from the works by Kapouleas (see \cite{Kapouleas1990}, \cite{Kapouleas1991}) who constructed a plethora of complete, CMC surfaces, of any genus, obtained by attaching together spheres and pieces of Delaunay surfaces. Concerning this kind of problems, among the others, let us quote also \cite{BreinerKapouleas2014}, \cite{Grosse-Brauckmann1993}, \cite{Grosse-BrauckmannKusnerSullivan2007},  \cite{MazzeoPacardPollack2001}, and the bibliography therein.

An important aspect in our study consists in the use of $n$ as a parameter in order to apply the Lyapunov-Schmidt reduction method. This idea was already implemented by J. Wei and S. Yan in \cite{WY} for a different problem regarding the nonlinear Schr\"odinger equation. We observe that for that problem symmetry assumption was removed in \cite{dPWY}. We suspect that the same can be done also for our problem but this issue seems quite delicate and difficult. 

Finally let us point out that the corresponding (actually simpler) problem about compact, oriented, two-dimensional surfaces of genus zero into Euclidean 3-space has been already studied, and a quite exhaustive literature is nowadays available. See, e.g., \cite{CM2004}, \cite{CaMu11}, \cite{TrWe} and the references therein.

The paper is organized as follows: in Section \ref{S:unduloids} we recall some preliminaries about unduloids.
%Some of the results contained there can be found also in other papers even if with different notation.
In Section \ref{S:toroidal-unduloids} we introduce the toroidal unduloids which provide the supporting surfaces for our construction, we derive the Jacobi operator about them, and we compute its limit as $n\to\infty$. Section \ref{S:linearized-problem} contains the study of the linear problem for the limit Jacobi operator $\mathfrak{L}_{a}$ with some fundamental uniform estimates. In Section \ref{S:reduced} we start to study problem (\ref{main-equation}) and, with an application of the contraction principle, we convert it into a reduced finite-dimensional problem. In Section \ref{S:variational} the variational argument is displayed and we arrive to prove the existence of regular solutions of (\ref{main-equation}). Finally the embeddedness property is discussed in Section \ref{S:embeddedness}.

\begin{comment}
For the sake of convenience, let us collect a list of notation used in the paper:
\begin{itemize}
\item
$X=X(t,\theta)\in\R^{3}$ parametrization
\item
$N=\frac{X_{t}\wedge X_{\theta}}{|X_{t}\wedge X_{\theta}|}$ normal or Gauss map
\item
${\E},{\F},{\G}$ and ${\L},{\M},{\N}$ coefficients of the first and the second fundamental forms, respectively
\item
$\mathfrak{M}(X)$ mean curvature of a surface with parametrization $X$
\item
$\mathfrak{L}$, $\mathfrak{F}$, $\mathfrak{B}$,  $\mathfrak{T}$ operators between H\"older-type function spaces $\mathscr{X}$, $\mathscr{Y}$
\item
$x,z,p,q,w,\varphi$ scalar functions of $t$ or of $(t,\theta)$.
\end{itemize}
\end{comment}

\section{Preliminaries on unduloids}
\label{S:unduloids}
Unduloids are constant mean curvature surfaces of revolution in $\R^{3}$ that are generated by the trace of a focus of an ellipse which rolls without sliding on the axis of revolution. They are members of the family of Delaunay surfaces. We define the neck-size of an unduloid as the distance of the unduloid from its axis of revolution.

Let $\Sigma_{a}$ be the unduloid with mean curvature 1, axis of revolution on the $x_{3}$-axis, and neck-size $a$. Hence necessarily $a\in\big(0,\frac{1}{2}\big]$. The limit case $a=0$ corresponds to a sequence of infinitely many tangent unit spheres placed on the $x_{3}$-axis, whereas for $a=\frac{1}{2}$, the unduloid reduces to a cylinder of radius $\frac{1}{2}$. We also assume that for $a\in\big(0,\frac{1}{2}\big)$ the largest circle obtained as a horizontal section of the unduloid $\Sigma_{a}$ belongs to the plane $\{x_{3}=0\}$. Therefore $\Sigma_{a}\cap\{x_{3}=0\}$ is a circle of radius $1-a$.

The unduloid $\Sigma_{a}$ admits a cylindrical parameterization of the form
$$
\widetilde{X}_{a}(z,\theta)=\left[\begin{array}{c}
\rho_{a}(z)\cos\theta\\
\rho_{a}(z)\sin\theta\\
z\end{array}\right]
$$
with $\rho_{a}(0)=\max_{z}\rho_{a}(z)$. The parameterization $\rho_{a}(z)$ solves the Cauchy problem
\begin{equation}
\label{eq:rho}
\left\{\begin{array}{l}-\rho_{a}''(1+(\rho_{a}')^{2})^{-3/2}+\rho_{a}^{-1}(1+(\rho_{a}')^{2})^{-1/2}=2\\
\rho_{a}(0)=1-a\\
\rho'_{a}(0)=0\end{array}\right.
\end{equation}
where the equation expresses the fact that $\Sigma_{a}$ has mean curvature 1. For $a=\frac{1}{2}$ the solution is the constant $\rho_{a}(0)=\frac{1}{2}$, i.e., the cylinder; for $a=0$ the solution is $\rho_{a}(z)=\sqrt{1-z^{2}}$, i.e., the sphere. For $a\in\big(0,\frac{1}{2}\big)$ problem (\ref{eq:rho}) admits a unique periodic solution $\rho_{a}$. The period of $\rho_{a}$, denoted $2h_{a}$, equals the height of a complete period of the unduloid. It is known that $h_{a}$ is a regular increasing function of $a\in\big(0,\frac{1}{2}\big)$ with
$$
\inf_{a}h_{a}=1\quad\text{and}\quad\sup_{a}h_{a}=\frac{\pi}{2}.
$$
Notice also that the unduloid $\Sigma_{a}$ is symmetric with respect to reflection about the $x_{1}x_{2}$-plane, that is, the parameterization $\rho_{a}$ is an even function.

Cylindrical coordinates are not always adequate to study the problem. Instead, isothermal coordinates are preferable. In fact, there exists a diffeomorphism $z=z_{a}(t)$, uniquely defined with the conditions $z_{a}(0)=0$ and $z_{a}'>0$, such that
\begin{equation}
\label{eq:X}
X_{a}(t,\theta):=\widetilde{X}_{a}(z_{a}(t),\theta)=\left[\begin{array}{c}
x_{a}(t)\cos\theta\\
x_{a}(t)\sin\theta\\
z_{a}(t)\end{array}\right]\quad\text{with}\quad x_{a}:=\rho_{a}\circ z_{a}
\end{equation}
turns out to be a conformal parametrization of the unduloid, i.e.,
$$
(X_{a})_{t}\cdot(X_{a})_{\theta}=0=\left|(X_{a})_{t}\right|^{2}-\left|(X_{a})_{\theta}\right|^{2}.
$$
The conformality of $X_{a}$ is equivalent to the equation
\begin{equation}
\label{eq:conformality}
x_{a}^{2}=(x'_{a})^{2}+(z'_{a})^{2}.
\end{equation}
The property of $\Sigma_{a}$ of having mean curvature 1 is expressed in terms of the conformal parametrization $X_{a}$ by means of the equation
\begin{equation}
\label{eq:CMC-X}
\Delta X_{a}=2(X_{a})_{t}\wedge(X_{a})_{\theta}
\end{equation}
or, equivalently, in terms of the following system of ode's
\begin{equation}
\label{eq:CMC}
\left\{\begin{array}{l}x_{a}''-x_{a}=-2x_{a}z'_{a}\\ z''_{a}=2x_{a}x'_{a}~\!.\end{array}\right.
\end{equation}

\begin{Lemma}
\label{L:c-a}
The mappings $x_{a}(t)$ and $z_{a}(t)$ are the only solutions to the following problems
$$
\left\{\begin{array}{l}
x''_{a}=(1-2\gamma_{a})x_{a}-2x_{a}^{3}\\ x_{a}(0)=1-a~\!,~x'_{a}(0)=0
\end{array}\right.
\text{~~and~~}
\left\{\begin{array}{l}
z'_{a}=\gamma_{a}+x_{a}^{2}\\ z_{a}(0)=0
\end{array}\right.
\text{~~where~~}\gamma_{a}:=a(1-a)~\!.
$$
\end{Lemma}

\Proof
The initial conditions for both problems were already observed. Let us prove
\begin{equation}
\label{eq:z'a}
z'_{a}=\gamma_{a}+x_{a}^{2}~\!.
\end{equation}
One can observe that if $\rho_{a}(z)$ solves (\ref{eq:rho}), then
$$
H(\rho_{a},\rho_{a}'):=\rho_{a}^{2}-\frac{\rho_{a}}{\sqrt{1+(\rho_{a}')^{2}}}
$$
is constant. In particular $H(\rho_{a}(z),\rho_{a}'(z))=H(\rho_{a}(0),\rho_{a}'(0))$, namely
\begin{equation}
\label{eq:conservation}
\rho_{a}^{2}-\frac{\rho_{a}}{\sqrt{1+(\rho_{a}')^{2}}}=-a(1-a)~\!.
\end{equation}
Differentiating the identity $x_{a}=\rho_{a}\circ z_{a}$ and using (\ref{eq:conservation}), one obtains (\ref{eq:z'a}). Finally the differential equation for $x_{a}$ follows from (\ref{eq:CMC}) and (\ref{eq:z'a}).
\qed

The mappings $x_{a}(t)$ and $z_{a}(t)$ are respectively even and odd, and satisfy
$$
x_{a}(t+2\tau_{a})=x_{a}(t)\quad\text{and}\quad z_{a}(t+2\tau_{a})=2h_{a}+z_{a}(t)~\!,
$$
for some $\tau_{a}>0$. The behaviour of $\tau_{a}$ and $h_{a}$ in a right neighborhood of $0$ is displayed by the next lemma.

\begin{Lemma}\label{L:tau-a-expansion}
There exist $C^{1}$ mappings $T,S\colon\big(0,\frac{1}{2}\big]\to\R$ bounded and with bounded derivatives in $\big(0,\frac{\pi}{2}\big]$ such that
$$
\tau_{a}=-2\log a+T(a)\quad\text{and}\quad h_{a}=1-2a\log a+a~\!S(a)~\!.
$$
\end{Lemma}

\Proof
It is known (see, e.g., \cite{Kapouleas1990} or \cite{MazzeoPacard2001}) that $\tau_{a}$ and $h_{a}$ can be expressed in terms of the complete elliptic integrals of the first and second kind
\begin{equation}
\label{eq:elliptic-integrals}
K(k):=\int_{0}^{\pi/2}\frac{dr}{\sqrt{1-k^{2}\sin^{2}r}}~\!,\quad
E(k):=\int_{0}^{\pi/2}\sqrt{1-k^{2}\sin^{2}r}~\!dr~\!,
\end{equation}
according to the formulas
\begin{equation}
\label{eq:taua-ha}
\tau_{a}=2K(k_{a})\quad\text{and}\quad h_{a}=\gamma_{a}\tau_{a}+(1-a)E(k_{a})\quad\text{where}\quad k_{a}^{2}=1-\frac{a^{2}}{(1-a)^{2}}~\!.
\end{equation}
The mappings $k\mapsto K(k)$ and $k\mapsto E(k)$ are of class $C^{\infty}$ and their expansions in a left neighborhood of $k=1$ are known and can be found in Sect.19.12 of \cite{DLMF}. More precisely, introducing the complementary modulus $\kappa=\sqrt{1-k^{2}}$, one can write
$$
K(k)=-\log\kappa+\widetilde{K}(\kappa)\quad\text{and}\quad E(k)=1-\frac{\kappa^{2}}{2}\log\kappa+\kappa^{2}\widetilde{E}(\kappa)
$$
for some regular mappings $\widetilde{K}(\kappa)$ and $\widetilde{E}(\kappa)$ which are bounded, with bounded derivatives, in $(0,1)$.
%\begin{gather*}
%K(k)=-\log\kappa+2\log 2-\frac{\kappa^{2}}{4}\log\kappa+\frac{2\log 2-1}{4}\kappa^{2}+O(\kappa^{4}\log\kappa)\\
%E(k)=1-\frac{\kappa^{2}}{2}\log\kappa+\frac{4\log 2-1}{4}\kappa^{2}+O(\kappa^{4}\log\kappa).
%\end{gather*}
Taking $\kappa=1-k_{a}^{2}=\frac{a}{1-a}$, we obtain
\begin{equation}
\label{eq:KE-expansion}
K(k_{a})=-\log a+\frac{1}{2}T\left(a\right)\quad\text{and}\quad
E(k_{a})=1-\frac{a^{2}}{2}\log a+a^{2}S\left(a\right)
\end{equation}
where $T\left(a\right)$ and $S\left(a\right)$ are regular, bounded mappings, with bounded derivatives, in $\big(0,\frac{1}{2}\big)$.
Then the conclusion follows in view of (\ref{eq:taua-ha}).
\qed

Finally let us discuss the limit behaviour of $x_{a}$ and $z_{a}$, as well as useful uniform bounds.

\begin{Lemma}
\label{L:xa-limit}
The families $\{x_{a}\}_{a\in(0,\frac{1}{2}]}$ and $\{z'_{a}\}_{a\in(0,\frac{1}{2}]}$ are bounded in $C^{3}(\R)$. Moreover $x_{a}(t)\to\sech t$ and $z_{a}(t)\to\tanh t$ in $C^{2}_{loc}(\R)$, as $a\to 0$.
In addition, for $a>0$ small enough
\begin{equation}
\label{eq:xa-bound}
0<x_{a}(t)\le(1-a)\sqrt{\sech t}\quad\forall t\in[-\tau_{a},\tau_{a}]~\!.
\end{equation}
\end{Lemma}

\Proof
We know that $a\le x_{a}\le 1-a$. Moreover, from (\ref{eq:conformality}) and
\begin{equation}
\label{eq:xa-system}
x_{a}''=(1-2\gamma_{a})x_{a}-2x_{a}^{3}
\end{equation}
a uniform estimate in $C^{2}$ holds. Then, differentiating  (\ref{eq:xa-system}), one obtains that the family $\{x_{a}\}_{a\in(0,\frac12]}$ is bounded also in $C^{3}$. In view of Lemma \ref{L:c-a}, the same holds also for the family $\{z'_{a}\}_{a\in(0,\frac12]}$.
Let us prove that $x_{a}(t)\to\sech t$ in $C^{2}_{loc}(\R)$, when $a\to 0$. Because of the uniform bound in $C^{3}$, along a sequence, $x_{a}(t)\to x_{0}(t)$ in $C^{2}_{loc}(\R)$ and, passing to the limit in the Cauchy problem for $x_{a}$ stated in Lemma \ref{L:c-a}, we obtain that $x_{0}(t)$ solves
$$
\left\{\begin{array}{l}
x''=x-2x^{3}\\ x(0)=1\\ x'(0)=0
\end{array}\right.
$$
namely $x_{0}(t)=\sech t$. The uniqueness of the limit function implies that the same conclusion holds true as $a\to 0$ and not just along a sequence. Using again Lemma \ref{L:c-a} and the condition $z_{a}(0)=0$ one infers that $z_{a}(t)\to\tanh t$ as $a\to 0$ in $C^{2}_{loc}(\R)$. In order to check (\ref{eq:xa-bound}), let us introduce the auxiliary function $f_{a}(t)=x_{a}(t)^{2}\cosh t$. By Lemma \ref{L:c-a}
$$
f_{a}''(t)=\left[(5-8\gamma_{a})x_{a}(t)^{2}-6x_{a}(t)^{4}-2\gamma_{a}^{2}\right]\cosh t+4x_{a}(t)x'_{a}(t)\sinh t~\!.
$$
We point out that $f_{a}'(0)=0$ and $f_{a}''(0)=-(1-a)^{2}(1-4a)$. Then $t=0$ is a strict local maximum for $f_{a}$, for $a\in\big(0,\frac{1}{4}\big)$. Moreover, since $\tau_{a}\sim-2\log a$ as $a\to 0$, we have that $f_{a}(\tau_{a})\sim\frac{1}{2}$,  $f_{a}'(\tau_{a})\sim\frac{1}{2}$, and $f_{a}''(\tau_{a})\sim\frac{3}{2}$, as $a\to 0$. In particular $t=\tau_{a}$ is a strict local maximum for $f_{a}$ restricted to $[0,\tau_{a}]$. Let us show that if $t_{a}$ is a critical point for $f_{a}$ in $(0,\tau_{a})$, then it must be a minimum point. Indeed if $f_{a}'(t_{a})=0$ then 
$$
2x_{a}(t_{a})x'_{a}(t_{a})\sinh t_{a}=-x_{a}(t_{a})^{2}\frac{(\sinh t_{a})^{2}}{\cosh t_{a}}>-x_{a}(t_{a})^{2}\cosh t_{a}
$$
and consequently
$$
f_{a}''(t_{a})\ge
\left[3-8\gamma_{a}-6x_{a}(t_{a})^{2}-2\frac{\gamma_{a}^{2}}{x_{a}(t_{a})^{2}}\right]x_{a}(t_{a})^{2}\cosh t_{a}\ge \left[1-8\gamma_{a}-6x_{a}(t_{a})^{2}\right]x_{a}(t_{a})^{2}\cosh t_{a}
$$
because $x_{a}(t_{a})\ge a>\gamma_{a}$. Since $x_{a}(t)\to \sech t$ in $C^{2}_{loc}(\R)$, also $f_{a}(t)\to \sech t$ in $C^{2}_{loc}(\R)$. Therefore $t_{a}\to\infty$ and $x_{a}(t_{a})\to 0$ as $a\to 0$. Hence $f''_{a}(t_{a})>0$ for $a>0$ small enough. Thus we can infer that $f_{a}$ admits no maximum point in $(0,\tau_{a})$. In conclusion, for $a>0$ small enough, $0<f_{a}(t)\le\max\{f_{a}(0),f_{a}(\tau_{a})\}=(1-a)^{2}$ for every $t\in[0,\tau_{a}]$. Since $f_{a}$ is even, (\ref{eq:xa-bound}) follows.
\qed

%\section{Analytical formulation of the problem}
\section{Toroidal unduloids and the Jacobi operator}
\label{S:toroidal-unduloids}

In this section we introduce the toroidal unduloids and we derive the linearization of the mean curvature operator about them.

For every $a\in\big(0,\frac{1}{2}\big]$ and for every $n\in\mathbb{N}$ sufficiently large, we construct a toroidal surface $\Sigma_{n,a}$ by cutting a section of the unduloid $\Sigma_{a}$ made by $n$ periods, hence with length $2nh_{a}$, and folding it along a circle placed on the $x_{2}x_{3}$-plane, centered at the origin, and with circumference $2nh_{a}$.
The corresponding parametrization is given by
\begin{equation*}
%\label{eq:toroidal-unduloid}
X_{n,a}(t,\theta):=\left[\begin{array}{c}x_{a}(t)\cos\theta\\ \left(\frac{nh_{a}}{\pi}+x_{a}(t)\sin\theta\right)\cos\frac{\pi z_{a}(t)}{nh_{a}}\\ \left(\frac{nh_{a}}{\pi}+x_{a}(t)\sin\theta\right)\sin\frac{\pi z_{a}(t)}{nh_{a}}\end{array}\right]~\!.
\end{equation*}
The map ${X}_{n,a}$ turns out to be $2n\tau_{a}$-periodic with respect to $t$ and $2\pi$-periodic with respect to $\theta$. The surface $\Sigma_{n,a}$ is symmetric with respect to rotations of an angle $\frac{2\pi}{n}$ about the $x_{1}$-axis, namely
$$
R_{\frac{2\pi}{n}}\Sigma_{n,a}=\Sigma_{n,a}
$$
where, in general,
\begin{equation}
\label{eq:rotation-matrix}
R_{\sigma}=\left[\begin{array}{ccc}1&0&0\\ 0&\cos\sigma&-\sin\sigma\\ 0&\sin\sigma&\cos\sigma\end{array}\right].
\end{equation}
This discrete symmetry reflects into the following condition on the parametrization
\begin{equation}
\label{eq:discrete-symmetry}
{X}_{n,a}(t+2\tau_{a},\theta)=R_{\frac{2\pi}{n}}{X}_{n,a}(t,\theta).
\end{equation}
Notice that
$$
{X}_{n,a}(t,\theta)-\frac{nh_{a}}{\pi}\mathbf{e}_{2}
%=\left[\begin{array}{c}\rho_{a}(z)\cos\theta\\ \frac{nh_{a}}{\pi}\big(\cos\frac{\pi z}{nh_{a}}-1\big)+\rho_{a}(z)\sin\theta\cos\frac{\pi z}{nh_{a}}\\ \frac{nh_{a}}{\pi}\sin\frac{\pi z}{nh_{a}}+\rho_{a}(z)\sin\theta\sin\frac{\pi z}{nh_{a}}\end{array}\right]
\to
%\left[\begin{array}{c}\rho_{a}(z)\cos\theta\\ \rho_{a}(z)\sin\theta\\ z\end{array}\right]=
{X}_{a}(t,\theta)\quad\text{as}\quad n\to\infty
$$
uniformly on compact sets of $\R\times\S^{1}$. In other words, we recover the unduloid $\Sigma_{a}$ by taking toroidal unduloids with $n$ periods in the limit $n\to\infty$, up to a suitable natural translation.

It is convenient to change $n$ into
$$
\eps:=\frac{\pi}{nh_{a}}\!~.
$$
Then the parameterization $X_{n,a}$ becomes
\begin{equation}
\label{eq:eps-toroidal-unduloid}
X_{\eps,a}(t,\theta)=\left[\begin{array}{c}x_{a}(t)\cos\theta\\ \left(\eps^{-1}+x_{a}(t)\sin\theta\right)\cos(\eps z_{a}(t))\\ \left(\eps^{-1}+x_{a}(t)\sin\theta\right)\sin(\eps z_{a}(t))\end{array}\right].
\end{equation}
Considering $\eps$ as a (small) continuous parameter, the mapping $X_{\eps,a}$ defines a parameterization of a generalized toroidal unduloid $\Sigma_{\eps,a}:=X_{\eps,a}(\R\times\mathbb{S}^{1})$. It holds that
\begin{equation}
\label{eq:derivatives}
X_{\eps,a}=R_{\eps z_{a}}U_{\eps,a}~\!,\quad(X_{\eps,a})_{t}=R_{\eps z_{a}}V_{\eps,a}~\!,\quad(X_{\eps,a})_{\theta}=R_{\eps z_{a}}W_{\eps,a}~\!,
\end{equation}
where, omitting the explicit dependence on $t$,
\begin{equation}
\label{eq:UVW}
U_{\eps,a}:=\left[\begin{array}{c}x_{a}\cos\theta\\
\eps^{-1}+x_{a}\sin\theta\\ 0\end{array}\right],\quad V_{\eps,a}:=\left[\begin{array}{c}x'_{a}\cos\theta\\ x'_{a}\sin\theta\\ z'_{a}\left(1+\eps x_{a}\sin\theta\right)\end{array}\right],\quad W_{\eps,a}:=\left[\begin{array}{c}-x_{a}\sin\theta\\ x_{a}\cos\theta\\
0\end{array}\right].
\end{equation}
Notice that $U_{\eps,a}$, $V_{\eps,a}$ and $W_{\eps,a}$ are $2\tau_{a}$-periodic with respect to $t$. Moreover
\begin{equation}
\label{eq:VW1}
V_{\eps,a}\cdot W_{\eps,a}=0~\!,\quad|V_{\eps,a}|^{2}=(x'_{a})^{2}+(z'_{a})^{2}\left(1+\eps x_{a}\sin\theta\right)^{2}~\!,\quad|W_{\eps,a}|^{2}=x_{a}^{2}.
\end{equation}
Hence $X_{\eps,a}$ satisfies
$$
\left\{\begin{array}{l}
(X_{\eps,a})_{t}\cdot(X_{\eps,a})_{\theta}=0\\
\left|(X_{\eps,a})_{t}\right|^{2}-\left|(X_{\eps,a})_{\theta}\right|^{2}=
(z'_{a})^{2}\left[(1+\eps x_{a}\sin\theta)^{2}-1\right]~\!.\end{array}\right.
$$
Thus $X_{\eps,a}$ is just an orthogonal parameterization of $\Sigma_{\eps,a}$ but not conformal. However
$$
\left|(X_{\eps,a})_{t}\right|^{2}-\left|(X_{\eps,a})_{\theta}\right|^{2}\to 0\quad\text{as}\quad \eps\to 0
$$
uniformly on $\R\times\S^{1}$. That is, for every fixed $a\in\big(0,\frac{1}{2}\big]$, the parameterization $X_{\eps,a}$ becomes conformal in the limit $\eps\to 0$.

Let
$$
{\normal}_{\eps,a}=\frac{(X_{\eps,a})_{t}\wedge(X_{\eps,a})_{\theta}}{\left|(X_{\eps,a})_{t}\wedge(X_{\eps,a})_{\theta}\right|}
$$
be the unit normal for $\Sigma_{\eps,a}$. From (\ref{eq:derivatives}) it follows that
\begin{gather}
\nonumber
(X_{\eps,a})_{t}\wedge(X_{\eps,a})_{\theta}=R_{\eps z_{a}}(V_{\eps,a}\wedge W_{a})\\
\label{eq:orthogonality}
\left|(X_{\eps,a})_{t}\wedge(X_{\eps,a})_{\theta}\right|=\left|(X_{\eps,a})_{t}\right|\left|(X_{\eps,a})_{\theta}\right|=\left|V_{\eps,a}\right|\left|W_{a}\right|\\
\nonumber
{\normal}_{\eps,a}(t+2\tau_{a},\theta)=R_{\frac{2\pi}{n}}{\normal}_{\eps,a}(t,\theta)\text{~~if~~$\eps=\frac{\pi}{nh_{a}}$ with $n\in\mathbb{N}$~\!.}
\end{gather}

%$$
%(X_{n,a})_{t}\wedge(X_{n,a})_{\theta}(t,\theta)=\left[\begin{array}{ccc}1&0&0\\ 0&\cos\frac{\pi z_{a}(t)}{nh_{a}}&-\sin\frac{\pi z_{a}(t)}{nh_{a}}\\ 0&\sin\frac{\pi z_{a}(t)}{nh_{a}}&\cos\frac{\pi z_{a}(t)}{nh_{a}}\end{array}\right]\left[\begin{array}{c}-x_{a}(t)z'_{a}(t)\left(1+\frac{\pi}{nh_{a}}x_{a}(t)\sin\theta\right)\cos\theta\\-x_{a}(t)z'_{a}(t)\left(1+\frac{\pi}{nh_{a}}x_{a}(t)\sin\theta\right)\sin\theta\\ x_{a}(t)x'_{a}(t)\end{array}\right],
%$$

For any regular function $\varphi\colon\R/_{2\tau_{a}}\times\R/_{2\pi}\to\R$ sufficiently small, the mapping $X_{\eps,a}+\varphi{\normal}_{\eps,a}$
provides a regular parameterization of a surface $\Sigma_{\varphi}$ given by a normal graph over $\Sigma_{\eps,a}$.

In general, given a surface $\Sigma$ with parameterization $X(t,\theta)\in\R^{3}$ and normal versor
\begin{equation}
\label{eq:normal-definition}
{\normal}=\frac{X_{t}\wedge X_{\theta}}{\left|X_{t}\wedge X_{\theta}\right|}~\!,
\end{equation}
the mean curvature of $\Sigma$ at a given point $X\in\Sigma$, denoted $\mathfrak{M}(X)$, can be expressed in terms of the coefficients of the first and second fundamental form (we use the Gaussian notation)
\begin{equation}
\label{eq:Gaussian-notation}
\begin{array}{c}
{\E}=|X_{t}|^{2}~\!,\quad {\F}=X_{t}\cdot X_{\theta}~\!,\quad {\G}=|X_{\theta}|^{2}~\!,\\
{\L}=X_{tt}\cdot {\normal}~\!,\quad {\M}=X_{t\theta}\cdot {\normal}~\!,\quad {\N}=X_{\theta\theta}\cdot {\normal}~\!,
\end{array}
\end{equation}
as
\begin{equation}
\label{eq:mean-curvature-def}
\mathfrak{M}=\frac{{\E}{\N}-2{\F}{\M}+{\G}{\L}}{2({\E}{\G}-{\F}^{2})}~\!.
\end{equation}
In particular, for a conformal parameterization, i.e., ${\F}=0={\E}-{\G}$, one has
$$
\mathfrak{M}=\frac{{\N}+{\L}}{2{\E}}~\!.
$$
This is the case of the parameterization $X_{a}$ of the unduloid with neck-size $a$ which is a surface with constant mean curvature 1. Indeed, denoting ${\N}_{a}=(X_{a})_{\theta\theta}\cdot {\normal}_{a}$, ${\L}_{a}=(X_{a})_{tt}\cdot {\normal}_{a}$, ${\E}_{a}=\left|(X_{a})_{t}\right|^{2}$, we have that
$$
\mathfrak{M}(X_{a})=\frac{{\N}_{a}+{\L}_{a}}{2{\E}_{a}}=\frac{-x_{a}''z_{a}'+z_{a}''x_{a}'+x_{a}z_{a}'}{2x_{a}^{3}}=1
$$
thanks to (\ref{eq:conformality})--(\ref{eq:CMC}).

Consider now the surface $\Sigma_{\varphi}$ parameterized by $X_{\eps,a}+\varphi{\normal}_{\eps,a}$, with $\varphi\colon\R/_{2\tau_{a}}\times\R/_{2\pi}\to\R$. This surface $\Sigma_{\varphi}$ is close to the (generalized) toroidal unduloid $\Sigma_{\eps,a}$ for $\varphi$ small. The mean curvature of $\Sigma_{\varphi}$ is given by
$$
\mathfrak{M}(X_{\eps,a}+\varphi{\normal}_{\eps,a})=\frac{{\E}_{\varphi}{\N}_{\varphi}-2{\F}_{\varphi}{\M}_{\varphi}+{\G}_{\varphi}{\L}_{\varphi}}{2({\E}_{\varphi}{\G}_{\varphi}-{\F}_{\varphi}^{2})}
$$
where
\begin{gather*}
{\E}_{\varphi}=|(X_{\eps,a}+\varphi{\normal}_{\eps,a})_{t}|^{2}~\!,\quad {\F}_{\varphi}=(X_{\eps,a}+\varphi{\normal}_{\eps,a})_{t}\cdot (X_{\eps,a}+\varphi{\normal}_{\eps,a})_{\theta}~\!,\quad {\G}_{\varphi}=|(X_{\eps,a}+\varphi{\normal}_{\eps,a})_{\theta}|^{2}~\!,\\
{\L}_{\varphi}=(X_{\eps,a}+\varphi{\normal}_{\eps,a})_{tt}\cdot {\normal}_{\varphi}~\!,\quad {\M}_{\varphi}=(X_{\eps,a}+\varphi{\normal}_{\eps,a})_{t\theta}\cdot {\normal}_{\varphi}~\!,\quad {\N}_{\varphi}=(X_{\eps,a}+\varphi{\normal}_{\eps,a})_{\theta\theta}\cdot {\normal}_{\varphi}~\!,\\
{\normal}_{\varphi}=\frac{(X_{\eps,a}+\varphi{\normal}_{\eps,a})_{t}\wedge (X_{\eps,a}+\varphi{\normal}_{\eps,a})_{\theta}}{\left|(X_{\eps,a}+\varphi{\normal}_{\eps,a})_{t}\wedge (X_{\eps,a}+\varphi{\normal}_{\eps,a})_{\theta}\right|}~\!.
\end{gather*}

The Jacobi operator on the (generalized) toroidal unduloid $X_{\eps,a}$ is defined as
\begin{equation}
\label{eq:L-eps-a}
\mathfrak{L}_{\eps,a}:=2x_{a}^{2}\left.\frac{\partial}{\partial\varphi}\mathfrak{M}(X_{\eps,a}+\varphi{\normal}_{\eps,a})\right|_{\varphi=0}~\!.
\end{equation}
As we will see in the sequel, the factor $2x_{a}^{2}$ is placed for the sake of convenience in computations. We are interested in finding the limiting behaviour of $\mathfrak{L}_{\eps,a}$ as $\eps\to 0$.

\begin{Lemma}
\label{L:L-na}
One has that
\begin{equation}
\label{eq:Leps-a}
\mathfrak{L}_{\eps,a}=b_{\eps,a}(t,\theta)\partial_{tt}+\partial_{\theta\theta}+c_{\eps,a}(t,\theta)+d_{\eps,a}(t,\theta)\partial_{t}+e_{\eps,a}(t,\theta)\partial_{\theta}
\end{equation}
where the coefficients are given by (\ref{eq:coefficients}) and by the related formulas. Moreover, fixing $\alpha\in(0,1)$, there exists a constant $C>0$ independent of $\eps$ and $a$, such that  for every $\eps>0$ small enough and for every $a\in\big(0,\frac{1}{2}\big]$ one has
\begin{equation}
\label{eq:Leps-a-estimate}
\begin{array}{c}
\|b_{\eps,a}-1\|_{C^{0,\alpha}(\R^{2})}\le C\eps~\!,\quad \left\|c_{\eps,a}-2\left(x_{a}^{2}+\frac{\gamma_{a}^{2}}{x_{a}^{2}}\right)\right\|_{C^{0,\alpha}(\R^{2})}\le C\eps~\!,\\
\quad\|d_{\eps,a}\|_{C^{0,\alpha}(\R^{2})}\le C\eps~\!,\quad\|e_{\eps,a}\|_{C^{0,\alpha}(\R^{2})}\le C\eps~\!.
\end{array}
\end{equation}
\end{Lemma}

\Proof
To simplify notation, we set $X=X_{\eps,a}$, ${\normal}={\normal}_{\eps,a}$. A mapping $\mathcal{A}(\varphi,\varphi_{t},\varphi_{\theta},\varphi_{tt},\varphi_{t\theta},\varphi_{\theta\theta})$ will be denoted just by $\mathcal{A}_{\varphi}$. Moreover we will compute first order expansion,  writing
$$
\mathcal{A}_{\varphi}\approx A+\varphi B_{1}+\varphi_{t} B_{2}+\varphi_{\theta} B_{3}+C_{1}\varphi_{tt}+C_{2}\varphi_{t\theta}+C_{3}\varphi_{\theta\theta}
$$
where
$A=\mathcal{A}_{\varphi}\big|_{\varphi=0}$, $B_{1}=\tfrac{\partial\mathcal{A}_{\varphi}}{\partial\varphi}\big|_{\varphi=0}$, $B_{2}=\tfrac{\partial\mathcal{A}_{\varphi}}{\partial\varphi_{t}}\big|_{\varphi=0}$, $B_{3}=\tfrac{\partial\mathcal{A}_{\varphi}}{\partial\varphi_{\theta}}\big|_{\varphi=0}$, 
$C_{1}=\tfrac{\partial\mathcal{A}_{\varphi}}{\partial\varphi_{tt}}\big|_{\varphi=0}$, $C_{2}=\tfrac{\partial\mathcal{A}_{\varphi}}{\partial\varphi_{t\theta}}\big|_{\varphi=0}$, $C_{3}=\tfrac{\partial\mathcal{A}_{\varphi}}{\partial\varphi_{\theta\theta}}\big|_{\varphi=0}$. 
Following this agreement, we have that
\begin{equation}
\label{eq:EGF}
{\E}_{\varphi}\approx|X_{t}|^{2}-2\left({X}_{tt}\cdot{\normal}\right)\varphi~\!,\quad
{\G}_{\varphi}\approx|X_{\theta}|^{2}-2\left({X}_{\theta\theta}\cdot{\normal}\right)\varphi~\!,\quad
{\F}_{\varphi}\approx-2\left({X}_{t\theta}\cdot{\normal}\right)\varphi~\!.
\end{equation}
Then
\begin{gather*}
{\E}_{\varphi}{\G}_{\varphi}-{\F}_{\varphi}^{2}\approx|X_{t}|^{2}|X_{\theta}|^{2}-2\left[\left({X}_{tt}\cdot{\normal}\right)|X_{\theta}|^{2}+\left({X}_{\theta\theta}\cdot{\normal}\right)|X_{t}|^{2}\right]\varphi\\
\frac{1}{{\E}_{\varphi}{\G}_{\varphi}-{\F}_{\varphi}^{2}}\approx\frac{1}{|X_{t}|^{2}|X_{\theta}|^{2}}\left[1+2\left(\frac{{X}_{tt}\cdot{\normal}}{|X_{t}|^{2}}+\frac{{X}_{\theta\theta}\cdot{\normal}}{|X_{\theta}|^{2}}\right)\varphi\right]~\!.
\end{gather*}
In order to write the expansion of ${\normal}_{\varphi}$ we start by computing
\begin{gather}
\nonumber
({X}+\varphi{\normal})_{t}\wedge  ({X}+\varphi{\normal})_{\theta}\approx X_{t}\wedge  X_{\theta}+\varphi\left(X_{t}\wedge {\normal}_{\theta}+{\normal}_{t}\wedge  X_{\theta}\right)+\varphi_{\theta}X_{t}\wedge {\normal}+\varphi_{t}{\normal}\wedge X_{\theta}
\\
\label{eq:useful}
X_{t}\wedge {\normal}=-\left|X_{t}\right|\widehat{X}_{\theta}~\!,\quad\quad\qquad\qquad
{\normal}\wedge X_{\theta}=-\left|X_{\theta}\right|\widehat{X}_{t}
\\
\label{316B}
X_{t}\wedge {\normal}_{\theta}=
\frac{\left|X_{t}\right|}{\left|X_{\theta}\right|}
\left[
\left(\widehat{X}_{t}\cdot{X}_{\theta\theta}\right)\widehat{X}_{t}+
\left(\widehat{X}_{\theta}\cdot{X}_{\theta\theta}\right)\widehat{X}_{\theta}-{X}_{\theta\theta}
\right]
\\
\label{316C}
{\normal}_{t}\wedge X_{\theta}=
\frac{\left|X_{\theta}\right|}{\left|X_{t}\right|}
\left[
\left(\widehat{X}_{\theta}\cdot{X}_{tt}\right)\widehat{X}_{\theta}+
\left(\widehat{X}_{t}\cdot{X}_{tt}\right)\widehat{X}_{t}-{X}_{tt}
\right]
\end{gather}
where
\begin{equation}
\label{eq:useful2}
\widehat{X}_{t}=\frac{X_{t}}{|X_{t}|}~\!,\quad\quad \widehat{X}_{\theta}=\frac{X_{\theta}}{|X_{\theta}|}~\!.
\end{equation}
Then
\begin{gather*}
\left|({X}+\varphi{\normal})_{t}\wedge  ({X}+\varphi{\normal})_{\theta}\right|\approx |{X}_{t}\wedge{X}_{\theta}|\left[1-\varphi\left(
\frac{X_{\theta\theta}}{|X_{\theta}|^{2}}+\frac{X_{tt}}{|X_{t}|^{2}}
\right)\cdot{\normal}\right]
\\
\frac{1}{\left|({X}+\varphi{\normal})_{t}\wedge({X}+\varphi{\normal})_{\theta}\right|}
\approx \frac{1}{|{X}_{t}\wedge{X}_{\theta}|}\left[1+\varphi\left(
\frac{X_{\theta\theta}}{|X_{\theta}|^{2}}+\frac{X_{tt}}{|X_{t}|^{2}}
\right)\cdot{\normal}\right]
\end{gather*}
and finally
$$
{\normal}_{\varphi}\approx{\normal}+
\left[\frac{X_{t}\wedge {\normal}_{\theta}+{\normal}_{t}\wedge  X_{\theta}}{|{X}_{t}\wedge{X}_{\theta}|}+{\normal}\left(
\frac{X_{\theta\theta}}{|X_{\theta}|^{2}}+\frac{X_{tt}}{|X_{t}|^{2}}
\right)\cdot{\normal}\right]\varphi
-\frac{X_{t}}{|X_{t}|^{2}}\varphi_{t}
-\frac{X_{\theta}}{|X_{\theta}|^{2}}\varphi_{\theta}~\!.
$$
Thus we have an expansion of the form
$$
{\normal}_{\varphi}\approx{\normal}+A\varphi+B\varphi_{t}+C\varphi_{\theta}
$$
with
$$
A=\frac{X_{t}\wedge {\normal}_{\theta}+{\normal}_{t}\wedge  X_{\theta}}{|{X}_{t}\wedge{X}_{\theta}|}+{\normal}\left(
\frac{X_{\theta\theta}}{|X_{\theta}|^{2}}+\frac{X_{tt}}{|X_{t}|^{2}}
\right)\cdot{\normal}~\!,\quad
B=-\frac{X_{t}}{|X_{t}|^{2}}
~\!,\quad
C=-\frac{X_{\theta}}{|X_{\theta}|^{2}}~\!.
$$
Then the following expansions hold:
\begin{equation*}
\begin{split}
{\L}_{\varphi}
%&=(X+\varphi{\normal})_{tt}\cdot {\normal}_{\varphi}\approx(X_{tt}+\varphi_{tt}{\normal}+2\varphi_{t}{\normal}_{t}+\varphi{\normal}_{tt})\cdot({\normal}+A\varphi+B\varphi_{t}+C\varphi_{\theta})\\
&\approx X_{tt}\cdot{\normal}
+\left({\normal}_{tt}\cdot{\normal}+X_{tt}\cdot A\right)\varphi
+\left(2{\normal}_{t}\cdot{\normal}+X_{tt}\cdot B\right)\varphi_{t}
+\left(X_{tt}\cdot C\right)\varphi_{\theta}
+\varphi_{tt}
\\
%\end{split}\end{equation*}
%
%\begin{equation*}\begin{split}
{\M}_{\varphi}
%&=(X+\varphi{\normal})_{t\theta}\cdot {\normal}_{\varphi}\approx(X_{t\theta}+\varphi_{t\theta}{\normal}+\varphi_{t}{\normal}_{\theta}+\varphi_{\theta}{\normal}_{t}+\varphi{\normal}_{t\theta})\cdot({\normal}+A\varphi+B\varphi_{t}+C\varphi_{\theta})\\
&\approx X_{t\theta}\cdot{\normal}
+\left({\normal}_{t\theta}\cdot{\normal}+X_{t\theta}\cdot A\right)\varphi
+\left({\normal}_{\theta}\cdot{\normal}+X_{t\theta}\cdot B\right)\varphi_{t}
+\left({\normal}_{t}\cdot{\normal}+X_{t\theta}\cdot C\right)\varphi_{\theta}
+\varphi_{t\theta}
\\
%\end{split}\end{equation*}
%
%\begin{equation*}\begin{split}
{\N}_{\varphi}
%&=(X+\varphi{\normal})_{\theta\theta}\cdot {\normal}_{\varphi}\approx(X_{\theta\theta}+\varphi_{\theta\theta}{\normal}+2\varphi_{\theta}{\normal}_{\theta}+\varphi{\normal}_{\theta\theta})\cdot({\normal}+A\varphi+B\varphi_{t}+C\varphi_{\theta})\\
&\approx X_{\theta\theta}\cdot{\normal}
+\left({\normal}_{\theta\theta}\cdot{\normal}+X_{\theta\theta}\cdot A\right)\varphi
+\left(X_{\theta\theta}\cdot B\right)\varphi_{t}
+\left(2{\normal}_{\theta}\cdot{\normal}+X_{\theta\theta}\cdot C\right)\varphi_{\theta}
+\varphi_{\theta\theta}
\end{split}
\end{equation*}
which, together with (\ref{eq:EGF}), yield
\begin{equation*}
\begin{split}
{\E}_{\varphi}{\N}_{\varphi}&\approx
%\left(|X_{t}|^{2}+2\left(X_{tt}\cdot{\normal}\right)\varphi\right)
%\left(X_{\theta\theta}\cdot{\normal}
%+\left[{\normal}_{\theta\theta}\cdot{\normal}+X_{\theta\theta}\cdot A\right]\varphi
%+\left[X_{\theta\theta}\cdot B\right]\varphi_{t}
%+\left[2{\normal}_{\theta}\cdot{\normal}+X_{\theta\theta}\cdot C\right]\varphi_{\theta}\right)+\varphi_{\theta\theta}\\
%&\approx
(X_{\theta\theta}\cdot{\normal})|X_{t}|^{2}
+\left[\left({\normal}_{\theta\theta}\cdot{\normal}+X_{\theta\theta}\cdot A\right)|X_{t}|^{2}-2\left(X_{tt}\cdot{\normal}\right)\left(X_{\theta\theta}\cdot{\normal}\right)\right]\varphi\\
&\qquad\qquad\qquad\qquad\qquad
+\left(X_{\theta\theta}\cdot B\right)|X_{t}|^{2}\varphi_{t}
+\left(2{\normal}_{\theta}\cdot{\normal}+X_{\theta\theta}\cdot C\right)|X_{t}|^{2}\varphi_{\theta}+|X_{t}|^{2}\varphi_{\theta\theta}
\\
%\end{split}\end{equation*}
%
%\begin{equation*}\begin{split}
{\G}_{\varphi}{\L}_{\varphi}&\approx
%\left(|X_{\theta}|^{2}+2\left(X_{\theta\theta}\cdot{\normal}\right)\varphi\right)\left(X_{tt}\cdot{\normal}+\left[{\normal}_{tt}\cdot{\normal}+X_{tt}\cdot A\right]\varphi+\left[2{\normal}_{t}\cdot{\normal}+X_{tt}\cdot B\right]\varphi_{t}+\left[X_{tt}\cdot C\right]\varphi_{\theta}+\varphi_{tt}\right)\\
%&\approx
\left(X_{tt}\cdot{\normal}\right)|X_{\theta}|^{2}
+\left[\left({\normal}_{tt}\cdot{\normal}+X_{tt}\cdot A\right)|X_{\theta}|^{2}-2\left(X_{\theta\theta}\cdot{\normal}\right)\left(X_{tt}\cdot{\normal}\right)\right]\varphi\\
&\qquad\qquad\qquad\qquad\qquad
+\left(2{\normal}_{t}\cdot{\normal}+X_{tt}\cdot B\right)|X_{\theta}|^{2}\varphi_{t}
+\left(X_{tt}\cdot C\right)|X_{\theta}|^{2}\varphi_{\theta}
+|X_{\theta}|^{2}\varphi_{tt}
\\
%\end{split}\end{equation*}
%
%\begin{equation*}\begin{split}
{\F}_{\varphi}{\M}_{\varphi}&\approx
%-\left[\left(X_{tt}+X_{\theta\theta}\right)\cdot {\normal}\right]\varphi\left(X_{t\theta}\cdot{\normal}+\left[{\normal}_{t\theta}\cdot{\normal}+X_{t\theta}\cdot A\right]\varphi+\left[{\normal}_{\theta}\cdot{\normal}+X_{t\theta}\cdot B\right]\varphi_{t}+\left[{\normal}_{t}\cdot{\normal}+X_{t\theta}\cdot C\right]\varphi_{\theta}+\varphi_{t\theta}\right)\\
%&\approx
-2\left(X_{t\theta}\cdot{\normal}\right)^{2}\varphi
\end{split}
\end{equation*}
and hence
\begin{equation*}\begin{split}
\frac{{\E}_{\varphi}{\N}_{\varphi}-2{\F}_{\varphi}{\M}_{\varphi}+{\G}_{\varphi}{\L}_{\varphi}}{{\E}_{\varphi}{\G}_{\varphi}-{\F}_{\varphi}^{2}}&\approx
\frac{X_{\theta\theta}\cdot{\normal}}{|X_{\theta}|^{2}}+\frac{X_{tt}\cdot{\normal}}{|X_{t}|^{2}}
%\\
%&\quad
+\bigg[
\frac{{\normal}_{\theta\theta}\cdot{\normal}+X_{\theta\theta}\cdot A}{|X_{\theta}|^{2}}
+\frac{{\normal}_{tt}\cdot{\normal}+X_{tt}\cdot A}{|X_{t}|^{2}}\\
&%\quad
~-\frac{4({X}_{tt}\cdot{\normal})({X}_{\theta\theta}\cdot{\normal})}{|X_{t}|^{2}|X_{\theta}|^{2}}
+\frac{4({X}_{t\theta}\cdot{\normal})^{2}}{|X_{t}|^{2}|X_{\theta}|^{2}}
+2\left(\frac{X_{\theta\theta}\cdot{\normal}}{|X_{\theta}|^{2}}+\frac{X_{tt}\cdot{\normal}}{|X_{t}|^{2}}\right)^{2}
\bigg]
\varphi
\\
&%\quad
~+\left[\frac{X_{\theta\theta}\cdot B}{|X_{\theta}|^{2}}+\frac{X_{tt}\cdot B}{|X_{t}|^{2}}+\frac{2{\normal}_{t}\cdot{\normal}}{|X_{t}|^{2}}
\right]
\varphi_{t}
%\\
%&\quad
+\left[\frac{X_{tt}\cdot C}{|X_{t}|^{2}}+\frac{X_{\theta\theta}\cdot C}{|X_{\theta}|^{2}}+\frac{2{\normal}_{\theta}\cdot{\normal}}{|X_{\theta}|^{2}}
\right]
\varphi_{\theta}
\\
&%\quad
~+\frac{\varphi_{tt}}{|X_{t}|^{2}}+\frac{\varphi_{\theta\theta}}{|X_{\theta}|^{2}}~\!.
\end{split}\end{equation*}
Since $|{\normal}|=1$ we have that ${\normal}_{t}\cdot{\normal}={\normal}_{\theta}\cdot{\normal}=0$, ${\normal}_{tt}\cdot{\normal}=-|{\normal}_{t}|^{2}$, ${\normal}_{\theta\theta}\cdot{\normal}=-|{\normal}_{\theta}|^{2}$. Setting
\begin{equation}
\label{eq:X''}
X'':=\frac{X_{\theta\theta}}{|X_{\theta}|^{2}}+\frac{X_{tt}}{|X_{t}|^{2}}
\end{equation}
we can write
\begin{equation*}\begin{split}
\frac{{\E}_{\varphi}{\N}_{\varphi}-2{\F}_{\varphi}{\M}_{\varphi}+{\G}_{\varphi}{\L}_{\varphi}}{{\E}_{\varphi}{\G}_{\varphi}-{\F}_{\varphi}^{2}}&\approx X''\cdot{\normal}+\bigg[X''\cdot A
-\frac{|{\normal}_{\theta}|^{2}}{|{X}_{\theta}|^{2}}
-\frac{|{\normal}_{t}|^{2}}{|{X}_{t}|^{2}}
+\frac{2({X}_{tt}\cdot{\normal})^{2}}{|X_{t}|^{4}}
+\frac{2({X}_{\theta\theta}\cdot{\normal})^{2}}{|X_{\theta}|^{4}}\\
&\phantom{\approx X''\cdot{\normal}+}+\frac{4({X}_{t\theta}\cdot{\normal})^{2}}{|X_{t}|^{2}|X_{\theta}|^{2}}\bigg]\varphi+\left(X''\cdot B\right)\varphi_{t}+
\left(X''\cdot C\right)\varphi_{\theta}+\frac{\varphi_{tt}}{|X_{t}|^{2}}+\frac{\varphi_{\theta\theta}}{|X_{\theta}|^{2}}~\!.
\end{split}\end{equation*}
Hence
\begin{equation}
\label{eq:M-na}
\mathfrak{M}(X_{\eps,a})=\frac{1}{2}X''\cdot{\normal}
%=\frac{1}{2}\left[\frac{(X_{\eps,a})_{tt}}{|(X_{\eps,a})_{t}|^{2}}+\frac{(X_{\eps,a})_{\theta\theta}}{|(X_{\eps,a})_{\theta}|^{2}}\right]\cdot{\normal}_{\eps,a}-1-K(X_{\eps,a})
\end{equation}
and
$$
\frac{\partial}{\partial\varphi}\mathfrak{M}(X_{\eps,a}+\varphi{\normal}_{\eps,a})\bigg|_{\varphi=0}=\frac{1}{2}\frac{\partial_{tt}}{|X_{t}|^{2}}+\frac{1}{2}\frac{\partial_{\theta\theta}}{|X_{\theta}|^{2}}+\frac{1}{2}
\left(X''\cdot B\right)\partial_{t}+\frac{1}{2}
\left(X''\cdot C\right)\partial_{\theta}+q
$$
where
$$
q=\frac{1}{2}\left[X''\cdot A
-\frac{|{\normal}_{\theta}|^{2}}{|{X}_{\theta}|^{2}}
-\frac{|{\normal}_{t}|^{2}}{|{X}_{t}|^{2}}
+\frac{2({X}_{tt}\cdot{\normal})^{2}}{|X_{t}|^{4}}
+\frac{2({X}_{\theta\theta}\cdot{\normal})^{2}}{|X_{\theta}|^{4}}
+\frac{4({X}_{t\theta}\cdot{\normal})^{2}}{|X_{t}|^{2}|X_{\theta}|^{2}}\right]~\!.
$$
Thus we obtained (\ref{eq:Leps-a}) with
\begin{equation}
\label{eq:coefficients}
b_{\eps,a}=\frac{x_{a}^{2}}{|X_{t}|^{2}}~\!,\quad c_{\eps,a}=2x_{a}^{2}q~\!,\quad d_{\eps,a}=x_{a}^{2}(X''\cdot B)~\!,\quad e_{\eps,a}=x_{a}^{2}(X''\cdot C)~\!.
\end{equation}
Now we compute the point-wise limit of the coefficients $b_{\eps,a}$,..., $e_{\eps,a}$ as $\eps\to 0$, for fixed $a\in\big(0,\frac{1}{2}\big]$. Firstly, one has that $\mathfrak{M}(X_{\eps,a})\to 1$ as $\eps\to 0$, that means that the generalized toroidal unduloid $\Sigma_{\eps,a}$ has mean curvature close to $1$ as $|\eps|\ll 1$. Moreover, as $\eps\to 0$
$$
|X_{\theta}|=x_{a}~\!,\quad |X_{t}|\to\sqrt{(x'_{a})^{2}+(z'_{a})^{2}}=x_{a}~\!,\quad |X_{t}\wedge X_{\theta}|\to x_{a}^{2}
$$
$$
X''\sim 2{\normal}~\!,\quad X''\cdot{\normal}\to 2~\!,\quad X''\cdot A\to 0,\quad X''\cdot B\to 0~\!,\quad X''\cdot C\to 0~\!.
$$
In particular $b_{\eps,a}\to 1$, $d_{\eps,a}\to 0$ and $e_{\eps,a}\to 0$ as $\eps\to 0$, point-wise, for every fixed $a$. Writing
$$
X_{t}=RV~\!,\quad X_{\theta}=RW\quad\text{where}\quad
R=R_{\eps z_{a}}~\!,\quad V=V_{\eps,a}~\!,\quad W=W_{\eps,a}
$$
and using
\begin{gather*}
RV_{1}\cdot RV_{2}=V_{1}\cdot V_{2}~\!,\quad RV_{1}\wedge RV_{2}=R(V_{1}\wedge V_{2})\quad(V_{1},V_{2}\in\R^{3})~\!,\\
[(\partial_{t}R)V]\cdot RV=0~\!,\quad [(\partial_{t}R)V]\wedge RV=\eps z'_{a}R\left[(\mathbf{e}_{1}\wedge V)\wedge V\right],
\end{gather*}
one obtains that
\begin{equation}
\label{eq:XzzN}
{X}_{\theta\theta}\cdot{\normal}=\frac{X_{t}\cdot(X_{\theta}\wedge X_{\theta\theta})}{|X_{t}\wedge X_{\theta}|}=\frac{W_{\theta}\cdot(V\wedge W)}{|V\wedge W|}\to z'_{a}
\end{equation}
\begin{equation}
\label{eq:XtzN}
{X}_{t\theta}\cdot{\normal}=\frac{X_{t}\cdot(X_{\theta}\wedge X_{\theta\theta})}{|X_{t}\wedge X_{\theta}|}=\frac{V_{\theta}\cdot(V\wedge W)}{|V\wedge W|}\to 0
\end{equation}
\begin{equation}
\label{eq:XttN}
{X}_{tt}\cdot{\normal}=\frac{X_{\theta}\cdot(X_{tt}\wedge X_{t})}{|X_{t}\wedge X_{\theta}|}=\frac{V_{t}\cdot(V\wedge W)}{|V\wedge W|}
+\eps z'_{a}\frac{W\cdot\left[(\mathbf{e}_{1}\wedge V)\wedge V\right]}{|V\wedge W|}\to\frac{x'_{a}z''_{a}-x''_{a}z'_{a}}{x_{a}}=2x_{a}^{2}-z'_{a}
\end{equation}
thanks to (\ref{eq:conformality}) and (\ref{eq:CMC}).
Moreover
\begin{equation*}\begin{split}
|{\normal}_{\theta}|^{2}&=\frac{|(V\wedge W)_{\theta}|^{2}}{|V\wedge W|^{2}}-\frac{\left[(V\wedge W)\cdot(V\wedge W)_{\theta}\right]^{2}}{|V\wedge W|^{4}}\to\frac{(z'_{a})^{2}}{x_{a}^{2}}
%\end{split}\end{equation*}\end{document}
\\
|{\normal}_{t}|^{2}&=
\frac{|{(\partial_{t}R)(V\wedge W)|^{2}}}{|V\wedge W|^{2}}+
\frac{|(V\wedge W)_{t}|^{2}}{|V\wedge W|^{2}}-
\frac{\left[(V\wedge W)\cdot(V\wedge W)_{t}\right]^{2}}{|V\wedge W|^{4}}+
\frac{2\left[(\partial_{t}R)(V\wedge W)\right]\cdot \left[R(V\wedge W)_{t}\right]}{|V\wedge W|^{2}}\\
&\to\frac{(z_{a}'')^{2}+(x_{a}'')^{2}+3(x_{a}')^{2}}{x_{a}^{2}}-\frac{4(x'_{a})^{2}}{x_{a}^{2}}=\frac{(z_{a}'')^{2}+(x_{a}'')^{2}-(x_{a}')^{2}}{x_{a}^{2}}
\end{split}\end{equation*}
and then, again by (\ref{eq:conformality}) and (\ref{eq:CMC}), and using also Lemma \ref{L:c-a},
$$
q\to\frac{1}{2}\left[-\frac{(z'_{a})^{2}}{x_{a}^{4}}-\frac{(z_{a}'')^{2}-(x_{a}')^{2}+(x_{a}'')^{2}}{x_{a}^{4}}+\frac{2(2x_{a}^{2}-z'_{a})^{2}}{x_{a}^{4}}+\frac{2(z'_{a})^{2}}{x_{a}^{4}}
\right]=1+\frac{\gamma_{a}^{2}}{x_{a}^{4}}~\!,
$$
that is $c_{\eps,a}\to x_{a}+\frac{\gamma_{a}^{2}}{x_{a}^{2}}$ as $\eps\to 0$, point-wise. Finally, let us check (\ref{eq:Leps-a-estimate}). Set
\begin{equation}
\label{eq:useful-notation}
\phi_{a}=x_{a}\sin\theta\quad\text{and}\quad w_{a}=\frac{z'_{a}}{x_{a}}~\!,
\end{equation}
and notice that
\begin{equation}
\label{eq:phix}
|\phi_{a}|\le x_{a}<1\quad\text{and}\quad 0<w_{a}\le 1
\end{equation}
(see (\ref{eq:conformality})). One can write
\begin{equation}
\label{eq:X-derivatives}
|X_{\theta}|^{2}=x_{a}^{2}~\!,\quad|X_{t}|^{2}=x_{a}^{2}(1+2\eps\phi_{a}w_{a}^{2}+\eps^{2}\phi_{a}^{2}w_{a}^{2})
\end{equation}
and then
$$
1-b_{\eps,a}=\frac{|X_{t}|^{2}-|X_{\theta}|^{2}}{|X_{t}|^{2}}=\frac{2\eps\phi_{a}w_{a}^{2}+\eps^{2}\phi_{a}^{2}w_{a}^{2}}{1+2\eps\phi_{a}w_{a}^{2}+\eps^{2}\phi_{a}^{2}w_{a}^{2}}~\!.
$$
Thanks to (\ref{eq:phix}) one obtains that
\begin{equation}
\label{eq:b1}
\|1-b_{\eps,a}\|_{C^{0}(\R^{2})}\le C\eps
\end{equation}
with $C>0$ independent of $\eps$ and $a$. Moreover
\begin{equation}
\label{eq:b2}
|(b_{\eps,a})_{\theta}|=\frac{2\eps x_{a}w_{a}^{2}\left(1+\eps\phi_{a}\right)\left|\cos\theta\right|}{(1+2\eps\phi_{a}w_{a}^{2}+\eps^{2}\phi_{a}^{2}w_{a}^{2})^{2}}\le C\eps
\end{equation}
in view of (\ref{eq:phix}). In addition, using also (\ref{eq:CMC}) and Lemma \ref{L:c-a}, we obtain
$$
|(b_{\eps,a})_{t}|=\frac{2\eps x'_{a}z'_{a}(\gamma_{a}-3x_{a}^{2}-2\eps x_{a}^{3})}{x_{a}^{2}(1+2\eps\phi_{a}w_{a}^{2}+\eps^{2}\phi_{a}^{2}w_{a}^{2})^{2}}
$$
and then
\begin{equation}
\label{eq:b3}
|(b_{\eps,a})_{t}|\le C\eps
\end{equation}
because of (\ref{eq:phix}). Hence (\ref{eq:b1})--(\ref{eq:b3}) imply the first estimate in (\ref{eq:Leps-a-estimate}). The other estimates can be proved in a similar way.
\qed

Let us conclude this section with some more estimates which will be useful in the following. As in Lemma \ref{L:L-na}, $\alpha$ is a fixed number in $(0,1)$.

\begin{Lemma}\label{L:M(X)-eps}
There exists a constant $C_{0}>0$ and regular mappings $\xi_{\eps,a}\colon\R/_{2\tau_{a}}\times\R/_{2\pi}\to\R$ such that for every $\eps>0$ small enough and for every $a\in\big(0,\tfrac{1}{2}\big]$ one has
$$
\mathfrak{M}(X_{\eps,a})=1+\eps\xi_{\eps,a}\quad\text{and}\quad \|\xi_{\eps,a}\|_{C^{0,\alpha}(\R^{2})}\le C_{0}~\!.
$$
Moreover, the mappings $\xi_{\eps,a}$ depend on the parameter $a$ in a continuous way.
\end{Lemma}

\Proof
Using (\ref{eq:X''}), (\ref{eq:M-na}), (\ref{eq:XzzN}) and (\ref{eq:XttN}), after computations, one finds
\begin{equation}
\label{eq:M-1-expansion}
\mathfrak{M}(X_{\eps,a})-1=\frac{1-\left(1+\eps\psi_{\eps,a}\right)^{3/2}}{\left(1+\eps\psi_{\eps,a}\right)^{3/2}}+\frac{\eps\left(\sin\theta\right)(5+w_{a}+w_{a}^{2}-w_{a}^{3}+2\eps\phi_{a}+2\eps\phi_{a}w_{a}^{2}+2\eps\phi_{a}w_{a}^{3}+\eps^{2}\phi_{a}^{2}w_{a}^{2})}{2\left(1+\eps\psi_{\eps,a}\right)^{3/2}}
\end{equation}
where $\phi_{a}$ and $w_{a}$ are as in (\ref{eq:useful-notation}) and $\psi_{\eps,a}=2\phi_{a}w_{a}^{2}+\eps\phi_{a}^{2}w_{a}^{2}$. By (\ref{eq:M-1-expansion}), one can derive the desired estimate with the aid of (\ref{eq:conformality}) and (\ref{eq:phix}).
\qed

\begin{Lemma}\label{L:Xt}
The mappings $(X_{\eps,a})_{t}$, $(X_{\eps,a})_{\theta}$, $(X_{\eps,a}\cdot N_{\eps,a})$, $(X_{\eps,a}\cdot N_{\eps,a})_{t}$, and $(X_{\eps,a}\cdot N_{\eps,a})_{\theta}$ are bounded in $C^{0,\alpha}(\R^{2})$ uniformly in $a\in\big(0,\tfrac{1}{2}\big]$ and $\eps>0$ small enough.
\end{Lemma}

\Proof
As far as concerns $(X_{\eps,a})_{t}$ and $(X_{\eps,a})_{\theta}$, one can control their norms in $C^{1}(\R^{2})$, hence in $C^{0,\alpha}(\R^{2})$, by means of (\ref{eq:X-derivatives}), (\ref{eq:phix}) and the bound $|w_{a}'|\le w_{a}$, which is a consequence of the second equation in (\ref{eq:wa}).
Moreover, with the notation already used in the previous lemma, one can write
$$
X_{\eps,a}\cdot N_{\eps,a}=-\frac{w_{a}(x_{a}+\eps\sin\theta)(1+\eps\phi_{a})}{\left(1+\eps\psi_{\eps,a}\right)^{1/2}}~\!.
$$
Using this expression, (\ref{eq:phix}) and the bound $|w_{a}'|\le w_{a}$, one can estimate the derivatives of $X_{\eps,a}\cdot N_{\eps,a}$ in $C^{1}(\R^{2})$ uniformly with respect to the parameters.
\qed

\section{Study of the linearized problem}
\label{S:linearized-problem}

In this section we aim to study some properties of the linear differential operator
\begin{equation}
\label{eq:pa-def}
{\mathfrak{L}}_{a}:=\Delta+2p_{a}\quad\text{where}\quad\Delta=\partial_{tt}+\partial_{\theta\theta}~\!,\quad p_{a}:=x_{a}^{2}+\frac{\gamma_{a}^{2}}{x_{a}^{2}}~\!,
\end{equation}
$x_{a}$ and $\gamma_{a}$ are given by Lemma \ref{L:c-a}, and $a$ is a parameter in $\big(0,\frac{1}{2}\big)$ corresponding to the neck-size of a given unduloid. In fact, the operator ${\mathfrak{L}}_{a}$ arises by linearizing the mean curvature operator about unduloids: according to Lemma \ref{L:L-na}, we have that
$$
\frac{\partial}{\partial\varphi}\mathfrak{M}(X_{\eps,a}+\varphi{\normal}_{\eps,a})\Big|_{\varphi=0}\to\frac{1}{2x_{a}^{2}}{\mathfrak{L}}_{a}\quad\text{as $\eps\to 0$.}
$$

For our purposes, a nice choice of the domain and the target space for the operator ${\mathfrak{L}}_{a}$ is given by the H\"older spaces
$$
\begin{array}{l}
\mathscr{X}_{a}:=\{\varphi\in C^{2,\alpha}(\R/_{2\tau_{a}}\times\R/_{2\pi})~|~\varphi(t,\cdot)=\varphi(-t,\cdot),~\varphi(\cdot,\pi/2-\theta)=\varphi(\cdot,\pi/2+\theta)\}~\!,\\
\mathscr{Y}_{a}:=\{\varphi\in C^{0,\alpha}(\R/_{2\tau_{a}}\times\R/_{2\pi})~|~\varphi(t,\cdot)=\varphi(-t,\cdot),~\varphi(\cdot,\pi/2-\theta)=\varphi(\cdot,\pi/2+\theta)\}~\!,
\end{array}
$$
respectively.

Our first goal is to find the kernel of ${\mathfrak{L}}_{a}$ in $\mathscr{X}_{a}$. The next two lemmata provide some useful information in this direction.

\begin{Lemma}\label{L:N-equation}
Each component of the Gauss map ${\normal}_{a}$ of the unduloid, defined by
$$
{\normal}_{a}:=\dfrac{(X_{a})_{t}\wedge(X_{a})_{\theta}}{|(X_{a})_{t}\wedge(X_{a})_{\theta}|}
$$ satisfies
\begin{equation}
\label{eq:L=0}
\Delta\varphi+2p_{a}\varphi=0\text{~~in~~}\R/_{2\tau_{a}}\times\R/_{2\pi}~\!.
\end{equation}
\end{Lemma}

\Proof
It is known (see \cite{DHS} Theorem 1, Sect. 5.1, pag. 369) that the Gauss map ${\normal}$ of a parametric regular surface with mean curvature $H$ and Gaussian curvature $K$ satisfies the differential equation
$$
\Delta{\normal}+2p(X){\normal}=-2\Lambda\nabla H(X)
$$
%where $\Delta=\partial_{tt}+\partial_{\theta\theta}$,
where $\Lambda:=|X_{t}\wedge X_{\theta}|$ and $p(X)$ is the \emph{density function} given by
$$
p(X):=2\Lambda H^{2}(X)-\Lambda K(X)-\Lambda\nabla H(X)\cdot{\normal}~\!.
$$
In our case $H\equiv 1$ (hence $\nabla H\equiv 0$), $\Lambda=|(X_{a})_{t}\wedge(X_{a})_{\theta}|=x_{a}^{2}$, and the Gaussian curvature is given (in terms of the Gaussian notation (\ref{eq:Gaussian-notation})) by
$$
K(X_{a})=\frac{{\L}{\N}-{\M}^{2}}{{\E}{\G}-{\F}^{2}}=1-\frac{\gamma_{a}^{2}}{x_{a}^{4}}~\!.
$$
Hence
$$
p(X_{a})=x_{a}^{2}+\frac{\gamma_{a}^{2}}{x_{a}^{2}}~\!.
\quad\square
$$

\begin{Lemma}\label{L:w-0a}
The functions ${\normal}_{a}\cdot Z_{a}$ where $Z_{a}:=\dfrac{\partial X_{a}}{\partial a}$ and each component of ${\normal}_{a}\wedge X_{a}$ satisfy (\ref{eq:L=0}).
\end{Lemma}

\Proof
Recall that
$$
X_{a}=\left[\begin{array}{c}x_{a}\cos\theta\\ x_{a}\sin\theta\\ z_{a}\end{array}\right]\quad\text{and}\quad{\normal}_{a}=\left[\begin{array}{c}-z'_{a}x_{a}^{-1}\cos\theta\\ -z'_{a}x_{a}^{-1}\sin\theta\\ x'_{a}x_{a}^{-1}\end{array}\right].
$$
Differentiating (\ref{eq:CMC-X}) with respect to $a$ we obtain that
\begin{equation}
\label{eq:Delta-Z}
\Delta Z_{a}=2\left[(Z_{a})_{t}\wedge(X_{a})_{\theta}+(X_{a})_{t}\wedge(Z_{a})_{\theta}\right]~\!.
\end{equation}
Moreover, using (\ref{eq:Delta-Z}), (\ref{eq:conformality}) and (\ref{eq:CMC}), we compute
$$
({\normal}_{a})_{t}\cdot(Z_{a})_{t}=\frac{z'_{a}-2x_{a}^{2}}{x_{a}}\frac{\partial x_{a}}{\partial a}~,\quad({\normal}_{a})_{\theta}\cdot(Z_{a})_{\theta}=-\frac{z'_{a}}{x_{a}}\frac{\partial x_{a}}{\partial a}~\!.
$$
Hence, by Lemma \ref{L:N-equation}
\begin{equation*}
\begin{split}
\Delta ({\normal}_{a}\cdot Z_{a})&=(\Delta{\normal}_{a})\cdot Z_{a}+2({\normal}_{a})_{t}\cdot(Z_{a})_{t}+2({\normal}_{a})_{\theta}\cdot(Z_{a})_{\theta}+{\normal}_{a}\cdot \Delta Z_{a}\\
&=-2 p_{a}{\normal}_{a}\cdot Z_{a}-4x_{a}\frac{\partial x_{a}}{\partial a}+2{\normal}_{a}\cdot\left[(Z_{a})_{t}\wedge(X_{a})_{\theta}+(X_{a})_{t}\wedge(Z_{a})_{\theta}\right]~\!.
\end{split}
\end{equation*}
An explicit computation yields that
$$
{\normal}_{a}\cdot\left[(Z_{a})_{t}\wedge(X_{a})_{\theta}+(X_{a})_{t}\wedge(Z_{a})_{\theta}\right]=2x_{a}\frac{\partial x_{a}}{\partial a}
$$
and then the conclusion follows for ${\normal}_{a}\cdot Z_{a}$. Secondly, using again Lemma \ref{L:N-equation}, we have that
$$
\Delta({\normal}_{a}\wedge X_{a})+2p_{a}{\normal}_{a}\wedge X_{a}=
2({\normal}_{a})_{t}\wedge(X_{a})_{t}+2({\normal}_{a})_{\theta}\wedge(X_{a})_{\theta}+2{\normal}_{a}\wedge \Delta X_{a}~\!.
$$
By (\ref{eq:CMC-X}) one has that ${\normal}_{a}\wedge \Delta X_{a}=0$. Moreover, by (\ref{eq:useful})--(\ref{eq:useful2}),
$
({\normal}_{a})_{t}\wedge(X_{a})_{t}+({\normal}_{a})_{\theta}\wedge(X_{a})_{\theta}=\big({\normal}_{a}\wedge(X_{a})_{t}\big)_{t}+\big({\normal}_{a}\wedge(X_{a})_{\theta}\big)_{\theta}-{\normal}_{a}\wedge \Delta X_{a}=(X_{a})_{\theta t}-(X_{a})_{t\theta}=0
$
and this completes the proof.
\qed

Considering the scalar function ${\normal}_{a}\cdot Z_{a}$ and the components of the vectors ${\normal}_{a}$ and ${\normal}_{a}\wedge X_{a}$, from Lemmata \ref{L:N-equation} and \ref{L:w-0a}, one plainly obtains:

\begin{Lemma}
\label{L:eigenfunctions}
For $j=0,1$, a mapping $\phi\colon\R\to\R$ is a $2\tau_{a}$-periodic solution of $\phi''+(2p_{a}(t)-j^{2})\phi=0$ if and only if $\phi\in\mathrm{span}\{w_{a,j}^{+},w_{a,j}^{-}\}$ where
\begin{equation}
\label{eq:wa0+}
w_{a,0}^{+}=\frac{\partial X_{a}}{\partial a}\cdot N_{a}=-\frac{z'_{a}}{x_{a}}\frac{\partial x_{a}}{\partial a}+\frac{x'_{a}}{x_{a}}\frac{\partial z_{a}}{\partial a}~\!,\quad w_{a,0}^{-}=\frac{x'_{a}}{x_{a}}~\!,\quad w_{a,1}^{+}=\frac{z'_{a}}{x_{a}}~\!,\quad w_{a,1}^{-}=x'_{a}+\frac{z_{a}z'_{a}}{x_{a}}~\!.
\end{equation}
In particular $w_{a,j}^{+}$ are even functions, whereas $w_{a,j}^{-}$ are odd.
\end{Lemma}

%Set
%$$
%\widetilde{\mathfrak{L}}'_{a}\varphi(t,\theta):=\sum_{j=0}^{1}\chi_{j}(\theta)\widetilde{\mathfrak{L}}_{a,j}\phi_{j}(t)\quad\text{and}\quad \widetilde{\mathfrak{L}}''_{a}\varphi(t,\theta):=\sum_{j=2}^{\infty}\chi_{j}(\theta)\widetilde{\mathfrak{L}}_{a,j}\phi_{j}(t).
%$$
%Moreover set
%$$
%\begin{array}{ll}
%\mathscr{X}'_{a}:=\{\varphi\in \mathscr{X}_{a}~|~\varphi(t,\theta)=\textstyle\sum_{j=0}^{1}\phi_{j}(t)\chi_{j}(\theta)\}~\!,&\quad \mathscr{X}''_{a}:=\{\varphi\in \mathscr{X}_{a}~|~\varphi(t,\theta)=\textstyle\sum_{j=2}^{\infty}\phi_{j}(t)\chi_{j}(\theta)\}~\!,\\
%\mathscr{Y}'_{a}:=\{\varphi\in \mathscr{Y}_{a}~|~\varphi(t,\theta)=\textstyle\sum_{j=0}^{1}\phi_{j}(t)\chi_{j}(\theta)\}~\!,&\quad \mathscr{Y}''_{a}:=\{\varphi\in \mathscr{Y}_{a}~|~\varphi(t,\theta)=\textstyle\sum_{j=2}^{\infty}\phi_{j}(t)\chi_{j}(\theta)\}~\!.
%\end{array}
%$$

The first main result of this section is the following.

\begin{Theorem}
\label{T:ker}
One has that $\mathrm{ker}(\mathfrak{L}_{a})=\mathrm{span}\{w_{a,0},w_{a,1}\}$ where
\begin{equation}
\label{eq:eigenfunctions}
w_{a,0}(t,\theta)=w_{a,0}^{+}(t)~\!,\qquad w_{a,1}(t,\theta)=w_{a,1}^{+}(t)\sin\theta
\end{equation}
and $w_{a,j}^{+}(t)$ are the functions given in Lemma \ref{L:eigenfunctions}.
\end{Theorem}

\Proof
Let $\varphi\in \mathscr{X}_{a}$ satisfy (\ref{eq:L=0}). Because of the periodicity and symmetry with respect to $\theta$ of the functions in $\mathscr{X}_{a}$, we can express $\varphi\in \mathscr{X}_{a}$ by its Fourier series
\begin{equation*}
%\label{eq:decomposition}
\varphi(t,\theta)=\sum_{j=0}^{\infty}\phi_{j}(t)\chi_{j}(\theta)
\end{equation*}
where
\begin{equation}
\label{eq:chi}
\chi_{j}(\theta)=\left\{\begin{array}{ll}\frac{1}{\sqrt{2\pi}}&\text{for $j=0$}\\
\frac{1}{\sqrt{\pi}}\sin(j\theta)&\text{for $j=1,2,...$}\end{array}
\right.\quad\text{and}\quad\phi_{j}(t)=\int_{-\pi}^{\pi}\varphi(t,\theta)\chi_{j}(\theta)~\!d\theta~\!.
\end{equation}
Then
\begin{equation*}
%\label{eq:L-decomposition}
(\Delta+2p_{a})\varphi(t,\theta)=\sum_{j=0}^{\infty}\chi_{j}(\theta)\mathfrak{L}_{a,j}\phi_{j}(t)\quad\text{where}\quad
\mathfrak{L}_{a,j}:=\partial_{tt}+\left(2 p_{a}(t)-j^{2}\right)~\!.
\end{equation*}
Hence
$$
0=\int_{-\pi}^{-\pi}(\Delta+2p_{a})\varphi(t,\theta)\chi_{j}(\theta)~\!d\theta={\mathfrak{L}}_{a,j}\phi_{j}(t)\quad\forall t\in\R,~\forall j\ge 0.
$$
Multiplying by $\phi_{j}$ and integrating one gets
$$
\int_{-\tau_{a}}^{\tau_{a}}|\phi_{j}'|^{2}+\left(j^{2}-2 p_{a}(t)\right)|\phi_{j}|^{2}~\!dt=0\quad\forall j\ge 0.
$$
Recalling the definition of $p_{a}:=x_{a}^{2}+\gamma_{a}^{2}x_{a}^{-2}$ and that $\min_{t} x_{a}(t)=\min_{z}\rho_{a}(z)=a$ and $\max_{t}x_{a}(t)=\max_{z}\rho_{a}(z)=1-a$, one has $\|p_{a}\|_{\infty}=\max_{s\in[a,1-a]}\left(s^{2}+\gamma_{a}^{2}s^{-2}\right)=a^{2}+(1-a)^{2}$. Therefore
\begin{equation*}
%\label{eq:ca-j-estimate}
j^{2}-2 p_{a}(t)\ge 4-2\left(a^{2}+(1-a)^{2}\right)>0\quad\forall j\ge 2~\!.
\end{equation*}
Hence one infers that $\phi_{j}=0$ for all $j\ge 2$. For $j=0,1$, according to Lemma \ref{L:eigenfunctions}, the equations ${\mathfrak{L}}_{a,j}\phi=0$ admit exactly two independent $2\tau_{a}$-periodic solutions which, up to a multiplicative constant, are $w_{a,j}^{\pm}$. In fact, only $w_{a,j}^{+}$ are even. Hence $\varphi$ is a linear combination of $w_{a,0}$ and $w_{a,1}$.
\qed

According to Theorem \ref{T:ker} it is convenient to introduce the decompositions:
\begin{equation*}
%\label{decomposition}
\begin{array}{l}
\mathscr{X}_{a}=\mathscr{X}'_{a}+\mathscr{X}''_{a}\\
\mathscr{Y}_{a}=\mathscr{Y}'_{a}+\mathscr{Y}''_{a}
\end{array}
\quad\text{where}\quad
\left\{\begin{array}{l}
\mathscr{X}'_{a}=\mathscr{Y}'_{a}=\mathrm{span}\{w_{a,0},w_{a,1}\}\\
\mathscr{X}''_{a}=\{\varphi\in\mathscr{X}_{a}~|~ \textstyle\int_{Q_{\tau_{a}}}\varphi w_{a,i}~\!dt~\!d\theta=0\text{ for }i=0,1\}\\
\mathscr{Y}''_{a}=\{f\in\mathscr{Y}_{a}~|~ \textstyle\int_{Q_{\tau_{a}}}fw_{a,i}~\!dt~\!d\theta=0\text{ for }i=0,1\}~\!.
\end{array}\right.
\end{equation*}

Now we show
\begin{Theorem}\label{T:linear-operator}
Fixing $f\in\mathscr{Y}_{a}$, problem
\begin{equation}
\label{eq:L=f}
\left\{\begin{array}{ll}\Delta\varphi+2p_{a}\varphi=f&\text{in }\R/_{2\tau_{a}}\times\R/_{2\pi}
\\
\varphi\in \mathscr{X}_{a}&\end{array}\right.
\end{equation}
has a solution if and only if $f\in\mathscr{Y}''_{a}$ and in such case the solution is unique in $\mathscr{X}''_{a}$.
\end{Theorem}

\Proof
If problem (\ref{eq:L=f}) has a solution $\varphi$, then, testing the equation with $w_{a,i}$ ($i=0,1$) and integrating by parts, one easily obtains that $f\in\mathscr{Y}''_{a}$. Assume now that $f\in\mathscr{Y}''_{a}$. A solution to (\ref{eq:L=f}) can be found as a minimizer of the problem
$$
\inf_{H_{a}}\mathcal{E}_{a}\quad\text{where}\quad\mathcal{E}_{a}(\varphi)=\int_{Q_{\tau_{a}}}\left(\frac{1}{2}|\nabla\varphi|^{2}-p_{a}\varphi^{2}+f\varphi\right)~\!dt~\!d\theta
$$
and $H_{a}$ is the Hilbert space obtained as the completion of $\mathscr{X}''_{a}$ with respect to the standard norm in $W^{1,2}(Q_{\tau_{a}})$. In fact, one can observe that $\mathcal{E}_{a}$ is {coercive} (because we detached the kernel of $\Delta+2p_{a}$ and any non zero eigenvalue is negative) and convex in $H_{a}$, and then the existence of a minimizer $\varphi\in H_{a}$ follows by well known arguments. Moreover the minimizer $\varphi$ is a weak solution of
$$
\Delta\varphi+2p_{a}\varphi-f=\lambda_{0}w_{a,0}+\lambda_{1}w_{a,1}\quad\text{in }\R/_{2\tau_{a}}\times\R/_{2\pi}
$$
for some Lagrange multipliers $\lambda_{0},\lambda_{1}\in\R$. In fact regularity theory applies and allows us to deduce that $\varphi\in\mathscr{X}_{a}$. Testing the equation with $w_{a,i}$ ($i=0,1$), integrating by parts, and exploiting the fact that $f\in\mathscr{Y}''_{a}$ and
$$
\int_{Q_{\tau_{a}}}w_{a,0}w_{a,1}~\!dt~\!d\theta=0~\!,
$$
one obtains that $\lambda_{0}=\lambda_{1}=0$, that is, $\varphi$ solves (\ref{eq:L=f}). The uniqueness of the solution in $\mathscr{X}''_{a}$ can be checked in a standard way.
\qed

Our next goal is to find uniform estimates on $\varphi\in\mathscr{X}''_{a}$ with respect to $f\in\mathscr{Y}''_{a}$, with constant independent on the parameter $a$. To this purpose, we introduce some weighted norms as follows. Fixing $\alpha\in(0,1)$, $\delta>0$, $\mu>0$, and $k\ge 0$ integer, for $f\in C^{k,\alpha}(\R/_{2\tau_{a}}\times\R/_{2\pi})$ set
\begin{equation}
\label{eq:weighted-norms}
\|f\|_{a,k,\mu}:=\sup_{s\in\R}x_{a}(s)^{-\mu}\|f\|_{C^{k,\alpha}([s-\delta,s+\delta]\times\S^{1})}
\end{equation}
where $x_{a}$ is given according to Lemma \ref{L:c-a} and $\|f\|_{C^{k,\alpha}}$ denotes the H\"{o}lder norm in $C^{k,\alpha}$ on the specified domain. We also denote $\|f\|_{C^{k}}$ the usual $C^{k}$ norm.

The second part of this Section is focused on the discussion of the following result.
\begin{Theorem}
\label{T:a0}
Fix $\mu\in(1,2)$ and $\delta>0$. Then there exist $a_{0}\in \big(0,\frac{1}{2}\big)$ and $C_{1}>0$, depending only on $\mu$ and $\delta$, such that
$$
\|\varphi\|_{a,2,\mu}\le C_{1}\|f\|_{a,0,\mu}\quad\forall f\in \mathscr{Y}''_{a}~\!,\quad\forall a\in(0,a_{0})
$$
where $\varphi$ is the unique solution in $\mathscr{X}''_{a}$ of (\ref{eq:L=f}), given by Theorem \ref{T:linear-operator}.
\end{Theorem}

We need some preliminary technical results, collected in the next lemma.

\begin{Lemma}
\label{L:xwp}
One has that $w_{a,0}^{+}(t)\to 1-t\tanh t$ and $w_{a,1}^{+}(t)\to\sech t$ in $C^{2}_{loc}(\R)$, as $a\to 0$. Moreover, for every $a>0$ small enough
\begin{equation}
\label{eq:wa-bound}
|w_{a,0}^{+}(t)|\le C(1+|t|)\quad\text{and}\quad 0<w_{a,1}^{+}(t)\le 1\quad\forall t\in\R
\end{equation}
for some constant $C>0$ independent of $a$. In addition, for every $t\in\R$ one has that $\lim_{a\to 0}p_{a}(t)=(\sech t)^{2}$, $\lim_{a\to 0}p_{a}(t+\tau_{a})=(\sech t)^{2}$ and, if $a_{n}\to 0$, $t_{n}\to\infty$, and $\tau_{a_{n}}-t_{n}\to\infty$, then $\lim_{n\to\infty}p_{a_{n}}(t+t_{n})=0$ for every $t\in\R$.
\end{Lemma}

\Proof
Since $w_{a,1}^{+}(t)=z'_{a}(t)/x_{a}(t)$ (see Lemma \ref{L:eigenfunctions}), Lemma \ref{L:xa-limit} implies that $w_{a,1}^{+}(t)\to\sech t$ point-wise, as $a\to 0$. In fact, we can prove convergence in $C^{2}_{loc}$. To this aim, let us observe that, by Lemma \ref{L:eigenfunctions}, the mapping $w_{a,1}^{+}=w_{a}$ solves
$$
\left\{\begin{array}{l}
w_{a}''+(2p_{a}-1)w_{a}=0\\ w_{a}(0)=1\\ w_{a}(t+2\tau_{a})=w_{a}(t)=w_{a}(-t)\end{array}\right.
$$
Using the definition of $w_{a}$, after differentiation, and taking into account that $x'_{a}<0$ on $(0,\tau_{a})$, one infers that $\max w_{a}=1$, attained at $0$ and $\tau_{a}$, $\min w_{a}=2\sqrt\gamma_{a}$ is attained at some $t_{a}\in(0,\tau_{a})$ where
\begin{equation}
\label{eq:xa(ta)}
x_{a}(t_{a})^{2}=\gamma_{a}~\!.
\end{equation}
In particular the second estimate in (\ref{eq:wa-bound}) is proved. Moreover $w_{a}$ is strictly decreasing in $[0,t_{a}]$ and strictly increasing in $[t_{a},\tau_{a}]$. Using (\ref{eq:conformality}) and the definitions of $p_{a}$ and $w_{a}$, one has that
\begin{equation}
\label{eq:wa}
\left\{\begin{array}{l}
{w}''_{a}=(1+4\gamma_{a}){w}_{a}-2{w}_{a}^{3}\\
({w}'_{a})^{2}=(1-{w}_{a}^{2})({w}_{a}^{2}-4\gamma_{a})~\!.
\end{array}\right.
\end{equation}
Thus one obtains a uniform bound in $C^{3}$ for $w_{a}$, which implies that $w_{a,1}^{+}(t)\to\sech t$ as $a\to 0$, in $C^{2}_{loc}$. Moreover, since $w_{a}(0)=w_{a}(\tau_{a})$ one also infers that $w_{a}$ is $\tau_{a}$-periodic, and since the first equation in (\ref{eq:wa}) is invariant under reflection, it holds that $w_{a}\big(\frac{\tau_{a}}{2}+t\big)=w_{a}\big(\frac{\tau_{a}}{2}-t\big)$. In particular the minimum point of $w_{a}$ is $t_{a}=\frac{\tau_{a}}{2}$. Hence, from (\ref{eq:xa(ta)}) it follows that
\begin{equation*}
%\label{eq:xa(taua/2)}
x_{a}\left(\frac{\tau_{a}}{2}\right)^{2}=\gamma_{a}~\!.
\end{equation*}
Now let us consider $y_{a}:=w_{a,0}^{+}$. Since $x_{a}(0)=1-a$, $z_{a}(0)=0$, using also Lemma \ref{L:c-a}, we compute
\begin{equation}
\label{eq:ya0}
y_{a}(0)=1~\!.
\end{equation}
Moreover, after some computations, where (\ref{eq:conformality}), (\ref{eq:CMC}) and again Lemma \ref{L:c-a} are often used, one can check that
\begin{equation}
\label{eq:ya'}
y'_{a}=(1-2a)\frac{x_{a}}{x'_{a}}-\left(x_{a}-\frac{\gamma_{a}}{x_{a}}\right)\frac{z'_{a}}{x'_{a}}y_{a}~\!.
\end{equation}
Furthermore, recall that, by Lemma \ref{L:eigenfunctions},
\begin{equation}
\label{eq:ya''}
y_{a}''=-2p_{a}y_{a}
\end{equation}
and that $y_{a}$ is a $2\tau_{a}$-periodic, even function. From (\ref{eq:ya0})--(\ref{eq:ya''}) one deduces that $y_{a}$ is strictly decreasing in $[0,\tau_{a}]$ and there exists $\overline{\tau}_{a}\in(0,\tau_{a})$ such that $y_{a}$ is positive and concave in $[0,\overline{\tau}_{a})$ and is negative and convex in $(\overline{\tau}_{a},\tau_{a}]$. Hence
\begin{equation}
\label{eq:ya-bound}
|y_{a}(t)|\le 1-2a+|y'_{a}(\overline{\tau}_{a})|~\!|t|\quad\forall t\in\R~\!.
\end{equation}
Taking the limit of (\ref{eq:ya'}) as $a\to 0$ one obtains the equation
$$
y'=\frac{y}{(\sinh t)(\cosh t)}-\frac{\cosh t}{\sinh t}
$$
whose general solution is $y(t)=K\tanh t+e^{-t}\sech t-t\tanh t$. In fact the only even solution occurs when $K=1$, yielding
$$
y_{0}(t)=1-t\tanh t~\!.
$$
Observe that $y_{0}(0)=1$, $y_{0}$ is strictly decreasing in $[0,\infty)$ and there exists $\overline{\tau}_{0}>1$ such that $y_{0}(\overline{\tau}_{0})=0$, $y_{0}$ is positive and concave in $[0,\overline{\tau}_{0})$ and is negative and convex in $(\overline{\tau}_{0},\infty)$. Moreover, by (\ref{eq:ya0})--(\ref{eq:ya''}), one has that $y_{a}\to y_{0}$, namely $w_{a,0}^{+}\to 1-t\tanh t$, in $C^{2}_{loc}(\R)$. In particular
\begin{equation}
\label{eq:ya'-bound}
y'_{a}(\overline{\tau}_{a})\to y'_{0}(\overline{\tau}_{0})~\!.
\end{equation}
Then (\ref{eq:ya-bound}) and (\ref{eq:ya'-bound}) imply the first estimate in (\ref{eq:wa-bound}) for $a>0$ small enough.

The limit of $p_{a}(t)$ as $a\to 0$ immediately follows from the definition of $p_{a}$, given in (\ref{eq:pa-def}). In order to study the remaining cases, let us set
$$
\tilde{w}_{a}(t):=w_{a}(t+\tau_{a})~\!.
$$
We observe that $\tilde{w}_{a}$ solves (\ref{eq:wa}), too. Hence a uniform bound in $C^{3}$ for $\tilde{w}_{a}$ holds true. Then one plainly obtains that any sequence $a_{n}\to 0$ admits a subsequence, still denoted $a_{n}$, such that $\tilde{w}_{a_{n}}\to \tilde{w}$ as $n\to\infty$ uniformly on compact sets, where $\tilde{w}$ solves
$$
\left\{\begin{array}{l}
\tilde{w}''=\tilde{w}-2\tilde{w}^{3}\\
(\tilde{w}')^{2}=\tilde{w}^{2}-\tilde{w}^{4}.
\end{array}\right.
$$
The only solutions to this problem are $\sech(t+\tilde{t})$ for some $\tilde{t}\in\R$ or the null function. Since $\tilde{w}_{a}(0)=1$ for every $a$, it must be $\tilde{w}(t)=\sech t$. The result for ${p}_{a}(t+\tau_{a})$ follows from the identity ${p}_{a}(t+\tau_{a})=\tilde{w}_{a}(t)^{2}-2\gamma_{a}$. Finally, consider the case of sequences $a_{n}\to 0$ and $t_{n}\to\infty$ with $\tau_{a_{n}}-t_{n}\to\infty$. Hence necessarily $t_{n}\in(0,\tau_{a_{n}})$. Repeating the same argument as before for the sequence
$$
\overline{w}_{a_{n}}(t):=w_{a_{n}}(t+t_{n})
$$
one concludes that, up to a subsequence, $\overline{w}_{a_{n}}\to \overline{w}$ as $n\to\infty$ uniformly on compact sets, where $\overline{w}$ can be $\sech(t+\bar{t})$ for some $\bar{t}\in\R$ or the null function. We claim that $\overline{w}_{a_{n}}(0)\to 0$ and thus, necessarily, $\overline{w}=0$. To check the claim, we argue by contradiction, assuming that $\lim w_{a_{n}}(t_{n})=\alpha\in(0,1]$. Since $w_{a_{n}}(t)\to\sech t$ and $w_{a_{n}}(t-\tau_{a_{n}})\to\sech t$ uniformly on compact sets, and by assumption, $t_{n}\to\infty$ and $\tau_{a_{n}}-t_{n}\to\infty$, it must be $\alpha<1$. Observe that, by (\ref{eq:wa}), $w_{a_{n}}$ is decreasing in $[0,t_{a_{n}}]$ and increasing in $[t_{a_{n}},\tau_{a_{n}}]$. Assume that $t_{n}\in[0,t_{a_{n}}]$ for infinitely many indices, and note that in fact $t_{n}\in(0,t_{a_{n}})$, because $\alpha\in(0,1)$. Fixing $\beta\in(\alpha,1)$, since $w_{a_{n}}(t)\to\sech t$ uniformly on compact sets, there exists $s_{n}\in(0,t_{n})$ such that $(s_{n})$ is a converging sequence and $w_{a_{n}}(s_{n})=\beta$. Then $w_{a_{n}}$ is decreasing in $[s_{n},t_{n}]$ and, by (\ref{eq:wa}),
\begin{equation*}
\begin{split}
w_{a_{n}}(s_{n})-w_{a_{n}}(t_{n}))&=\int_{s_{n}}^{t_{n}}\sqrt{(1-w_{a_{n}}(s)^{2})(w_{a_{n}}(s)^{2}-4\gamma_{a_{n}})}~\!ds\\
&\ge(t_{n}-s_{n})\inf_{r\in[\alpha_{n},\beta]}\sqrt{(1-r^{2})(r^{2}-4\gamma_{a_{n}})}
\end{split}
\end{equation*}
where $\alpha_{n}=w_{a_{n}}(t_{n})$. This yields a contradiction, because the right hand side diverges whereas the left hand side is bounded. Hence the claim is proved. A similar argument holds if $t_{n}\in[t_{a_{n}},\tau_{a_{n}}]$ for infinitely many $n$'s. In conclusion $\overline{w}_{a_{n}}\to 0$ and hence ${p}_{a_{n}}(t+t_{n})=\overline{w}_{a_{n}}(t)^{2}-2\gamma_{a_{n}}\to 0$, too.
\qed

\begin{Lemma}
\label{L:weight}
For $\mu>0$ one has that
\begin{gather}
\label{eq:weight-local}
x_{a}(t)^{-\mu}\le e^{\delta\mu}x_{a}(s)^{-\mu}\quad\forall s\in\R~\!,\quad\forall t\in[s-\delta,s+\delta]~\!,\\
\label{eq:weight-global}
x_{a}(t_{0})^{-\mu}\le e^{\mu t}x_{a}(t+t_{0})^{-\mu}\quad\forall t_{0}\in[0,\tau_{a}]~\!,\quad\forall t\ge -t_{0}~\!.
\end{gather}
\end{Lemma}
\Proof
By (\ref{eq:conformality}) $|x'_{a}|\le x_{a}$ and this implies that the mapping $e^{-t}x_{a}(t)$ is decreasing whereas $e^{t}x_{a}(t)$ is increasing. Then, since $\mu>0$,
\begin{equation}
\label{eq:weight-partial}
x_{a}(t)^{-\mu}\le e^{\mu|t-s|}x_{a}(s)^{-\mu}\quad\forall s,t\in\R~\!.
\end{equation}
In particular (\ref{eq:weight-local}) holds true. Similarly, one observes when $t_{0}\ge 0$ one has that $e^{t_{0}}x_{a}(t_{0})\ge x_{a}(0)$. Moreover the mapping $e^{-t}x_{a}(t+t_{0})$ is decreasing and then
$$
e^{-t}x_{a}(t+t_{0})\le e^{-t_{0}}x_{a}(0)\le x_{a}(t_{0})\quad\forall t\ge -t_{0}
$$
and hence (\ref{eq:weight-global}) holds true, as $\mu>0$.
\qed

\noindent
\emph{Proof of Theorem \ref{T:a0}.}
Arguing by contradiction, assume that there exist sequences of Delaunay parameters $a_{n}\to 0$ and of mappings $f_{n}\in \mathscr{Y}''_{a_{n}}$ with $\|f_{n}\|_{a_{n},0,\mu}\le 1$ such that $\|\varphi_{n}\|_{a_{n},2,\mu}\to\infty$, where $\varphi_{n}$ is the unique solution in $\mathscr{X}''_{a_{n}}$ of
\begin{equation*}
%\label{eq:varphi-fn}
\Delta\varphi+2p_{a_{n}}(t)\varphi=f_{n}~\!.
\end{equation*}
Firstly let us prove that
\begin{equation}
\label{eq:diverges}
\|x_{a_{n}}^{-\mu}\varphi_{n}\|_{C^{0}([0,\tau_{a_{n}}]\times\S^{1})}\to\infty~\!.
\end{equation}
Indeed, by regularity theory, one knows that
\begin{equation}
\label{eq:regularity}
\|\varphi_{n}\|_{C^{2,\alpha}([s-\delta,s+\delta]\times\S^{1})}\le C\big(\|\varphi_{n}\|_{C^{0}([s-\delta,s+\delta]\times\S^{1})}+\|f_{n}\|_{C^{0,\alpha}([s-\delta,s+\delta]\times\S^{1})}\big)
\end{equation}
where $C$ is a positive constant independent of $s$ and $n$. From (\ref{eq:weight-local}) and (\ref{eq:regularity}) it follows that
\begin{equation*}
%\label{eq:weighted-regularity}
x_{a_{n}}(s)^{-\mu}\|\varphi_{n}\|_{C^{2,\alpha}([s-\delta,s+\delta]\times\S^{1})}\le Ce^{\delta\mu}\|x_{a_{n}}^{-\mu}\varphi_{n}\|_{C^{0}([s-\delta,s+\delta]\times\S^{1})}+ Cx_{a_{n}}(s)^{-\mu}\|f_{n}\|_{C^{0,\alpha}([s-\delta,s+\delta]\times\S^{1})}~\!.
\end{equation*}
Hence, by the definition of the norms (\ref{eq:weighted-norms}) and since we deal with mappings $2\tau_{a}$-periodic and even with respect to $t$, we infer that
$$
\|\varphi_{n}\|_{a_{n},2,\mu}\le Ce^{\delta\mu}\|x_{a_{n}}^{-\mu}\varphi_{n}\|_{C^{0}([0,\tau_{a_{n}}]\times\S^{1})}+C\|f_{n}\|_{a_{n},0,\mu}
$$
and thus (\ref{eq:diverges}) is proved. Let $(t_{n},\theta_{n})\in[0,\tau_{a_{n}}]\times\S^{1}$ be such that
$$
x_{a_{n}}(t_{n})^{-\mu}|\varphi_{n}(t_{n},\theta_{n})|=\|x_{a_{n}}^{-\mu}\varphi_{n}\|_{C^{0}(\R\times\S^{1})}=:A_{n}~\!.
$$
Up to a subsequence, three alternative cases can occur:
\begin{itemize}
\item[\textbf{1.}] $t_{n}\to t_{0}$ for some $t_{0}\in[0,\infty)$, or
\item[\textbf{2.}] $\tau_{a_{n}}-t_{n}\to t_{0}$ for some $t_{0}\in[0,\infty)$, or
\item[\textbf{3.}] $t_{n}\to\infty$ and $\tau_{a_{n}}-t_{n}\to \infty$.
\end{itemize}
%We can also assume that $\theta_{n}\to\theta_{0}\in[0,2\pi)$.
We are going to obtain a contradiction in each case.
\medskip

\noindent
\textbf{Case 1.} $t_{n}\to t_{0}\in[0,\infty)$. Since $x_{a_{n}}(t)\to\sech t$ uniformly on compact sets, there exists a constant $C>0$ such that $x_{a_{n}}(t_{n})^{-\mu}\le C$ for every $n\in\mathbb{N}$. Define
$$
\tilde\varphi_{n}(t,\theta)=A_{n}^{-1}x_{a_{n}}(t_{n})^{-\mu}\varphi_{n}(t,\theta)\quad\text{and}\quad \tilde{f}_{n}(t,\theta)=A_{n}^{-1}x_{a_{n}}(t_{n})^{-\mu}f_{n}(t,\theta)~\!.
$$
Then
\begin{eqnarray}
\label{eq:first-1}
&|\tilde\varphi_{n}(t,\theta)|\le Cx_{a_{n}}(t)^{\mu}\\
\label{eq:first-2}
&\tilde\varphi_{n}\in\mathscr{X}''_{a_{n}}\\
\label{eq:first-3}
&\Delta\tilde\varphi_{n}+2p_{a_{n}}\tilde\varphi_{n}=\tilde{f}_{n}\\
\label{eq:first-4}
&|\tilde\varphi_{n}(t_{n},\theta_{n})|=1~\!.%,~~\tilde f_{n}\in \mathscr{Y}''_{a_{n}}.
\end{eqnarray}
In addition,
\begin{equation}
\label{eq:f-estimate-1}
\|\tilde f_{n}\|_{a_{n},0,\mu}=A_{n}^{-1}x_{a_{n}}^{-\mu}(t_{n})\|f_{n}\|_{a_{n},0,\mu}\le C A_{n}^{-1}.
\end{equation}
By regularity theory (use (\ref{eq:regularity}) with $\varphi_{n}$ and $f_{n}$ replaced by $\tilde{\varphi}_{n}$ and $\tilde{f}_{n}$, respectively), using (\ref{eq:weight-local}), (\ref{eq:first-1}) and (\ref{eq:f-estimate-1}), we obtain
\begin{equation*}
%\label{eq:phi-estimate-1}
x_{a_{n}}(s)^{-\mu}\|\tilde{\varphi}_{n}\|_{C^{2,\alpha}([s-\delta,s+\delta]\times\S^{1})}\le Ce^{\delta\mu}+CA_{n}^{-1}\quad\forall s\in[-\tau_{a_{n}}+\delta,\tau_{a_{n}}-\delta].
\end{equation*}
Hence, since $A_{n}\to\infty$, up to a subsequence, $\tilde\varphi_{n}$ converges in $C^{2}_{loc}(\R\times\S^{1})$ to some $\tilde\varphi$ which, by (\ref{eq:first-1})--(\ref{eq:f-estimate-1}) and by Lemmata \ref{L:xa-limit} and \ref{L:xwp}, solves
\begin{eqnarray}
\label{eq:phi-tilde-1-1}
%&|\tilde\varphi(t,\theta)|\le(\sech t_{0})^{\mu}(\cosh t)^{\mu}\quad\forall (t,\theta)\in\R^{2}\\
&|\tilde\varphi(t,\theta)|\le C(\sech t)^{\mu}\quad\forall (t,\theta)\in\R\times\S^{1}\\
\label{eq:phi-tilde-1-2}
&\tilde\varphi(t,\cdot)=\tilde\varphi(-t,\cdot)~\!,~~\tilde\varphi(\cdot,\frac{\pi}{2}-\theta)=\tilde\varphi(\cdot,\frac{\pi}{2}+\theta)\\
\label{eq:phi-tilde-1-3}
&\Delta\tilde\varphi+2(\sech t)^{2}\tilde\varphi=0\quad\text{in }\R\times\S^{1}\\
\label{eq:phi-tilde-1-4}
&|\tilde\varphi(t_{0},\theta_{0})|=1\quad\text{for some }\theta_{0}\in\S^{1}~\!.
\end{eqnarray}
We will show that (\ref{eq:phi-tilde-1-1})--(\ref{eq:phi-tilde-1-3}) imply that $\tilde\varphi\equiv 0$, contradicting (\ref{eq:phi-tilde-1-4}). To obtain that $\tilde\varphi\equiv 0$ we consider the Fourier decomposition of $\tilde\varphi$, which, in view of (\ref{eq:phi-tilde-1-2}), takes the form
\begin{equation}
\label{eq:Fourier-tilde}
\tilde\varphi(t,\theta)=\sum_{j=0}^{\infty}\psi_{j}(t)\chi_{j}(\theta)\quad\forall (t,\theta)\in\R\times\S^{1}
\end{equation}
where, for every integer $j\ge 0$ and for every $t\in\R$
$$
\psi_{j}(t)=\int_{-\pi}^{\pi}\chi_{j}(\theta)\tilde\varphi(t,\theta)~\!d\theta
$$
and $\chi_{j}$ is defined in (\ref{eq:chi}). Because of (\ref{eq:phi-tilde-1-2}) and (\ref{eq:phi-tilde-1-3}), $\psi_{j}(t)$ solves
\begin{equation}
\label{eq:psi-j-problem}
\left\{\begin{array}{l}
\psi_{j}''+\big[2(\sech t)^{2}-j^{2}\big]\psi_{j}=0\\ \psi_{j}(t)=\psi_{j}(-t)~\!.\end{array}\right.
\end{equation}
Let us prove that
\begin{equation}\label{j>1}
\psi_{j}(t)=0~~\forall t\in\R,~\forall j\ge 2~\!.
\end{equation}
In fact, fixing $j\ge 2$, assume by contradiction that $\psi_{j}(0)>0$. Then by (\ref{eq:psi-j-problem}), $\psi_{j}$ is convex, has a global minimum at $0$, in particular $\psi_{j}(t)>0$ for all $t\ge 0$. On the other hand, (\ref{eq:phi-tilde-1-1}) it follows that $|\psi_{j}(t)|\le C(\sech t)^{\mu}$, which contradicts the previous statement. Therefore we conclude that
$\psi_{j}(0)=0$ and, since $\psi'_{j}(0)=0$ ($\psi_{j}$ is even), by uniqueness of the Cauchy problem, (\ref{j>1}) holds true. For $j=0$ the only solutions to (\ref{eq:psi-j-problem}) are $c_{0}(1-t\tanh t)$, whereas for $j=1$ they are $c_{1}\sech t$. Hence
$$%\begin{equation}
%\label{eq:quasi-zero}
\tilde\varphi(t,\theta)=c_{0}(1-t\tanh t)+c_{1}(\sech t)\sin\theta~\!.
$$%\end{equation}
The uniform bound (\ref{eq:phi-tilde-1-1}) and the condition $\mu>1$ plainly imply that $c_{0}=c_{1}=0$.
\medskip

\noindent
\textbf{Case 2.} $\tau_{a_{n}}-t_{n}\to t_{0}\in[0,\infty)$.
%Since $w_{a_{n}}$ is $\tau_{a_{n}}$-periodic and even, and $w_{a_{n}}(t)\to\sech t$ uniformly on compact sets, also in this case there exists a constant $C>0$ such that $w_{a_{n}}(t_{n})^{-\mu}\le C$ for every $n\in\mathbb{N}$.
Define
$$
\tilde\varphi_{n}(t,\theta)=A_{n}^{-1}w_{a_{n}}(t_{n})^{-\mu}\varphi_{n}(t+\tau_{a_{n}},\theta)\quad\text{and}\quad \tilde{f}_{n}(t,\theta)=A_{n}^{-1}w_{a_{n}}(t_{n})^{-\mu}f_{n}(t+\tau_{a_{n}},\theta)~\!.
$$
Using (\ref{eq:weight-global}) with $t_{0}=\tau_{a_{n}}$ and $a=a_{n}$, and (\ref{eq:weight-partial}), one has that
\begin{equation}
\label{eq:one-side-exp}
|\tilde\varphi_{n}(t,\theta)|\le A_{n}^{-1}\left(\frac{x_{a_{n}}(t+\tau_{a_{n}})}{x_{a_{n}}(t_{n})}\right)^{\mu}\|x_{a_{n}}^{-\mu}\varphi_{n}\|_{C^{0}(\R\times\S^{1})}\le e^{\mu t}e^{\mu(\tau_{a_{n}}-t_{n})}\le Ce^{\mu t}\quad\forall (t,\theta)\in Q_{\tau_{a_{n}}}
\end{equation}
and with similar computations, $\|\tilde{f}_{n}\|_{a_{n},0,\mu}\le C A_{n}^{-1}$. Moreover
\begin{eqnarray}
%\label{eq:first-2-1}
%&|\tilde\varphi_{n}(t,\theta)|\le C w_{a_{n}}(t)^{\mu}\\
\nonumber%\label{eq:first-2-2}
&\tilde\varphi_{n}\in\mathscr{X}''_{a_{n}}\\
\label{eq:first-2-3}
&\Delta\tilde\varphi_{n}+2p_{a_{n}}(t+\tau_{a_{n}})\tilde\varphi_{n}=\tilde{f}_{n}\\
\nonumber%\label{eq:first-2-4}
&|\tilde\varphi_{n}(\tau_{a_{n}}-t_{n},\theta_{n})|=1~\!.%,~~\tilde f_{n}\in \mathscr{Y}''_{a_{n}}.
\end{eqnarray}
Since $\tilde\varphi_{n}$ is even with respect to $t$, from (\ref{eq:one-side-exp}) it follows that
\begin{equation*}
%\label{eq:first-2-1}
|\tilde\varphi_{n}(t,\theta)|\le C e^{-\mu|t|}\quad\forall (t,\theta)\in Q_{\tau_{a_{n}}}~\!.
\end{equation*}
Now we can argue as in Step 1. The only difference is the coefficient $p_{a_{n}}(t+\tau_{a_{n}})$ in the differential equation (\ref{eq:first-2-3}). Actually, by Lemma \ref{L:xwp}, $p_{a_{n}}(t+\tau_{a_{n}})\to(\sech t)^{2}$. Hence, up to a subsequence, $\tilde\varphi_{n}$ converges in $C^{2}_{loc}(\R\times\S^{1})$ to some $\tilde\varphi$ which solves (\ref{eq:phi-tilde-1-1})--(\ref{eq:phi-tilde-1-4}). Then we obtain a contradiction as in Step 1.
\medskip

\noindent
\textbf{Case 3.} $t_{n}\to\infty$ and $\tau_{a_{n}}-t_{n}\to \infty$, and $\theta_{n}\to\theta_{0}$. Define
$$
\tilde\varphi_{n}(t,\theta)=A_{n}^{-1}x_{a_{n}}(t_{n})^{-\mu}\varphi_{n}(t+t_{n},\theta)\quad\text{and}\quad \tilde{f}_{n}(t,\theta)=A_{n}^{-1}x_{a_{n}}(t_{n})^{-\mu}f_{n}(t+t_{n},\theta)~\!.
$$
Hence
\begin{gather}
\label{eq:equation-case3}
\Delta\tilde\varphi_{n}+2p_{a_{n}}(t+t_{n})\tilde\varphi_{n}=\tilde{f}_{n}\text{~~~in }[-t_{n},\tau_{a_{n}}-t_{n}]\times\S^{1}
\\
\tilde\varphi_{n}\big(\cdot,\tfrac{\pi}{2}-\theta\big)=\tilde\varphi_{n}\big(\cdot,\tfrac{\pi}{2}+\theta\big)\quad\forall\theta\in\S^{1}
\\
\label{eq:normalization-case3}
|\tilde\varphi_{n}(0,\theta_{n})|=1~\!.%,~~\tilde f_{n}\in \mathscr{Y}''_{a_{n}}.
\end{gather}
Using the definition of $A_{n}$ and (\ref{eq:weight-global}) with $t_{n}$ and $a_{n}$ instead of $t_{0}$ and $a$, respectively, we deduce that
\begin{equation}
\label{eq:exponential-growth}
|\tilde{\varphi}_{n}(t,\theta)|\le e^{\mu t}\quad\forall(t,\theta)\in[-t_{n},\tau_{a_{n}}-t_{n}]\times\S^{1}~\!.
\end{equation}
Now we need a $C^{2,\alpha}_{loc}$ estimate on the sequence $(\tilde\varphi_{n})$ in order to pass to the limit in (\ref{eq:equation-case3})--(\ref{eq:normalization-case3}). By regularity theory, for every $s\in[-t_{n}+\delta,\tau_{a_{n}}-t_{n}-\delta]$ we have that
\begin{equation}
\label{eq:C2alpha-case3}
\|\tilde\varphi_{n}\|_{C^{2,\alpha}([s-\delta,s+\delta]\times\S^{1})}\le C\left(\|\tilde\varphi_{n}\|_{C^{0}([s-\delta,s+\delta]\times\S^{1})}+ \|\tilde{f}_{n}\|_{C^{0,\alpha}([s-\delta,s+\delta]\times\S^{1})}\right)
\end{equation}
where $C$ is a positive constant independent of $n$ and $s$. Using (\ref{eq:weight-global}) we estimate
\begin{eqnarray}
\nonumber
\|\tilde{f}_{n}\|_{C^{0,\alpha}([s-\delta,s+\delta]\times\S^{1})}&=&A_{n}^{-1}x_{a_{n}}(t_{n})^{-\mu}\|f_{n}\|_{C^{0,\alpha}([s+t_{n}-\delta,s+t_{n}+\delta]\times\S^{1})}\\
\nonumber
&\le& A_{n}^{-1}e^{\mu|s|}x_{a_{n}}(s+t_{n})^{-\mu}\|f_{n}\|_{C^{0,\alpha}([s+t_{n}-\delta,s+t_{n}+\delta]\times\S^{1})}\\
\label{eq:f-tilde-estimate-3}
&\le& A_{n}^{-1}e^{\mu|s|}\|f_{n}\|_{a_{n},0,\mu}\le A_{n}^{-1}e^{\mu|s|}~\!.
\end{eqnarray}
Then, by (\ref{eq:exponential-growth}), from (\ref{eq:C2alpha-case3}) it follows that for every $s\in[-t_{n}+\delta,\tau_{a_{n}}-t_{n}-\delta]$
\begin{equation}
\label{eq:phi-tilde-estimate-3}
\|\tilde\varphi_{n}\|_{C^{2,\alpha}([s-\delta,s+\delta]\times\S^{1})}\le C(e^{\mu(|s|+\delta)}+A_{n}^{-1}e^{\mu|s|})\le C_{1}e^{\mu|s|}~\!,
\end{equation}
with $C_{1}$ positive constant independent of $n$ and $s$. Using (\ref{eq:f-tilde-estimate-3}) and (\ref{eq:phi-tilde-estimate-3}), we can extract subsequences $\tilde{f}_{n}\to 0$ in $C^{0}_{loc}(\R\times\S^{1})$ and $\tilde{\varphi}_{n}\to\tilde{\varphi}$ in $C^{2}_{loc}(\R\times\S^{1})$, with $\tilde{\varphi}$ satisfying
\begin{eqnarray}
\label{eq:phi-tilde-3-1}
&\Delta\tilde{\varphi}=0\quad\text{in~~}\R\times\S^{1}&\\
\label{eq:phi-tilde-3-2}
&\tilde\varphi\big(\cdot,\tfrac{\pi}{2}-\theta\big)=\tilde\varphi\big(\cdot,\tfrac{\pi}{2}+\theta\big)\quad\forall\theta\in\S^{1}&
\\
\label{eq:phi-tilde-3-3}
&|\tilde{\varphi}(0,\theta_{0})|=1\quad\text{for some }\theta_{0}\in\S^{1}~\!.&
\end{eqnarray}
Moreover, in view of (\ref{eq:exponential-growth}),
\begin{equation}
\label{eq:phi-tilde-3-4}
|\tilde{\varphi}(t,\theta)|\le e^{\mu t}\quad\forall(t,\theta)\in\R\times\S^{1}~\!.
\end{equation}
We will show that (\ref{eq:phi-tilde-3-1}), (\ref{eq:phi-tilde-3-2}) and (\ref{eq:phi-tilde-3-4}) imply that $\tilde\varphi\equiv 0$, in contradiction with (\ref{eq:phi-tilde-3-3}). To this goal, we consider the Fourier decomposition of $\tilde\varphi$, given by (\ref{eq:Fourier-tilde}) where now the $j$-th component
$$
\psi_{j}(t)=\int_{-\pi}^{\pi}\chi_{j}(\theta)\tilde\varphi(t,\theta)~\!d\theta
$$
satisfies $\psi_{j}''=j^{2}\psi_{j}$ on $\R$. Then $\psi_{j}(t)=A_{j}e^{jt}+B_{j}e^{-jt}$ for some $A_{j},B_{j}\in\R$, for $j\in\mathbb{N}$, and $\psi_{0}(t)=A_{0}+B_{0}t$. In fact, (\ref{eq:phi-tilde-3-4}) implies also that $|\psi_{j}(t)|\le Ce^{\mu t}$ for every $t\in\R$. Since $\mu\in(1,2)$, the only possibility is that $A_{j}=B_{j}=0$. This holds for every $j\ge 0$. Hence $\tilde\varphi=0$. Thus we proved that none of the three alternative cases on the sequence $(t_{n})$ can occur and this completes the proof.
\qed

\section{The reduced problem}\label{S:reduced}

We start to tackle the problem
\begin{equation}
\label{eq:main-equation}
\left\{\begin{array}{l}
\varphi\in\mathscr{X}_{a}\\
\mathfrak{M}(X_{\eps,a}+\varphi{\normal}_{\eps,a})=H(X_{\eps,a}+\varphi{\normal}_{\eps,a})\end{array}\right.
\end{equation}
where $X_{\eps,a}$ is the parametrization of the generalized toroidal unduloid defined in (\ref{eq:eps-toroidal-unduloid}) with corresponding Gauss map ${\normal}_{\eps,a}$, $\mathfrak{M}$ denotes the mean curvature operator, given by (\ref{eq:mean-curvature-def}), and $H\colon\R^{3}\setminus\{0\}\to\R$ is a prescribed, radial mapping of class $C^{2}$ satisfying $(H_{1})$--$(H_{2})$. The goal of this Section is to show the following result.

\begin{Theorem}
\label{T:reduction}
Under the above assumptions on $H$, there exist $a_{0}>0$, $\eps_{0}>0$ and mappings $(\eps,a)\mapsto\lambda^{0}_{\eps,a}\in\R$, $(\eps,a)\mapsto\lambda^{1}_{\eps,a}\in\R$, $(\eps,a)\mapsto\varphi_{\eps,a}\in\mathscr{X}_{a}$, defined for $\eps\in(0,\eps_{0})$ and $a\in(a,a_{0})$, continuous with respect to $a$ and such that
\begin{equation}
\label{eq:reduced-problem}
\begin{array}{c}
\mathfrak{M}(X_{\eps,a}+\varphi_{\eps,a}{\normal}_{\eps,a})-H(X_{\eps,a}+\varphi_{\eps,a}{\normal}_{\eps,a})=\frac{1}{2x_{a}^{2}}\left(\lambda^{0}_{\eps,a}w_{a,0}+\lambda^{1}_{\eps,a} w_{a,1}\right)\\
\displaystyle\int_{Q_{\tau_{a}}}\varphi_{\eps,a}w_{a,0}~\!dt~\!d\theta=\int_{Q_{\tau_{a}}}\varphi_{\eps,a}w_{a,1}~\!dt~\!d\theta=0
\end{array}
\end{equation}
where $w_{a,0}$ and $w_{a,1}$ are defined in (\ref{eq:eigenfunctions}). Moreover there exists $R>0$ such that for every $\eps\in(0,\eps_{0})$ and $a\in(0,a_{0})$
\begin{equation}
\label{eq:varphi-small}
\|\varphi_{\eps,a}\|_{2,a,\mu}\le R\eps^{\widetilde\gamma}\quad\text{where}\quad\widetilde\gamma=\min\{1,\gamma\}~\!.
\end{equation}
\end{Theorem}

\begin{Remark}
\label{R:point-wise-bound}
Condition (\ref{eq:varphi-small}) is quite strong. In fact, according to the definition (\ref{eq:weighted-norms}) of the weighted norm, one has that
$$
|\varphi_{\eps,a}(s,\theta)|\le R\eps^{\widetilde\gamma}x_{a}(s)^{\mu}
\quad\forall(s,\theta)\in\R^{2}
$$
and the same point-wise bound holds true also for the first and second derivatives of $\varphi_{\eps,a}$.
\end{Remark}

\Proof
Let $a_{0}>0$ given by Theorem \ref{T:a0}. We take $a\in(0,a_{0})$ and $\eps>0$ small enough (with a bound independent of $a$) in order that (\ref{eq:Leps-a-estimate}) and Lemmata \ref{L:M(X)-eps} and \ref{L:Xt} can be applied. Let us introduce the operators $\mathfrak{B}^{0}_{\eps,a},\mathfrak{B}^{1}_{\eps,a}\colon\mathscr{X}_{a}\to\mathscr{Y}_{a}$ defined by
\begin{gather}
\label{eq:MM-decomposition}
\mathfrak{M}(X_{\eps,a}+\varphi{\normal}_{\eps,a})=
\mathfrak{M}(X_{\eps,a})+\frac{1}{2x_{a}^{2}}\mathfrak{L}_{a}\varphi+\frac{1}{2x_{a}^{2}}(\mathfrak{L}_{\eps,a}-\mathfrak{L}_{a})\varphi+\frac{1}{2x_{a}^{2}}\mathfrak{B}^{0}_{\eps,a}(\varphi)\\
\nonumber
H(X_{\eps,a}+\varphi{\normal}_{\eps,a})=H(X_{\eps,a})+[\nabla H(X_{\eps,a})\cdot\normal_{\eps,a}]\varphi+\mathfrak{B}^{1}_{\eps,a}(\varphi)
\end{gather}
where $\mathfrak{L}_{\eps,a}$ and $\mathfrak{L}_{a}$ are defined by (\ref{eq:L-eps-a}) and (\ref{eq:pa-def}), respectively. Hence problem (\ref{eq:main-equation}) is equivalent to find a solution $\varphi\in \mathscr{X}_{a}$ of
\begin{equation}
\label{eq:Lphi=F(phi)}
\mathfrak{L}_{a}\varphi=\mathfrak{F}_{\eps,a}(\varphi)
\end{equation}
where
$$
\mathfrak{F}_{\eps,a}(\varphi):=-2x_{a}^{2}[\mathfrak{M}(X_{\eps,a})-H(X_{\eps,a})]-(\mathfrak{L}_{\eps,a}-\mathfrak{L}_{a})\varphi+2x_{a}^{2}[\nabla H(X_{\eps,a})\cdot\normal_{\eps,a}]\varphi-\mathfrak{B}^{0}_{\eps,a}(\varphi)+2x_{a}^{2}\mathfrak{B}^{1}_{\eps,a}(\varphi)~\!.
$$
We aim to rewrite (\ref{eq:Lphi=F(phi)}) as a fixed point problem. To this purpose we consider the inverse of $\mathfrak{L}_{a}$, which actually is well defined only on $\mathscr{Y}_{a}''$. Let us denote $\widehat{\mathfrak{F}}_{\eps,a}(\varphi)$ the projection of ${\mathfrak{F}}_{\eps,a}(\varphi)$ on $\mathscr{Y}_{a}''$ and let us set
$$
\mathfrak{T}_{\eps,a}(\varphi):=\mathfrak{L}_{a}^{-1}\circ\widehat{\mathfrak{F}}_{\eps,a}(\varphi)~\!.
$$
Thus $\mathfrak{T}_{\eps,a}$ maps $\mathscr{X}''_{a}$ into $\mathscr{X}''_{a}$. Moreover a function $\varphi\in \mathscr{X}_{a}$ solves
\begin{equation}
\label{eq:fixed-point-problem}
\left\{\begin{array}{l}\varphi\in \mathscr{X}''_{a}\\ \varphi=\mathfrak{T}_{\eps,a}(\varphi)\end{array}\right.
\end{equation}
if and only if it solves (\ref{eq:reduced-problem}). We are going to prove that there exists $R>0$ and $\eps_{0}>0$ such that for every $a\in(0,a_{0})$ and for every $\eps\in(0,\eps_{0})$, denoting
$$
\mathscr{B}_{\eps,a}=\{\varphi\in\mathscr{X}''_{a}~|~\|\varphi\|_{a,2,\mu}\le R\eps^{\widetilde\gamma}\}\quad\text{where~~}\widetilde\gamma=\min\{\gamma,1\}~\!,
$$
one has
\begin{itemize}
\item[(i)] $\mathfrak{T}_{\eps,a}(\mathscr{B}_{\eps,a})\subset\mathscr{B}_{\eps,a}$,
\item[(ii)] $\mathfrak{T}_{\eps,a}$ is a contraction in $\mathscr{B}_{\eps,a}$ with a Lipschitz constant independent of $a$.
\end{itemize}
If (i) and (ii) hold true, then the contraction principle implies the existence of a solution to (\ref{eq:fixed-point-problem}) for every $a\in(0,a_{0})$ and $\eps\in(0,\eps_{0})$, with continuous dependence on $a$, namely the conclusion of the theorem.

In order to prove properties (i) and (ii) we need some preliminary estimates that will be stated and proved below. In such estimates we will denote by $C$ a positive constant which may change from line to line but does not depend on $a\in(0,a_{0})$ or $\eps>0$ (sufficiently small, for every estimate, independently of $a$). Moreover we will often use the property on the product of H\"older continuous functions
$$
\|fg\|_{C^{0,\alpha}(\overline\Omega)}\le\|f\|_{C^{0,\alpha}(\overline\Omega)}\|g\|_{C^{0,\alpha}(\overline\Omega)}\quad\forall f,g\in C^{0,\alpha}(\overline\Omega)~\!.
$$
\textbf{Estimate 1.}
There exists a constant $C_{H}>0$, depending only on $H$, such that
\begin{equation}
\label{eq:main-term}
\|x_{a}^{2}(\mathfrak{M}(X_{\eps,a})-H(X_{\eps,a}))\|_{a,0,\mu}\le C_{H}\eps^{\widetilde\gamma}~~\forall\eps>0\text{ small,}~~\forall a\in(0,a_{0})~\!,\quad\text{where~~}\widetilde\gamma=\min\{1,\gamma\}~\!.
\end{equation}
We can estimate
\begin{equation*}
\begin{split}
\|x_{a}^{2}(\mathfrak{M}(X_{\eps,a})-H(X_{\eps,a}))\|_{a,0,\mu}&=\sup_{s\in[0,\tau_{a}]}x_{a}(s)^{-\mu}\|x_{a}^{2}(\mathfrak{M}(X_{\eps,a})-H(X_{\eps,a}))\|_{C^{0,\alpha}([s-\delta,s+\delta]\times\S^{1})}\\
&\le\|\mathfrak{M}(X_{\eps,a})-H(X_{\eps,a}))\|_{C^{0,\alpha}(Q_{\tau_{a}})}\sup_{s\in[0,\tau_{a}]}x_{a}(s)^{-\mu}\|x_{a}\|^{2}_{C^{0,\alpha}([s-\delta,s+\delta])}~\!.
\end{split}
\end{equation*}
Since $|x'_{a}|\le x_{a}$ (see (\ref{eq:conformality})), using (\ref{eq:weight-local}), for every $s\in\R$ one has
$$
x_{a}(s)^{-\mu}\|x_{a}\|^{2}_{C^{0,\alpha}([s-\delta,s+\delta])}\le x_{a}(s)^{-\mu}\|x_{a}\|^{2}_{C^{1}([s-\delta,s+\delta])}\le 4x_{a}(s)^{-\mu}\|x_{a}\|^{2}_{C^{0}([s-\delta,s+\delta])}\le C
$$
because $\mu<2$ and $0<x_{a}<1$.
Now we estimate
$$
\|\mathfrak{M}(X_{\eps,a})-H(X_{\eps,a}))\|_{C^{0,\alpha}(Q_{\tau_{a}})}\le
\|\mathfrak{M}(X_{\eps,a})-1\|_{C^{0,\alpha}(Q_{\tau_{a}})}+
\|H(X_{\eps,a})-1\|_{C^{0,\alpha}(Q_{\tau_{a}})}~\!.
$$
According to Lemma \ref{L:M(X)-eps}, there exists a constant $C_{0}>0$ such that
\begin{equation}
\label{eq:M-1}
\|\mathfrak{M}(X_{\eps,a})-1\|_{C^{0,\alpha}(\R^{2})}\le C_{0}\eps\quad\forall\eps>0\text{~~small,}\quad\forall a\in\big(0,\tfrac{1}{2}\big]~\!.
\end{equation}
Moreover
$$
|H(X_{\eps,a})-1|\le\frac{C}{|X_{\eps,a}|^{\gamma}}\le C\eps^{\gamma}
$$
because $|X_{\eps,a}|\ge|\eps^{-1}-x_{a}|$ and $0<x_{a}<1$. In addition we know that
\begin{equation}
\label{eq:nablaH}
|\nabla H(X)|\le\frac{C}{|X|^{\gamma+1}}\quad\forall X\in\R^{3}
\end{equation}
and then
$$
\left|\nabla_{(t,\theta)}\left(H(X_{\eps,a})-1\right)\right|^{2}\le|\nabla_{X}H(X_{\eps,a})|^{2}\left(|(X_{\eps,a})_{t}|^{2}+|(X_{\eps,a})_{\theta}|^{2}\right)\le \frac{Cx_{a}^{2}}{|X_{\eps,a}|^{2(\gamma+1)}}
$$
from which
$$
\left|\nabla_{(t,\theta)}\left(H(X_{\eps,a})-1\right)\right|\le C\eps^{\gamma+1}~\!.
$$
Therefore
\begin{equation}
\label{eq:H-1}
\|H(X_{\eps,a})-1\|_{C^{0,\alpha}(Q_{\tau_{a}})}\le C\eps^{\gamma}
\end{equation}
for some constant $C>0$ depending only on $H$. Thus (\ref{eq:main-term}) follows from (\ref{eq:M-1}) and (\ref{eq:H-1}).
\medskip

\noindent
\textbf{Choice of $R$.} We fix
\begin{equation*}
%\label{eq:R-fixed}
R=2C_{1}C_{H}
\end{equation*}
where $C_{1}$ is given by Theorem \ref{T:a0} and $C_{H}$ is given by the estimate (\ref{eq:main-term}).
\medskip

\noindent
\textbf{Estimate 2.} There exists a constant $C>0$ independent of $\eps$ and $a$, such that
\begin{equation}
\label{eq:estimate-2}
\|(\mathfrak{L}_{\eps,a}-\mathfrak{L}_{a})\varphi\|_{a,0,\mu}\le C\eps\|\varphi\|_{a,2,\mu}\quad\forall\eps>0\text{ small,}\quad\forall a\in(0,a_{0})~\!,\quad\forall\varphi\in\mathscr{X}_{a}~\!.
\end{equation}
By Lemma \ref{L:L-na} and by the definition of $\mathfrak{L}_{a}$ (see (\ref{eq:pa-def})), arguing as for estimate 1, one has that
\begin{equation*}
\begin{split}
\|(\mathfrak{L}_{\eps,a}-\mathfrak{L}_{a})\varphi\|_{a,0,\mu}&\le\|b_{\eps,a}-1\|_{C^{0,\alpha}(\R^{2})}\|\varphi_{tt}\|_{a,0,\mu}+\|c_{\eps,a}-p_{a}\|_{C^{0,\alpha}(\R^{2})}\|\varphi\|_{a,0,\mu}\\
&\qquad+\|d_{\eps,a}\|_{C^{0,\alpha}(\R^{2})}\|\varphi_{t}\|_{a,0,\mu}+\|e_{\eps,a}\|_{C^{0,\alpha}(\R^{2})}\|\varphi_{\theta}\|_{a,0,\mu}\le C\eps\|\varphi\|_{a,2,\mu}~\!.
\end{split}
\end{equation*}
\medskip

\noindent
\textbf{Estimate 3.} There exists a constant $C>0$ independent of $\eps$ and $a$, such that
\begin{equation}
\label{eq:estimate-3}
\|x_{a}^{2}(\nabla H(X_{\eps,a})\cdot\normal_{\eps,a})\varphi\|_{a,0,\mu}\le C\eps^{\gamma+1}\|\varphi\|_{a,0,\mu}\quad\forall\eps>0\text{ small,}\quad\forall a\in(0,a_{0})~\!,\quad\forall\varphi\in\mathscr{X}_{a}~\!.
\end{equation}
Arguing as for the first estimate, we plainly arrive to
$$
\|x_{a}^{2}(\nabla H(X_{\eps,a})\cdot\normal_{\eps,a})\varphi\|_{a,0,\mu}\le\|x_{a}^{2}\|_{C^{0,\alpha}(\R^{2})}\|\nabla H(X_{\eps,a})\cdot\normal_{\eps,a}\|_{C^{0,\alpha}(\R^{2})}\|\varphi\|_{a,0,\mu}~\!.
$$
Since $|x'_{a}|\le x_{a}<1$, one has that
\begin{equation}
\label{eq:xa2}
\|x_{a}^{2}\|_{C^{0,\alpha}(\R^{2})}\le C
\end{equation}
with $C$ independent of $a$. Moreover, using (\ref{eq:nablaH}), we also obtain
\begin{equation}
\label{eq:est3-0}
|\nabla H(X_{\eps,a})\cdot\normal_{\eps,a}|\le |\nabla H(X_{\eps,a})|\le\frac{C}{|X_{\eps,a}|^{\gamma+1}}\le C\eps^{\gamma+1}
\end{equation}
for $\eps>0$ small and uniformly in $a$. Moreover
\begin{equation}
\label{eq:nablaHXN}
\left(\nabla H(X_{\eps,a})\cdot\normal_{\eps,a}\right)_{t}=\left(H''(|X_{\eps,a}|)-\frac{H'(|X_{\eps,a}|)}{|X_{\eps,a}|}\right)|X_{\eps,a}|_{t}\frac{X_{\eps,a}\cdot N_{\eps,a}}{|X_{\eps,a}|}+H'(|X_{\eps,a}|)\frac{(X_{\eps,a}\cdot N_{\eps,a})_{t}}{|X_{\eps,a}|}~\!.
\end{equation}
Now we apply (\ref{eq:nablaH}),
\begin{equation*}
%\label{eq:H''}
|H''(|X|)|\le \frac{C}{|X|^{\gamma+2}}\quad\forall X\in\R^{3}~\!,
\end{equation*}
Lemma \ref{L:Xt}, and elementary estimates (in particular note that $\big||X_{\eps,a}|_{t}\big|\le\left|(X_{\eps,a})_{t}\right|$), and we find that
\begin{equation}
\label{eq:est3-1}
\left|\left(\nabla H(X_{\eps,a})\cdot\normal_{\eps,a}\right)_{t}\right|\le \frac{C}{|X_{\eps,a}|^{\gamma+2}}\le C\eps^{\gamma+2}
\end{equation}
for $\eps>0$ small and uniformly in $a$. A similar estimate holds true also for $\left|\left(\nabla H(X_{\eps,a})\cdot\normal_{\eps,a}\right)_{\theta}\right|$.  Then (\ref{eq:xa2}), (\ref{eq:est3-0}) and (\ref{eq:est3-1}) and its analogous for the derivative with respect to $\theta$ imply (\ref{eq:estimate-3}).
\medskip

\noindent
\textbf{Estimate 4.} There exists a constant $C>0$ independent of $\eps$ and $a$ such that
\begin{equation}
\label{eq:estimate-4}
\|\mathfrak{B}^{0}_{\eps,a}(\varphi)\|_{a,0,\mu}\le C\eps^{\widetilde\gamma}\|\varphi\|_{a,0,\mu}\quad\forall\eps>0\text{ small,}\quad\forall a\in(0,a_{0})~\!,\quad\forall\varphi\in\mathscr{B}_{\eps,a}~\!.
\end{equation}
By definition, $\mathfrak{B}^{0}_{\eps,a}$ is a second order differential operator of $\varphi$ whose expansion about $\varphi=0$ does not contain constant or linear terms (with respect to $\varphi$ and its first and second derivatives), and the coefficients of every term have partial derivatives bounded uniformly in $a$ and $\eps$. Therefore
$$
\|\mathfrak{B}^{0}_{\eps,a}(\varphi)\|_{C^{0,\alpha}([s-\delta,s+\delta]\times\S^{1})}\le C\|\varphi\|^{2}_{C^{2,\alpha}([s-\delta,s+\delta]\times\S^{1})}
$$
with $C>0$ independent of $\eps$, $a$ and $s$. Then
\begin{equation*}
\begin{split}
\|\mathfrak{B}^{0}_{\eps,a}(\varphi)\|_{a,0,\mu}&=\sup_{s\in[0,\tau_{a}]}x_{a}(s)^{-\mu}\|\mathfrak{B}^{0}_{\eps,a}(\varphi)\|_{C^{0,\alpha}([s-\delta,s+\delta]\times\S^{1})}\\
&\le C\sup_{s\in[0,\tau_{a}]}x_{a}(s)^{-2\mu}\|\varphi\|^{2}_{C^{2,\alpha}([s-\delta,s+\delta]\times\S^{1})}\le C\|\varphi\|_{a,2,\mu}^{2}
\end{split}
\end{equation*}
because $0<x_{a}\le 1$. Then (\ref{eq:estimate-4}) follows for $\varphi\in\mathscr{B}_{\eps,a}$.
\medskip

\noindent
\textbf{Estimate 5.}
There exists a constant $C>0$ independent of $\eps$ and $a$ such that
\begin{equation}
\label{eq:estimate-5}
\|x_{a}^{2}\mathfrak{B}^{1}_{\eps,a}(\varphi)\|_{a,0,\mu}\le C\eps^{\gamma+1}\|\varphi\|_{a,0,\mu}\quad\forall\eps>0\text{ small,}\quad\forall a\in(0,a_{0})~\!,\quad\forall\varphi\in\mathscr{B}_{\eps,a}~\!.
\end{equation}
Since
$$
\|x_{a}^{2}\mathfrak{B}^{1}_{\eps,a}(\varphi)\|_{a,0,\mu}\le \|x_{a}^{2}\|_{C^{0,\alpha}(\R^{2})}\sup_{s\in[0,\tau_{a}]}x_{a}(s)^{-\mu}\|\mathfrak{B}^{1}_{\eps,a}(\varphi)\|_{C^{0,\alpha}([s-\delta,s+\delta]\times\S^{1})}~\!,
$$
thanks to (\ref{eq:xa2}), it is enough to estimate the $C^{0,\alpha}$ norm of $\mathfrak{B}^{1}_{\eps,a}(\varphi)$ in terms of the corresponding norm of $\varphi$. Observe that
$$
\mathfrak{B}^{1}_{\eps,a}(\varphi)=\varphi~\!\mathfrak{I}_{\eps,a}(\varphi)\quad\text{where}\quad\mathfrak{I}_{\eps,a}(\varphi)=\int_{0}^{1}[\nabla H(X_{\eps,a}+r\varphi N_{\eps,a})-\nabla H(X_{\eps,a})]\cdot N_{\eps,a}~\!dr~\!.
$$
By (\ref{eq:nablaH}) and since $|X_{\eps,a}+r\varphi N_{\eps,a}|\ge \eps^{-1}-R\eps^{\widetilde\gamma}$ for every $r\in[0,1]$, we obtain
\begin{equation}
\label{eq:Iea}
|\mathfrak{I}_{\eps,a}(\varphi)|\le C\eps^{\gamma+1}~\!.
\end{equation}
Moreover, we can estimate $\left|\left(\nabla H(X_{\eps,a}+r\varphi N_{\eps,a})\cdot N_{\eps,a}\right)_{t}\right|$ by means of an expansion like (\ref{eq:nablaHXN}), with $X_{\eps,a}+r\varphi N_{\eps,a}$ instead of $X_{\eps,a}$. Writing
\begin{gather*}
((X_{\eps,a}+r\varphi N_{\eps,a})\cdot N_{\eps,a})_{t}=(X_{\eps,a}\cdot N_{\eps,a})_{t}+r\varphi_{t}~\!,
\\
|X_{\eps,a}+r\varphi N_{\eps,a}|_{t}=\frac{X_{\eps,a}+r\varphi N_{\eps,a}}{|X_{\eps,a}+r\varphi N_{\eps,a}|}\cdot((X_{\eps,a})_{t}+r\varphi_{t} N_{\eps,a})+r\varphi\frac{(X_{\eps,a}\cdot N_{\eps,a})_{t}}{|X_{\eps,a}+r\varphi N_{\eps,a}|}~\!,
\end{gather*}
observing that $|\varphi|,|\varphi_{t}|\le\|\varphi\|_{a,1,\mu}\le R\eps^{\widetilde\gamma}$, and using Lemma \ref{L:Xt} and elementary estimates, we obtain
$$
\left|\left(\nabla H(X_{\eps,a}+r\varphi N_{\eps,a})\cdot\normal_{\eps,a}\right)_{t}\right|\le C\eps^{\gamma+2}
$$
and a similar estimate for the derivative with respect to $\theta$. From these and (\ref{eq:Iea}) we infer that $\|\mathfrak{I}_{\eps,a}(\varphi)\|_{C^{1}(\R^{2})}\le C\eps^{\gamma+1}$, hence $\|\mathfrak{B}^{1}_{\eps,a}(\varphi)\|_{C^{0,\alpha}([s-\delta,s+\delta]\times\S^{1})}\le C\eps^{\gamma+1}\|\varphi\|_{C^{0,\alpha}([s-\delta,s+\delta]\times\S^{1})}$ and then (\ref{eq:estimate-5}).
\medskip

\noindent
\textbf{Proof of (i).} By Theorem \ref{T:a0} and thanks to Estimates 1 to 5, for every $\varphi\in\mathscr{B}_{\eps,a}$ one has that
$$
\|\mathfrak{T}_{\eps,a}(\varphi)\|_{a,2\mu}\le C_{1}\|\mathfrak{F}_{\eps,a}(\varphi)\|_{a,2,\mu}\le C_{1}C_{H}\eps^{\widetilde\gamma}+o(1)\|\varphi\|_{a,2,\mu}
$$
with $o(1)\to 0$ as $\eps\to 0$ uniformly with respect to $a$ and $\varphi$. Therefore there exists $\eps_{0}$ (independent of $a$) such that for every $\eps\in(0,\eps_{0})$ and for every $a\in(0,a_{0})$ one has $\|\mathfrak{T}_{\eps,a}(\varphi)\|_{a,2\mu}\le 2 C_{1}C_{H}\eps^{\widetilde\gamma}$, namely, $\mathfrak{T}_{\eps,a}(\varphi)\in\mathscr{B}_{\eps,a}$, because of the choice of $R$.
\medskip

\noindent
\textbf{Estimate 6.} There exists $C>0$ independent of $\eps$ and $a$, such that for every $a\in(0,a_{0})$ and for every $\eps\in(0,1]$ one has
\begin{equation}
\label{eq:estimate-6}
\|\mathfrak{B}^{0}_{\eps,a}(\varphi_{1})-\mathfrak{B}^{0}_{\eps,a}(\varphi_{2})\|_{a,0,\mu}\le C\eps^{\widetilde\gamma}\|\varphi_{1}-\varphi_{2}\|_{a,2,\mu}\quad\forall\varphi_{1},\varphi_{2}\in\mathscr{B}_{\eps,a}~\!.
\end{equation}
It is a consequence of the estimate
$$
\|\mathfrak{B}^{0}_{\eps,a}(\varphi_{1})-\mathfrak{B}^{0}_{\eps,a}(\varphi_{2})\|_{C^{0,\alpha}}\le C\left(\|\varphi_{1}\|_{C^{2,\alpha}}+\|\varphi_{2}\|_{C^{2,\alpha}}\right)\|\varphi_{1}-\varphi_{2}\|_{C^{2,\alpha}}
$$
where $C>0$ is a constant independent of $\eps$, $a$, $\varphi_{1}$ and $\varphi_{2}$. Such an estimate is justified by the same argument used for Estimate 4. Reasoning in a similar way, one gets (\ref{eq:estimate-6}).
\medskip

\noindent
\textbf{Estimate 7.}
There exists $C>0$ independent of $\eps$ and $a$, such that for every $a\in(0,a_{0})$ and for every $\eps\in(0,1]$ one has
\begin{equation*}
%\label{eq:estimate-7}
\|x_{a}^{2}\mathfrak{B}^{1}_{\eps,a}(\varphi_{1})-x_{a}^{2}\mathfrak{B}^{1}_{\eps,a}(\varphi_{2})\|_{a,0,\mu}\le \eps^{\gamma+1}\|\varphi_{1}-\varphi_{2}\|_{a,0,\mu}\quad\forall\varphi_{1},\varphi_{2}\in\mathscr{B}_{\eps,a}~\!.
\end{equation*}
Arguing as for Estimate 5, it is enough to find an estimate of the form
$$
\|\mathfrak{B}^{1}_{\eps,a}(\varphi_{1})-\mathfrak{B}^{1}_{\eps,a}(\varphi_{2})\|_{C^{0,\alpha}([s-\delta,s+\delta]\times\S^{1})}\le C\eps^{\gamma+1}\|\varphi_{1}-\varphi_{2}\|_{C^{0,\alpha}([s-\delta,s+\delta]\times\S^{1})}~\!.
$$
Since
$$
\mathfrak{B}^{1}_{\eps,a}(\varphi_{1})-\mathfrak{B}^{1}_{\eps,a}(\varphi_{2})=(\varphi_{1}-\varphi_{2})\int_{0}^{1}[\nabla H(X_{\eps,a}+(r\varphi_{1}+(1-r)\varphi_{2})N_{\eps,a})-\nabla H(X_{\eps,a})]\cdot N_{\eps,a}~\!dr~\!,
$$
reasoning as before, it suffices to find a uniform bound in $C^{0,\alpha}(\R^{2})$ for $\nabla H(X_{\eps,a}+\varphi N_{\eps,a})\cdot N_{\eps,a}$ with $\varphi\in\mathscr{B}_{\eps,a}$ and this can be obtained by means of already considered estimates.
\medskip

\noindent
\textbf{Proof of (ii).} By Theorem \ref{T:a0} and thanks to Estimates 2, 3, 6 and 7, we have that
$$
\|\mathfrak{T}_{\eps,a}(\varphi_{1})-\mathfrak{T}_{\eps,a}(\varphi_{2})\|_{a,2,\mu}\le C_{0}\|\mathfrak{F}_{\eps,a}(\varphi_{1})-\mathfrak{F}_{\eps,a}(\varphi_{2})\|_{a,0,\mu}\le o(1)\|\varphi_{1}-\varphi_{2}\|_{a,2,\mu}\quad\forall\varphi_{1},\varphi_{2}\in\mathscr{B}_{\eps,a}
$$
with $o(1)\to 0$ as $\eps\to 0$ uniformly with respect to $a\in(0,a_{0})$,  and $\varphi_{1},\varphi_{2}\in\mathscr{B}_{\eps,a}$. In particular, taking $\eps>0$ small enough, $\mathfrak{T}_{\eps,a}$ turns out to be a contraction in $\mathscr{B}_{\eps,a}$ with a Lipschitz constant independent of $a$.
\qed

\section{The variational argument}
\label{S:variational}

The main result of this Section is:

\begin{Theorem}
\label{T:immersed}
Let $H\colon\R^{3}\setminus\{0\}\to\R$ be a radially symmetric mapping of class $C^{2}$ satisfying $(H_{1})$--$(H_{2})$. If $A<0$ then there exist $n_{A}\in\mathbb{N}$, $\underbar{b},\overline{b}>0$, and a sequence $a_{n}\in(0,a_{0})$ such that for every $n\ge n_{A}$
\begin{gather}
\label{eq:Xnan}
\mathfrak{M}(X_{n,a_{n}}+\varphi_{n,a_{n}}N_{n,a_{n}})=H(X_{n,a_{n}}+\varphi_{n,a_{n}}N_{n,a_{n}})\\
\label{an-behaviour}
a_{n}=\frac{b_{n}}{n^{\gamma}\log n}\text{~~with~~}\underbar{b}\le b_{n}\le\overline{b}
\end{gather}
where $X_{n,a}$ is the toroidal unduloid defined in (\ref{eq:toroidal-unduloid}), $N_{n,a}$ is its corresponding Gauss map, and $\varphi_{n,a}$ is the scalar mapping given by Theorem \ref{T:reduction}, taking $\eps=\eps_{n,a}=\pi/(nh_{a})$. Moreover
\begin{equation}
\label{eq:phi-n-weighted-estimate}
|\varphi_{n,a_{n}}|\le C\eps_{n,a_{n}}^{\widetilde\gamma}x_{a_{n}}^{\mu}\quad\text{on~~}\R^{2}~\!,
\end{equation}
where $\widetilde\gamma=\min\{1,\gamma\}$. Furthermore $n_{A}\to\infty$ as $A\to 0$.
\end{Theorem}

In other words, Theorem \ref{T:immersed} states the existence of a sequence of immersed tori with mean curvature $H$ which are close to a corresponding sequence of toroidal unduloids with a large number $n$ of moduli and whose neck-size $a_{n}$ becomes smaller and smaller with $n$ according to the relationship (\ref{an-behaviour}).

The starting point in the proof of Theorem \ref{T:immersed} is the result of Theorem \ref{T:reduction} which gives the existence of a sequence of  mappings $a\mapsto \varphi_{n,a}\in\mathscr{X}_{a}''$ depending on $a\in(0,a_{0})$ in a continuous way, and scalar mappings $a\mapsto\lambda_{n,a}^{i}$ ($i=0,1$) such that for $n\ge n_{0}$, with $n_{0}\in\mathbb{N}$ independent of $a$,
\begin{equation}
\label{eq:reduced-equation}
\mathfrak{M}(X_{n,a}+\varphi_{n,a}N_{n,a})=H(X_{n,a}+\varphi_{n,a}N_{n,a})=\frac{1}{2x_{a}^{2}}\left(\lambda_{n,a}^{0}w_{a,0}+\lambda_{n,a}^{1}w_{a,1}\right)
\end{equation}
where $w_{a,0}$ and $w_{a,1}$ are defined in (\ref{eq:eigenfunctions}).
Now, for every $n$ large enough we aim to find $a_{n}\in(0,a_{0})$ for which the coefficients $\lambda^{0}_{n,a}$ and $\lambda^{1}_{n,a}$ in (\ref{eq:reduced-equation}) vanish when $a=a_{n}$. This will be reached by exploiting the variational character of the problem.

Firstly we introduce the area and the volume functionals for parametric surfaces of the type of the torus, i.e., (images of) doubly periodic mappings $X\colon\R/_{2\tau}\times\R/_{2\pi}\to\R^{3}$ of class $C^{2}$, such that $X_{t}\wedge X_{\theta}\ne 0$ everywhere. It is convenient to keep an arbitrary, but fixed period $2\tau$ with respect to $t$. The rectangle $Q_{\tau}:=[-\tau,\tau]\times[-\pi,\pi]$ will be called the parameter domain of the surface.

The area of a toroidal surface $\Sigma$ parameterized by $X$ and the algebraic volume enclosed by $\Sigma$ are computed, respectively, in terms of the integrals
$$
\mathcal{A}(X;Q_{\tau})=\mathcal{A}(X):=\int_{Q_{\tau}}|X_{t}\wedge X_{\theta}|~\!dt~\!d\theta~\!,\quad\mathcal{V}(X;Q_{\tau})=\mathcal{V}(X):=\frac{1}{3}\int_{Q_{\tau}}X\cdot X_{t}\wedge X_{\theta}~\!dt~\!d\theta~\!.
$$
As an example, and for future reference, let us compute the area of and the enclosed volume by the toroidal unduloid $\Sigma_{n,a}$ introduced in Section \ref{S:toroidal-unduloids}, as well as their expansions at $a=0$.

\begin{Lemma}
\label{L:area-volume-expansions}
There exist $C^{1}$ mappings $S,G,G_{n}\colon(0,a_{0})\to\R$ ($n\in\mathbb{N}$ large) with $C^{1}$ norms bounded uniformly in $n$, such that
\begin{gather}
\label{eq:area-expansion}
\mathcal{A}(X_{n,a})=4\pi n\left[1-a-\frac{a^{2}}{2}\log a+a^{2}S(a)+n^{-2}G_{n}(a)\right]\\
\label{eq:volume-expansion}
\mathcal{V}(X_{n,a})=-\frac{2\pi n}{3}\left[2-3a+a^{2}G(a)\right]~\!.
\end{gather}
\end{Lemma}

\Proof
Using (\ref{eq:orthogonality}) and (\ref{eq:X-derivatives}), and the discrete symmetry (\ref{eq:discrete-symmetry}), the area of $\Sigma_{n,a}$ is
\begin{equation}
\label{eq:Area1}
\mathcal{A}(X_{n,a})=\int_{Q_{n\tau_{a}}}|(X_{n,a})_{t}\wedge(X_{n,a})_{\theta}|~\!dt\!~d\theta=n\int_{Q_{\tau_{a}}}x_{a}^{2}\sqrt{1+\eps_{n}\psi_{n,a}}~\!dt\!~d\theta
\end{equation}
where
\begin{equation}
\label{eq:eps-psi}
\eps_{n}=\frac{\pi}{nh_{a}}\quad\text{and}\quad\psi_{n,a}=2x_{a}w_{a}^{2}\sin\theta+\eps_{n}x_{a}^{2}w_{a}^{2}\sin^{2}\theta~\!.
\end{equation}
Then we write
$$
\sqrt{1+\eps_{n}\psi_{n,a}}=1+\frac{\eps_{n}}{2}\psi_{n,a}-\frac{\eps_{n}^{2}}{2}g_{n,a}\quad\text{where}\quad
g_{n,a}=\frac{\psi_{n,a}^{2}}{2+\eps_{n}\psi_{n,a}+2\sqrt{1+\eps_{n}\psi_{n,a}}}~\!.
$$
Thus (\ref{eq:Area1}) becomes
\begin{equation}
\label{eq:Area2}
\mathcal{A}(X_{n,a})=n\left[2\pi\int_{-\tau_{a}}^{\tau_{a}}x_{a}^{2}~\!dt+\frac{4\pi}{n^{2}}{G}_{n}(a)\right]\quad\text{where}\quad{G}_{n}(a)=-\frac{\pi}{8h_{a}^{2}}\int_{Q_{\tau_{a}}}x_{a}^{2}g_{n,a}~\!dt\!~d\theta~\!.
\end{equation}
Observe that the area of a complete single period $\Sigma_{a}$ of the unduloid with neck-size $a$ is given by
\begin{equation}
\label{eq:Area3}
\text{area}(\Sigma_{a})=\int_{Q_{\tau_{a}}}|(X_{a})_{t}\wedge(X_{a})_{\theta}|~\!dt\!~d\theta=2\pi\int_{-\tau_{a}}^{\tau_{a}}x_{a}^{2}~\!dt~\!.
\end{equation}
As computed in \cite{HMO}, the area of $\Sigma_{a}$ can be expressed in terms of the complete elliptic integrals of second kind introduced in (\ref{eq:elliptic-integrals}), according to the formula
$$
\text{area}(\Sigma_{a})=4\pi(1-a)E(k_{a})
$$
with $k_{a}$ as in (\ref{eq:taua-ha}). Then, by (\ref{eq:taua-ha}),
\begin{equation}
\label{eq:Area4}
\text{area}(\Sigma_{a})=4\pi(1-a)h_{a}~\!.
\end{equation}
Hence from (\ref{eq:Area2})--(\ref{eq:Area4}), taking into account of the expansion (\ref{eq:KE-expansion}), (\ref{eq:area-expansion}) follows. It remains to show that ${G}_{n}(a)$ is bounded, and differentiable, with bounded derivative in a right neighborhood of $0$. This is true for the mapping $a\mapsto h_{a}^{-1}$ (see Lemma \ref{L:tau-a-expansion}). We have to check it for the integral
$$
\overline{G}_{n}(a)=\int_{Q_{\tau_{a}}}x_{a}^{2}g_{n,a}~\!dt\!~d\theta~\!.
$$
Observe that $|\psi_{n,a}|\le C$ and also $0<g_{n,a}\le C$ for some constant $C$ independent of $n$, and $a$, if $n$ is large enough. Then $0<\overline{G}_{n}(a)\le C\int_{Q_{\tau_{a}}}x_{a}^{2}~\!dt\!~d\theta=C~\!\text{area}(\Sigma_{a})\le C$. Now we check that $\overline{G}_{n}(a)$ is differentiable with bounded derivative uniformly with respect to $n$ large. In general, given a family of regular mappings $f_{a}\colon\R/_{2\tau_{a}}\times\R/_{2\pi}\to\R$ with a $C^{1}$ dependence on the parameter $a\in(0,a_{0})$, one can check that the integral function $a\mapsto\int_{Q_{\tau_{a}}}f_{a}~\!dt~\!d\theta$ is differentiable and
\begin{equation}
\label{eq:derivative-a-integral}
\frac{\partial}{\partial a}\int_{Q_{\tau_{a}}}f_{a}~\!dt~\!d\theta
=\int_{Q_{\tau_{a}}}\frac{\partial f_{a}}{\partial a}~\!dt~\!d\theta
+2\frac{\partial\tau_{a}}{\partial a}\int_{-\pi}^{\pi}f_{a}(\tau_{a},\theta)~\!d\theta~\!.
\end{equation}
We aim to apply (\ref{eq:derivative-a-integral}) with $f_{a}=x_{a}^{2}g_{n,a}$. One can write
\begin{equation}
\label{eq:a-derivative(xg)}
\frac{\partial}{\partial a}\left(x_{a}^{2}g_{n,a}\right)=x_{a}^{2}\left[g^{1}_{n,a}\frac{\partial}{\partial a}\left(x_{a}\psi_{n,a}\right)+g^{2}_{n,a}\frac{\partial}{\partial a}\left(\eps_{n}\psi_{n,a}\right)\right]
\end{equation}
for some regular and periodic mappings $g^{1}_{n,a}$ and $g^{2}_{n,a}$ such that $\|g^{i}_{n,a}\|_{C^{0}(Q_{\tau_{a}})}\le C$ ($i=1,2$) with $C$ independent of $n$ and $a$. Notice that also $\eps_{n}$ depends on $a$ (see its definition in (\ref{eq:eps-psi})) and by Lemma \ref{L:tau-a-expansion}, one has that $\big|\frac{\partial\eps_{n}}{\partial a}\big|\le C$ uniformly in $a$. In order to estimate $\frac{\partial\psi_{n,a}}{\partial a}$ and $\frac{\partial}{\partial a}(x_{a}\psi_{n,a})$ one can use the definition of $\psi_{n,a}$ given in (\ref{eq:eps-psi}), apply the identity
$$
x'_{a}\frac{\partial w_{a}}{\partial a}=w_{a,0}^{+}+\left(1-\frac{x_{a}'}{x_{a}}\right)w_{a}\frac{\partial x_{a}}{\partial a}
$$
(coming from (\ref{eq:wa0+}) and (\ref{eq:useful-notation})) and take advantage of the estimate (\ref{eq:wa-bound}). After some computations, one finds that $\big\|\frac{\partial\psi_{n,a}}{\partial a}\big\|_{C^{0}(Q_{\tau_{a}})}$ and $\big\|\frac{\partial}{\partial a}(x_{a}\psi_{n,a})\big\|_{C^{0}(Q_{\tau_{a}})}$ are bounded uniformly in $a$ and $n$. Hence, from (\ref{eq:a-derivative(xg)}) it follows that
$$
\left|\frac{\partial}{\partial a}\left(x_{a}^{2}g_{n,a}\right)\right|\le Cx_{a}^{2}
$$
and then
$$
\left|\int_{Q_{\tau_{a}}}\frac{\partial}{\partial a}\left(x_{a}^{2}g_{n,a}\right)~\!dt\!~d\theta\right|\le C~\!.
$$
Moreover, since $0<g_{n,a}\le C$ and $x_{a}(\tau_{a})=a$, by Lemma \ref{L:tau-a-expansion}, we have that
$$
\left|\frac{\partial\tau_{a}}{\partial a}\int_{-\pi}^{\pi}x_{a}(\tau_{a})^{2}g_{n,a}(\tau_{a},\theta)~\!d\theta\right|\le C
$$
for some constant independent of $a\in(0,a_{0})$. Thus we proved that
$$
\left|\frac{\partial\overline{G}_{n}}{\partial a}(a)\right|\le C
$$
uniformly in $a\in(0,a_{0})$ and $n$ large enough. This completes the proof of (\ref{eq:area-expansion}). Now let us examine the volume integral. By (\ref{eq:derivatives}), (\ref{eq:VW1}) one has that
\begin{eqnarray}
\nonumber
\mathcal{V}(X_{n,a})&=&\frac{1}{3}\int_{Q_{n\tau_{a}}}X_{n,a}\cdot(X_{n,a})_{t}\wedge(X_{n,a})_{\theta}~\!dt\!~d\theta\\
\nonumber
&=&-\frac{1}{3}\int_{Q_{n\tau_{a}}}x_{a}^{2}z_{a}'~\!dt\!~d\theta-\frac{1}{3}\int_{Q_{n\tau_{a}}}x_{a}^{2}z_{a}'\sin^{2}\theta~\!dt\!~d\theta\\
\label{eq:V1}
&=&-n\pi\int_{-\tau_{a}}^{\tau_{a}}x_{a}^{2}z_{a}'~\!dt
\end{eqnarray}
which equals the volume enclosed by a section of the unduloid made by $n$ moduli closed with two parallel disks at the ends. In fact
\begin{equation}
\label{eq:V2}
\mathrm{vol}(\Sigma_{a})=-\frac{1}{3}\int_{Q_{\tau_{a}}}X_{a}\cdot(X_{a})_{t}\wedge(X_{a})_{\theta}~\!dt\!~d\theta=\frac{2\pi}{3}\int_{-\tau_{a}}^{\tau_{a}}(x_{a}^{2}z_{a}'-x_{a}x_{a}'z_{a})~\!dt=\pi\int_{-\tau_{a}}^{\tau_{a}}x_{a}^{2}z_{a}'~\!dt
\end{equation}
(use (\ref{eq:X}) and an integration by parts). Also $\text{vol}(\Sigma_{a})$ can be expressed in terms of the complete elliptic integrals (\ref{eq:elliptic-integrals}). More precisely, as proved in \cite{HMO},
$$
\text{vol}(\Sigma_{a})=\frac{2\pi}{3}(1-a)\left((2-a+a^{2})E(k_{a})-a^{2}K(k_{a})\right)
$$
with $k_{a}$ as in (\ref{eq:taua-ha}). Then, by (\ref{eq:taua-ha}),
\begin{equation}
\label{eq:V3}
\text{vol}(\Sigma_{a})=\frac{2\pi}{3}(1-a)\left((2-a+a^{2})h_{a}-\frac{a^{2}\tau_{a}}{2}\right)~\!.
\end{equation}
Finally, using (\ref{eq:KE-expansion}), from (\ref{eq:V1})--(\ref{eq:V3}) one obtains (\ref{eq:volume-expansion}).
\qed

Fixing $H\in C^{1}(\R^{3})$, the $H$-energy of a toroidal surface with parameter domain $Q_{\tau}$ and parameterization $X$ is defined by
\begin{equation}
\label{eq:Henergy}
\mathcal{E}_{H}(X;Q_{\tau})=\mathcal{E}_{H}(X):=\int_{Q_{\tau}}|X_{t}\wedge X_{\theta}|~\!dt~\!d\theta+2\int_{Q_{\tau}}\mathcal{H}(X)\cdot X_{t}\wedge X_{\theta}~\!dt~\!d\theta
\end{equation}
where $\mathcal{H}\colon\R^{3}\to\R^{3}$ is a vector field such that $\mathrm{div}~\!\mathcal{H}=H$.

\begin{Remark}
\label{R:mH}
For $H\equiv 0$ the $H$-energy reduces to the area. For $H\equiv 1$ one can take $\mathcal{H}(X)=\frac{1}{3}X$. Then $\mathcal{E}_{1}(X)=\mathcal{A}(X)+2\mathcal{V}(X)$. For $H\colon\R^{3}\setminus\{0\}\to\R$ sufficiently regular, a natural choice of the vector field $\mathcal{H}$ is $\mathcal{H}(X)=m_{H}(X)X$, with $m_{H}(X)=\int_{0}^{1}H(sX)s^{2}~\!ds$.
\end{Remark}

Now let us discuss some facts about differentiability of the $H$-energy. A variation of $X$ is a mapping $Y\colon\R/_{2\tau}\times\R/_{2\pi}\times(-\varepsilon_{0},\varepsilon_{0})\to\R^{3}$ of class $C^{2}$ such that $Y(\zeta,0)=X(\zeta)$, where $\zeta=(t,\theta)$. The vector field
$$
Z(\zeta)=\frac{\partial}{\partial\varepsilon}Y(\zeta,\varepsilon)\Big|_{\varepsilon=0}
$$
is called the first variation of the family of surfaces $Y(\cdot,\varepsilon)$. The first variation of the $H$-energy $\mathcal{E}_{H}$ at $X$ along the variation $Y(\cdot,\varepsilon)$ is defined as
$$
\mathcal{E}_{H}(X;Q_{\tau})[Z]=\mathcal{E}_{H}(X)[Z]:=\frac{\partial}{\partial\varepsilon}\mathcal{E}_{H}(Y(\cdot,\varepsilon);Q_{\tau})\Big|_{\varepsilon=0}~\!.
$$

%Let us recall the following result:

%\begin{Lemma}
% Let $X$ be a parameterization of a toroidal surface with parameter domain $Q_{\tau}$ and let $Y(\cdot,\varepsilon)$ be a variation of $X$. Then the first variation of the area functional $\mathcal{A}(X)$ at $X$ along the variation $Y(\cdot,\varepsilon)$ is given by
% $$\frac{d}{d\varepsilon}\mathcal{A}(Y(\cdot,\varepsilon))\Big|_{\varepsilon=0}=-2\int_{Q_{\tau}}\mathfrak{M}(X)Z\cdot X_{t}\wedge X_{\theta}~\!dt~\!d\theta$$
%where $Z$ is the first variation of the family of surfaces $Y(\cdot,\varepsilon)$ and $\mathcal{M}(X)$ equals the mean curvature of the surface $X$.
%\end{Lemma}
%For the proof, see, e.g., [Dierkes, Hildebrandt, Sauvigny: Minimal Surfaces, pp. 54--56].
%\medskip

\begin{Lemma}
\label{L:first-variation}
Let $X$ be a parameterization of a toroidal surface with parameter domain $Q_{\tau}$ and let $Y(\cdot,\varepsilon)$ be a variation of $X$. Then the first variation of the $H$-energy $\mathcal{E}_{H}$ at $X$ along the variation $Y(\cdot,\varepsilon)$ is given by
$$
\mathcal{E}'_{H}(X;Q_{\tau})[Z]=\mathcal{E}'_{H}(X)[Z]=2\int_{Q_{\tau}}[H(X)-\mathfrak{M}(X)]Z\cdot X_{t}\wedge X_{\theta}~\!dt~\!d\theta
$$
where $Z$ is the first variation of the family of surfaces $Y(\cdot,\varepsilon)$ and $\mathfrak{M}(X)$ equals the mean curvature of the surface $X$.
\end{Lemma}

\Proof
The formula of the first variation of the area term is proved in \cite{DHS}, pp. 54--56. For the computation of the first variation of the remainder term one can plainly adapt an argument already used in \cite{CM2004}.
\qed

Let us apply Lemma \ref{L:first-variation} in order to obtain that:
\begin{Lemma}
\label{L:a-derivative-E1}
One has that
$$
\int_{Q_{\tau_{a}}}\left[\mathfrak{M}(X_{n,a})-1\right]x_{a}^{2}w_{a,0}\left(1+\eps_{n} x_{a}\sin\theta\right)~\!dt~\!d\theta=-\frac{1}{2n}\frac{\partial}{\partial a}\left[\mathcal{E}_{1}(X_{n,a})\right]~\!.
$$
\end{Lemma}

\Proof
Set $\hat X_{n,a}(s,\theta)=X_{n,a}(s\tau_{a},\theta)$ in order to have a common parameter domain. Thus one obtains that $\mathcal{E}_{H}(X_{n,a};Q_{n\tau_{a}})=\mathcal{E}_{H}(\hat X_{n,a};Q_{n})$ and, by Lemma \ref{L:first-variation},
\begin{equation*}
\begin{split}
\frac{\partial}{\partial a}\left[\mathcal{E}_{1}(X_{n,a})\right]&=\mathcal{E}'_{1}(\hat X_{n,a};Q_{n})\left[\frac{\partial \hat X_{n,a}}{\partial a}\right]\\
&=2\int_{Q_{n}}\big[1-\mathfrak{M}(\hat X_{n,a})\big]\left[\frac{\partial X_{n,a}}{\partial a}(s\tau_{a},\theta)+s\frac{\partial\tau_{a}}{\partial a}\frac{\partial X_{n,a}}{\partial t}(s\tau_{a},\theta)\right]\cdot(\hat X_{n,a})_{t}\wedge(\hat X_{n,a})_{\theta}~\!ds~\!d\theta\\
&=2\int_{Q_{n\tau_{a}}}\left[1-\mathfrak{M}(X_{n,a})\right]\frac{\partial X_{n,a}}{\partial a}\cdot(X_{n,a})_{t}\wedge(X_{n,a})_{\theta}~\!dt~\!d\theta\\
&=2n\int_{Q_{\tau_{a}}}\left[1-\mathfrak{M}(X_{n,a})\right]\frac{\partial X_{n,a}}{\partial a}\cdot(X_{n,a})_{t}\wedge(X_{n,a})_{\theta}~\!dt~\!d\theta
\end{split}
\end{equation*}
thanks to the discrete symmetry (\ref{eq:discrete-symmetry}). Finally the conclusion follows from the identity
$$
\frac{\partial X_{n,a}}{\partial a}\cdot (X_{n,a})_{t}\wedge (X_{n,a})_{\theta}=x_{a}^{2}w_{a,0}\left(1+\eps_{n} x_{a}\sin\theta\right)~\!.\quad\square
$$

A key result in the proof of Theorem \ref{T:immersed} is the following result.

\begin{Lemma}
\label{L:lambda0=0}
For $n$ large enough there exist continuous mappings $f_{n},g_{n}\colon(0,a_{0})\to\R$ which are bounded in $(0,a_{0})$, uniformly in $n$, such that $\lambda^{0}_{n,a}=0$ if and only if
\begin{equation}
\label{eq:a(n)}
a\log a-A\Big(\frac{\pi}{nh_{a}}\Big)^{\gamma}+a~\!f_{n}(a)+\frac{g_{n}(a)}{n^{\gamma+\nu}}=0
\end{equation}
for some $\nu>0$.
\end{Lemma}

\Proof
By Theorem  \ref{T:reduction}, one has
\begin{equation}
\label{eq:error}
\mathfrak{M}(X_{n,a}+\varphi_{n,a}N_{n,a})-{H}(X_{n,a}+\varphi_{n,a}N_{n,a})=\frac{1}{2x_{a}^{2}}\left(\lambda^{0}_{n,a}w_{a,0}+\lambda^{1}_{n,a}w_{a,1}\right)~\!.
\end{equation}
Multiplying (\ref{eq:error}) by $2x_{a}^{2}w_{a,0}$ and integrating over $Q_{\tau_{a}}$, we obtain
$$
\lambda^{0}_{n,a}\int_{Q_{\tau_{a}}}w_{a,0}^{2}~\!dt~\!d\theta=2\int_{Q_{\tau_{a}}}\left[\mathfrak{M}(X_{n,a}+\varphi_{n,a}N_{n,a})-{H}(X_{n,a}+\varphi_{n,a}N_{n,a})\right]x_{a}^{2}w_{a,0}~\!dt~\!d\theta~\!.
$$
Hence $\lambda^{0}_{n,a}=0$ if and only if
\begin{equation}
\label{eq:Mna=Hna}
\underbrace{\int_{Q_{\tau_{a}}}\left[\mathfrak{M}(X_{n,a}+\varphi_{n,a}N_{n,a})-1\right]x_{a}^{2}w_{a,0}~\!dt~\!d\theta}_{M_{n}(a)}=\underbrace{\int_{Q_{\tau_{a}}}\left[H(X_{n,a}+\varphi_{n,a}N_{n,a})-1\right]x_{a}^{2}w_{a,0}~\!dt~\!d\theta}_{H_{n}(a)}~\!.
\end{equation}
\textbf{Estimate of $M_{n}(a)$.} Let us decompose
\begin{equation*}
\begin{split}
M_{n}&(a)=\underbrace{\int_{Q_{\tau_{a}}}\left[\mathfrak{M}(X_{n,a}+\varphi_{n,a}N_{n,a})-\mathfrak{M}(X_{n,a})\right]x_{a}^{2}w_{a,0}~\!dt~\!d\theta}_{M^{1}_{n}(a)}\\
&+\underbrace{\int_{Q_{\tau_{a}}}\left[\mathfrak{M}(X_{n,a})-1\right]x_{a}^{2}w_{a,0}\left(1+\eps_{n} x_{a}\sin\theta\right)~\!dt~\!d\theta}_{M^{2}_{n}(a)}-\eps_{n}\underbrace{\int_{Q_{\tau_{a}}}\left[\mathfrak{M}(X_{n,a})-1\right]x_{a}^{3}w_{a,0}\sin\theta~\!dt~\!d\theta}_{M^{3}_{n}(a)}
\end{split}
\end{equation*}
Using Lemmata \ref{L:area-volume-expansions} and \ref{L:a-derivative-E1}, we obtain that
$$
M^{2}_{n}(a)=2\pi\left[a\log a+a~\!\Phi(a)+n^{-2}\Phi_{n}(a)\right]
$$
where $\Phi$ and $\Phi_{n}$ are continuous mappings in $a$, bounded in $(0,a_{0})$, uniformly in $n$. By Lemma \ref{L:M(X)-eps}, we get
$$
M^{3}_{n}(a)=\eps_{n}\overline{\Phi}_{n}(a)
$$
where $\overline{\Phi}_{n}$ are continuous mappings in $a$, bounded in $(0,a_{0})$, uniformly in $n$. In order to evaluate $M_{n}^{1}(a)$ we use (\ref{eq:MM-decomposition}) as well as some estimates proved in Section \ref{S:reduced}. In particular, we have that
$$
M_{n}^{1}(a)=\frac{1}{2}\int_{Q_{\tau_{a}}}(\mathfrak{L}_{a}\varphi_{n,a})w_{a,0}~\!dt~\!ds+\frac{1}{2}\int_{Q_{\tau_{a}}}[(\mathfrak{L}_{\eps_{n},a}-\mathfrak{L}_{a})\varphi_{n,a}]w_{a,0}~\!dt~\!ds+\frac{1}{2}\int_{Q_{\tau_{a}}}\mathfrak{B}^{0}_{\eps_{n},a}(\varphi_{n,a})w_{a,0}~\!dt~\!ds~\!.
$$
Integrating by parts twice and using the fact that $w_{a,0}\in\mathrm{ker}(\mathfrak{L}_{a})$ we obtain
$$
\int_{Q_{\tau_{a}}}(\mathfrak{L}_{a}\varphi_{n,a})w_{a,0}~\!dt~\!ds=\int_{Q_{\tau_{a}}}\varphi_{n,a}(\mathfrak{L}_{a}w_{a,0})~\!dt~\!ds=0~\!.
$$
By (\ref{eq:estimate-2}), (\ref{eq:varphi-small}) and (\ref{eq:weighted-norms}), we estimate
$$
\left|\int_{Q_{\tau_{a}}}[(\mathfrak{L}_{\eps_{n},a}-\mathfrak{L}_{a})\varphi_{n,a}]w_{a,0}~\!dt~\!ds\right|\le C\eps_{n}^{\widetilde\gamma+1}\int_{Q_{\tau_{a}}}x_{a}^{\mu}|w_{a,0}|~\!dt~\!d\theta~\!.
$$
Moreover, by (\ref{eq:estimate-4}), (\ref{eq:varphi-small}) and (\ref{eq:weighted-norms}), we get
$$
\left|\int_{Q_{\tau_{a}}}\mathfrak{B}^{0}_{\eps_{n},a}(\varphi_{n,a})w_{a,0}~\!dt~\!ds\right|\le C\eps_{n}^{2\widetilde\gamma}\int_{Q_{\tau_{a}}}x_{a}^{\mu}|w_{a,0}|~\!dt~\!d\theta~\!.
$$
By Lemmata \ref{L:xa-limit} and \ref{L:xwp}, there exists a constant $C$ such that
\begin{equation*}
%\label{eq:4pi-bis}
\int_{Q_{\tau_{a}}}x_{a}^{\mu}|w_{a,0}|~\!dt~\!d\theta=2\pi\int_{-\tau_{a}}^{\tau_{a}}x_{a}^{\mu}|w_{a,0}^{+}|~\!dt\le C\quad\forall a\in(0,a_{0})~\!.
\end{equation*}
Then
\begin{equation}
\label{eq:Mna-estimate}
M_{n}(a)=2\pi\left[a\log a+a~\!\Phi(a)+n^{-2}\Phi_{n}(a)+n^{-2\widetilde\gamma}\widetilde{\Phi}_{n}(a)\right]
\end{equation}
where $\Phi$ and $\widetilde{\Phi}_{n}$ are continuous mappings in $a$, bounded in $(0,a_{0})$, uniformly in $n$.
\medskip

\noindent
\textbf{Estimate of $H_{n}(a)$.} By $(H_{1})$ we can write
$$
H(X)=1+A|X|^{-\gamma}+|X|^{-\gamma-\beta}\widetilde{H}(X)
$$
with $\widetilde{H}$ continuous and bounded function on $\R^{3}$.
Setting
$$
X_{n,a}^{0}=\eps_{n}^{-1}\left[\begin{array}{c}0\\ \cos(\eps_{n}z_{a})\\ \sin(\eps_{n}z_{a})\end{array}\right]~\!,
$$
we write the following decomposition
\begin{equation*}
\begin{split}
H_{n}(a)&=A\underbrace{\int_{Q_{\tau_{a}}}\frac{x_{a}^{2}w_{a,0}}{|X_{n,a}^{0}|^{\gamma}}~\!dt~\!d\theta}_{H^{1}_{n}(a)}+A\underbrace{\int_{Q_{\tau_{a}}}\left[|X_{n,a}+\varphi_{n,a}N_{n,a}|^{-\gamma}-|X_{n,a}^{0}|^{-\gamma}\right]x_{a}^{2}w_{a,0}~\!dt~\!d\theta}_{H^{2}_{n}(a)}\\
&\quad+\underbrace{\int_{Q_{\tau_{a}}}
\frac{\widetilde{H}(X_{n,a}+\varphi_{n,a}N_{n,a})}{|X_{n,a}+\varphi_{n,a}N_{n,a}|^{\gamma+\beta}}x_{a}^{2}w_{a,0}~\!dt~\!d\theta}_{H^{3}_{n}(a)}~\!.
\end{split}
\end{equation*}
The next auxiliary result holds true:
\begin{Lemma}
\label{L:4pi}
There exists a continuous and bounded function $\Psi\colon(0,a_{0})\to\R$ such that
\begin{equation}
\label{eq:4pi}
\int_{Q_{\tau_{a}}}x_{a}^{2}w_{a,0}~\!dt~\!d\theta
=2\pi+a~\!\Psi(a)~\!.
\end{equation}
Moreover, there exists a constant $C>0$ such that
\begin{equation}
\label{eq:4pi-ter}
\int_{Q_{\tau_{a}}}x_{a}^{2}|w_{a,0}|~\!dt~\!d\theta\le C\quad\forall a\in(0,a_{0})~\!.
\end{equation}
\end{Lemma}

\noindent
Postponing the proof of Lemma \ref{L:4pi}, let us complete that one of Lemma \ref{L:lambda0=0}. Since $|X_{n,a}^{0}|=\eps_{n}^{-1}$, by (\ref{eq:4pi}),
$$
H^{1}_{n}(a)=2\pi\eps_{n}^{-\gamma}(1+a~\!\Psi(a))~\!.
$$
Since $|X_{n,a}-X_{n,a}^{0}|=x_{a}<1$, with elementary arguments based on the mean value theorem (applied to the mapping $r\mapsto |X_{n,a}^{0}+r(X_{n,a}-X_{n,a}^{0}+\varphi_{n,a}N_{n,a})|^{-\gamma}$, with $r\in[0,1]$) we can estimate
$$
\big||X_{n,a}+\varphi_{n,a}N_{n,a}|^{-\gamma}-|X_{n,a}^{0}|^{-\gamma}\big|\le C\eps_{n}^{\gamma+1}
$$
and then
$$
\left|H^{2}_{n}(a)\right|\le C\eps_{n}^{\gamma+1}\int_{Q_{\tau_{a}}}x_{a}^{2}|w_{a,0}|~\!dt~\!d\theta\le C\eps_{n}^{\gamma+1}
$$
by (\ref{eq:4pi-ter}). Moreover, again by (\ref{eq:4pi-ter}), we have that $\left|H^{3}_{n}(a)\right|\le C\eps_{n}^{\gamma+\beta}$. Hence we obtain
\begin{equation}
\label{eq:Hna-estimate}
H_{n}(a)=2\pi \eps_{n}^{\gamma}(A+a~\!\Psi(a))+\eps_{n}^{\gamma+\min\{1,\beta\}}\Psi_{n}(a)
\end{equation}
where $\Psi$ and $\Psi_{n}$ are continuous mappings in $a$, bounded in $(0,a_{0})$, uniformly in $n$. From (\ref{eq:Mna=Hna}), (\ref{eq:Mna-estimate}) and (\ref{eq:Hna-estimate}), since $2\widetilde{\gamma}>\gamma$ for $\gamma\in(0,2)$, (\ref{eq:a(n)}) follows with $\nu=\min\{1,\beta\}$.
\qed

\noindent
\emph{Proof of Lemma \ref{L:4pi}.}
Using the expression (\ref{eq:wa0+}) of $w_{a,0}=w_{a,0}^{+}$, by (\ref{eq:derivative-a-integral}) and an integration by parts, one finds
$$
2\int_{-\tau_{a}}^{\tau_{a}}x_{a}^{2}w_{a,0}~\!dt=x_{a}(\tau_{a})^{2}\left[\frac{\partial z_{a}}{\partial a}(\tau_{a})-\frac{\partial z_{a}}{\partial a}(-\tau_{a})\right]-\frac{\partial}{\partial a}\int_{-\tau_{a}}^{\tau_{a}}x_{a}^{2}z'_{a}~\!dt~\!.
$$
Since $z_{a}$ is odd and $z'_{a}=\gamma_{a}+x_{a}^{2}$ (see Lemma \ref{L:c-a}), one computes
\begin{equation*}\begin{split}
\frac{\partial z_{a}}{\partial a}(\tau_{a})-\frac{\partial z_{a}}{\partial a}(-\tau_{a})&=2\frac{\partial z_{a}}{\partial a}(\tau_{a})=2\frac{\partial}{\partial a}\int_{0}^{\tau_{a}}(\gamma_{a}+x_{a}^{2})~\!dt-2(\gamma_{a}+x_{a}(\tau_{a})^{2})\frac{\partial\tau_{a}}{\partial a}\\
&=(1-2a)2\tau_{a}-2a^{2}\frac{\partial\tau_{a}}{\partial a}+\frac{\partial}{\partial a}\int_{-\tau_{a}}^{\tau_{a}}x_{a}^{2}~\!dt
\end{split}\end{equation*}
because $\gamma_{a}=a(1-a)$ and $x_{a}(\tau_{a})=a$. Therefore
$$
\int_{Q_{\tau_{a}}}x_{a}^{2}w_{a,0}~\!dt~\!d\theta=
2\pi a^{2}(1-2a)\tau_{a}-2\pi a^{4}\frac{\partial\tau_{a}}{\partial a}+\frac{a^{2}}{2}\frac{\partial}{\partial a}\int_{Q_{\tau_{a}}}x_{a}^{2}~\!dt~\!d\theta-\frac{1}{2}\frac{\partial}{\partial a}\int_{Q_{\tau_{a}}}x_{a}^{2}z'_{a}~\!dt~\!d\theta~\!.
$$
By (\ref{eq:Area3})--(\ref{eq:Area4}), (\ref{eq:V2})--(\ref{eq:V3}), and with the aid of Lemma \ref{L:tau-a-expansion}, we can estimate the last two terms in the previous equation and we arrive to (\ref{eq:4pi}). Finally (\ref{eq:4pi-ter}) follows from Lemmata \ref{L:xa-limit} and \ref{L:xwp}.
\qed

\noindent
\emph{Proof of Theorem \ref{T:immersed}.} We split the proof in two steps.
\medskip

\noindent
\textbf{Step 1.} There exist $n_{A}\in\mathbb{N}$, with $n_{A}\to\infty$ as $A\to 0$, and a sequence $(a_{n})_{n\ge n_{A}}\subset (0,a_{0})$ satisfying (\ref{an-behaviour}) and such that $\lambda^{0}_{n,a_{n}}=0$ for every $n\ge n_{A}$.
\medskip

\noindent
According to Lemma \ref{L:lambda0=0}, one has to solve equation
$$
F_{n}(a)=0\quad\text{where~~}F_{n}(a)=a\log a-A\Big(\frac{\pi}{nh_{a}}\Big)^{\gamma}+a~\!f_{n}(a)+\frac{g_{n}(a)}{n^{\gamma+\nu}}
$$
and $f_{n},g_{n}$ are continuous functions in $(0,a_{0})$ uniformly bounded in $n$. Setting
$$
a=\frac{b}{n^{\gamma}\log n}
$$
one has that
$$
a\in(0,a_{0})~\Leftrightarrow~b\in(0,\overline{b}_{n})\quad\text{where~~}\overline{b}_{n}:=a_{0}n^{\gamma}\log n~\!.
$$
Moreover, setting
$$
\vartheta_{n}(b):=\frac{b\log b-b\log(\log n)}{\log n}+\frac{b}{\log n}~\!f_{n}\left(\frac{b}{n^{\gamma}\log n}\right)+\frac{1}{n^{\nu}}g_{n}\left(\frac{b}{n^{\gamma}\log n}\right)~\!,\quad
k_{n}(b):=\left(\frac{\pi}{h_{b/(n^{\gamma}\log n)}}\right)^{\gamma}~\!,
$$
one has that
$$
\left\{\begin{array}{l}F_{n}(a)=0\\
a\in(0,a_{0})\end{array}\right.~\Leftrightarrow~
\left\{\begin{array}{l}
\widetilde{F}_{n}(b)=0\\
b\in(0,\overline{b}_{n}).\end{array}\right.\quad\text{where}\quad \widetilde{F}_{n}(b):=-\gamma b-Ak_{n}(b)+\vartheta_{n}(b)~\!.
$$
Observe that $\vartheta_{n},k_{n}\colon(0,\overline{b}_{n})\to\R$ are continuous functions. Moreover there exist $n_{0}\in\mathbb{N}$, $k_{0},k_{1}>0$ such that
$$
k_{1}\le k_{n}(b)\le k_{0}\quad\forall n\ge n_{0}~\!,\quad\forall b\in(0,\overline{b}_{n})~\!,\quad\text{and}\quad
k_{n}(b)\to k_{0}\quad\text{as~~}b\to 0~\!.
$$
In particular,
$$
-\gamma b-Ak_{n}(b)\le-\gamma b+|A|k_{0}\quad\forall n\ge n_{0}~\!,\quad\forall b\in(0,\overline{b}_{n})~\!.
$$
Fix $b_{0}>0$ such that $-\gamma b_{0}+|A|k_{0}=-\delta_{0}<0$. Since $\lim_{n\to\infty}\vartheta_{n}(b)=0$ for every $b>0$, there exists $n_{1}\ge n_{0}$ such that $\vartheta_{n}(b_{0})<\delta_{0}$. Hence
\begin{equation*}
%\label{eq:negativo}
\widetilde{F}_{n}(b_{0})<0\quad\forall n\ge n_{1}~\!.
\end{equation*}
Letting
$$
C:=\sup_{\scriptstyle n\ge n_{1}\atop\scriptstyle a\in(0,a_{0})}|g_{n}(a)|
$$
since $C<\infty$ and $A<0$, we can find $n_{A}\ge n_{1}$ such that
\begin{equation}
\label{eq:nA}
-Ak_{0}-\frac{C}{n^{\nu}}>0\quad\forall n\ge n_{A}~\!.
\end{equation}
Then
\begin{equation*}
%\label{eq:positivo}
\liminf_{b\to 0}\widetilde{F}_{n}(b)>0\quad\forall n\ge n_{A}~\!.
\end{equation*}
As $\widetilde{F}_{n}$ is continuous in $(0,b_{0}]$, for every $n\ge n_{A}$ there exists $b_{n}\in(0,b_{0})$ such that $\widetilde{F}_{n}(b_{n})=0$. Then $\lambda_{n,a_{n}}^{0}=0$ for $a_{n}=\frac{b_{n}}{n^{\gamma}\log n}$ and $n\ge n_{A}$. Notice that $n_{A}\to\infty$ as $A\to 0$, by (\ref{eq:nA}). Moreover $\liminf b_{n}>0$, because otherwise, for a subsequence, $\widetilde{F}_{n}(b_{n})\to-Ak_{0}\ne 0$.
\medskip

\noindent
\textbf{Step 2.} One has that $\lambda^{1}_{n,a_{n}}=0$ for every $n\ge n_{A}$.
\medskip

\noindent
Set $X=X_{\eps_{n},a_{n}}$, $N={\normal}_{\eps_{n},a_{n}}$, $\varphi=\varphi_{\eps_{n},a_{n}}$ and $Y=X+\varphi N$. Denote by $R_{\sigma}$ the rotation of an angle $\sigma$ about the $x_{1}$-axis, defined as in (\ref{eq:rotation-matrix}). Since $H$ is radially symmetric, the $H$-energy of $Y$, defined by (\ref{eq:Henergy}), does not change under any rotation $R_{\sigma}$. Hence, in view of Lemma \ref{L:first-variation}, one has
\begin{equation}
\label{eq:L1-1}
0=\frac{d}{d\sigma}\left[\mathcal{E}_{H}(R_{\sigma}Y)\right]=\mathcal{E}'_{H}(R_{\sigma}Y)\left[\frac{ dR_{\sigma}}{d\sigma}Y\right]=2\int_{Q_{n\tau_{a}}}E_{H}(R_{\sigma}Y)\left[\frac{ dR_{\sigma}}{d\sigma}Y\right]\cdot (R_{\sigma}Y)_{t}\wedge (R_{\sigma}Y)_{\theta}~\!dt~\!d\theta
\end{equation}
where $E_{H}(R_{\sigma}Y)=\mathfrak{M}(R_{\sigma}Y)-{H}(R_{\sigma}Y)$. One can compute
$$
\frac{dR_{\sigma}}{d\sigma}=\left[\begin{array}{ccc}0&0&0\\ 0&-\sin\sigma&-\cos\sigma\\ 0&\cos\sigma&-\sin\sigma\end{array}\right]
$$
and
\begin{equation}
\label{eq:L1-2}
\left[\frac{ dR_{\sigma}}{d\sigma}Y\right]\cdot (R_{\sigma}Y)_{t}\wedge (R_{\sigma}Y)_{\theta}=(\mathbf{e}_{1}\wedge Y)\cdot Y_{t}\wedge Y_{\theta}=\eps_{n}^{-1}x_{a_{n}}x'_{a_{n}}+
x^{2}_{a_{n}}x'_{a_{n}}\sin\theta+\rho_{n}
\end{equation}
where $\rho_{n}$ is a polynomial of order 3 in $\varphi,\varphi_{t},\varphi_{\theta}$, with no term of order zero and such that\begin{itemize}
\item
the coefficients of the first order terms are:\\
$(e_{1}\wedge X)\cdot(X_{t}\wedge N_{\theta})=(e_{1}\cdot X)_{t}(X\cdot N)_{\theta}-\tfrac{1}{2}(e_{1}\cdot N)_{\theta}|X|^{2}_{t}$,
\\
$(e_{1}\wedge X)\cdot(N_{t}\wedge X_{\theta})=\tfrac{1}{2}(e_{1}\cdot N)_{t}|X|_{\theta}^{2}-(e_{1}\cdot X)_{\theta}(X\cdot N)_{t}$,
\\
$(e_{1}\wedge X)\cdot(X_{t}\wedge N)=(e_{1}\cdot X)_{t}(X\cdot N)-\tfrac{1}{2}(e_{1}\cdot N)|X|^{2}_{t}$,
\\
$(e_{1}\wedge X)\cdot(N\wedge X_{\theta})=\tfrac{1}{2}(e_{1}\cdot N)|X|_{\theta}^{2}-(e_{1}\cdot X)_{\theta}(X\cdot N)$;
\item
the coefficients of the second order terms are:\\
$(e_{1}\wedge X)\cdot(N\wedge N_{\theta})=(e_{1}\cdot N)(X\cdot N)_{\theta}-(e_{1}\cdot N)_{\theta}(X\cdot N)$,
\\
$(e_{1}\wedge X)\cdot(N_{t}\wedge N)=(e_{1}\cdot N)_{t}(X\cdot N)-(e_{1}\cdot N)(X\cdot N)_{t}$,
\\
$(e_{1}\wedge X)\cdot(N_{t}\wedge N_{\theta})=(e_{1}\cdot N)_{t}(X\cdot N)_{\theta}-(e_{1}\cdot N)_{\theta}(X\cdot N)_{t}$,
\\
$(e_{1}\wedge N)\cdot(N_{t}\wedge X_{\theta})=(e_{1}\cdot N)_{t}(N\cdot X)_{\theta}$,
\\
$(e_{1}\wedge N)\cdot(X_{t}\wedge N_{\theta})=-(e_{1}\cdot N)_{\theta}(N\cdot X)_{t}$;
\item
the coefficients of the third order terms are:\\
$(e_{1}\wedge N)\cdot(N_{t}\wedge N_{\theta})=0$,\\
$(e_{1}\wedge N)\cdot(N\wedge N_{\theta})=-(e_{1}\cdot N)_{\theta}$,\\
$\quad(e_{1}\wedge N)\cdot(N_{t}\wedge N)=(e_{1}\cdot N)_{t}$.
\end{itemize}
By the estimates stated in Lemma \ref{L:Xt}, thanks to the following computations
\begin{equation*}
\begin{split}
&e_{1}\cdot N=-x_{a_{n}}z_{a_{n}}'(1+\eps_{n} x_{a_{n}}\sin\theta)\cos\theta\\
&(e_{1}\cdot N)_{t}=-(x_{a_{n}}'z_{a_{n}}'+x_{a_{n}}z_{a_{n}}'')(1+\eps_{n} x_{a_{n}}\sin\theta)\cos\theta-\eps_{n} x_{a_{n}}x_{a_{n}}'z_{a_{n}}'\sin\theta\cos\theta\\
&(e_{1}\cdot N)_{\theta}=x_{a_{n}}z_{a_{n}}'(1+\eps_{n} x_{a_{n}}\sin\theta)\sin\theta-\eps_{n} x_{a_{n}}^{2}z_{a_{n}}'\cos^{2}\theta~\!,
\end{split}
\end{equation*}
we infer that all the coefficients of $g_{n}$ are bounded by $Cx_{a_{n}}^{2}$ for some constant $C$ independent of $n$. Then, in view of Remark \ref{R:point-wise-bound}, we have that
\begin{equation}
\label{eq:rhon}
|\rho_{n}|\le C\eps_{n}^{\widetilde\gamma}x_{a_{n}}^{2+\mu}\quad\text{on~~}\R^{2}~\!.
\end{equation}
Since, by Step 1, $E_{H}(Y)=\frac{1}{2x_{a_{n}}^{2}}\lambda^{1}_{n,a_{n}}w_{a_{n},1}$ and recalling (\ref{eq:eigenfunctions}) and (\ref{eq:useful-notation}), from (\ref{eq:L1-1})--(\ref{eq:L1-2}) it follows that
$$
0=\lambda^{1}_{n,a_{n}}\int_{Q_{n\tau_{a_{n}}}}
\frac{w_{a_{n}}(t+\nu_{n}\sigma)}{x_{a_{n}}^{2}(t+\nu_{n}\sigma)}
\sin\theta\left(\eps_{n}^{-1}(1+\eps_{n} x_{a_{n}}(t)\sin\theta)x_{a_{n}}(t)x'_{a_{n}}(t)+\rho_{n}(t,\theta)\right)~\!dt~\!d\theta\quad\forall\sigma\in\R~\!.
$$
Setting $s=\nu_{n}\sigma$, denoting $T_{s}f(t)=f(t+s)$, and exploiting the discrete symmetry (\ref{eq:discrete-symmetry}), we obtain
\begin{equation}
\label{eq:L1-3}
0=\lambda^{1}_{n,a_{n}}\left[\int_{Q_{\tau_{a_{n}}}} \frac{T_{s}w_{a_{n}}}{T_{s}x_{a_{n}}^{2}}x_{a_{n}}^{2}x'_{a_{n}}\sin^{2}\theta~\!dt~\!d\theta+\int_{Q_{\tau_{a_{n}}}} \frac{T_{s}w_{a_{n}}}{T_{s}x_{a_{n}}^{2}}\rho_{n}\sin\theta~\!dt~\!d\theta\right]\quad\forall s\in\R~\!.
\end{equation}
Arguing as in the proof of Lemma \ref{L:weight}, we can see that $x_{a}(t)\le e^{|s|}x_{a}(t+s)$ for every $t,s\in\R$ and then
\begin{equation}
\label{eq:xas}
\left|\frac{x_{a_{n}}(t)^{2}}{T_{s}x_{a_{n}}(t)^{2}}\right|\le e^{2|s|}\quad\forall t,s\in\R~\!.
\end{equation}
Moreover $0<w_{a}\le 1$ (Lemma \ref{L:xwp}) and $|x'_{a}|\le x_{a}\le\sqrt{\sech t}$ in $[-\tau_{a},\tau_{a}]$ (Lemma \ref{L:xa-limit}). Hence, using again Lemmata \ref{L:xa-limit} and \ref{L:xwp}, by the dominated convergence theorem we obtain
\begin{equation}
\label{eq:L1-4}
\begin{split}
\lim_{n\to\infty}\int_{Q_{\tau_{a_{n}}}} \frac{T_{s}w_{a_{n}}}{T_{s}x_{a_{n}}^{2}}x_{a_{n}}^{2}x'_{a_{n}}\sin^{2}\theta~\!dt~\!d\theta&=2\pi\lim_{n\to\infty}\int_{0}^{\tau_{a_{n}}} \frac{T_{s}w_{a_{n}}}{T_{s}x_{a_{n}}^{2}}x_{a_{n}}^{2}x'_{a_{n}}~\!dt\\
&=-\pi\int_{-\infty}^{\infty}\frac{(\sech t)^{3}\tanh t}{\sech(s+t)}~\!dt=-\frac{2\pi}{3}\sinh s~\!.
\end{split}
\end{equation}
On the other hand, by (\ref{eq:rhon}), (\ref{eq:xas}), and the previously used estimates on $w_{a}$ and $x_{a}$, we infer that
$$
\left|\frac{T_{s}w_{a_{n}}}{T_{s}x_{a_{n}}^{2}}\rho_{n}\sin\theta\right|\le C\eps_{n}^{\widetilde\gamma}e^{2|s|}(\sech t)^{\frac{\mu}{2}}\quad\text{in }Q_{\tau_{a_{n}}}
$$
and then
\begin{equation}
\label{eq:L1-5}
\lim_{n\to\infty}\int_{Q_{\tau_{a_{n}}}} \frac{T_{s}w_{a_{n}}}{T_{s}x_{a_{n}}^{2}}\rho_{n}\sin\theta~\!dt~\!d\theta=0~\!.
\end{equation}
Thus (\ref{eq:L1-3}), (\ref{eq:L1-4}) and (\ref{eq:L1-5}) imply that $\lambda^{1}_{n,a_{n}}=0$ for every $n\ge n_{A}$ (note that $s$ can be taken arbitrarily large).
\qed

\begin{Remark}
The proof of Theorem \ref{T:immersed} works even for every $\eps_{q}=\frac{\pi}{h_{a}}q$ with $q$ rational, sufficiently small. Hence we can show the existence of a countable family of immersed tori close to generalized toroidal unduloids $X_{\eps_{q},a_{q}}$ for suitable $a_{q}>0$ small enough.
\end{Remark}

\section{The embeddedness property}
\label{S:embeddedness}

In this section we complete the proof of Theorem \ref{T:main} showing that:

\begin{Theorem}
\label{T:embedded}
If $H\colon\R^{3}\setminus\{0\}\to\R$ is a $C^{2}$ function satisfying $(H_{1})$--$(H_{2})$, with $A<0$, then for $n$ large enough the maps $X_{n,a_{n}}+\varphi_{n,a_{n}}N_{n,a_{n}}$ found in Theorem \ref{T:immersed} are regular parametrization of a compact, oriented, embedded surface with genus~1.
\end{Theorem}

To show the validity of Theorem \ref{T:embedded} we introduce suitable tubular neighborhoods of the reference surfaces $\Sigma_{n,a}=X_{n,a}(\R^{2})$. Such neighborhoods are foliations whose leaves are surfaces parametrized by $X_{n,a}+r x_{a}N_{n,a}$, where $r$ is a real parameter small enough. Observe that, for fixed $r$, the normal variation, ruled by the factor $x_{a}$, is small in correspondence of the necks, whereas it can be larger at the bulges.
The next result states that we can take $r$ in a neighborhood of $0$ which is independent of the parameters $a$ and $n$.

\begin{Lemma}
There exists $r_{0}>0$ such that for every $r\in[-r_{0},r_{0}]$, for every $a\in(0,a_{0})$ and for every $n\ge n_{0}$, the map $X_{n,a}+r x_{a}N_{n,a}$ is a regular parametrization of a compact, oriented, embedded surface with genus~1.
\end{Lemma}

\Proof
Fix $a\in(0,a_{0})$ and $n\ge n_{0}$. For every $r\in\R$, the map $X_{n,a}+r x_{a}N_{n,a}\colon\R^{2}\to\R^{3}$ is smooth and doubly-periodic with respect to the rectangle $Q_{n\tau_{a}}$. Moreover we can write
$$
X_{n,a}+r x_{a}N_{n,a}=R_{\eps z_{a}}\left(\left[\begin{array}{c}0\\ \eps^{-1}\\ rx_{a}^{2}x'_{a}\end{array}\right]+x_{a}\left[\begin{array}{c}(1-rx_{a}z'_{a}-\eps rx_{a}^{2}z'_{a}\sin\theta)\cos\theta\\
(1-rx_{a}z'_{a}-\eps rx_{a}^{2}z'_{a}\sin\theta)\sin\theta\\ 0
\end{array}\right]\right)
$$
where, as in the previous Sections, $\eps=\frac{\pi}{nh_{a}}$ and $R_{\sigma}$ is the rotation matrix of an angle $\sigma$ about the $x_{1}$-axis, defined in (\ref{eq:rotation-matrix}). Therefore, for fixed $t$, the image of $\theta\mapsto X_{n,a}(t,\theta)+r x_{a}(t)N_{n,a}(t,\theta)$ is a closed curve which, up to a translation and a dilation, is defined by the parametric equations of the form
$$
\left\{\begin{array}{l}
x_{1}=(r_{1}+r_{2}\sin\theta)\cos\theta\\
x_{2}=(r_{1}+r_{2}\sin\theta)\sin\theta\\
x_{3}=0
\end{array}\right.
$$
where $r_{1}=1-rx_{a}(t)z'_{a}(t)$ and $r_{2}=-\eps rx_{a}(t)^{2}z'_{a}(t)$. Observe that the function $\theta\mapsto(1+\alpha\sin\theta)e^{i\theta}$ is a diffeomorphism between $\S^{1}$ and a star-shaped (with respect to the origin), hence simple, Jordan curve if (and only if) $|\alpha|<1$. Since $a\le x_{a}\le 1-a$, and $0<z'_{a}<x_{a}$, we can find $r_{0}>0$ independent of $a$ and $n$ such that if $|r|\le r_{0}$ then
$$
\left|\frac{r_{2}}{r_{1}}\right|=\frac{\pi}{nh_{a}}\left|\frac{rx_{a}(t)^{2}z'_{a}(t)}{1-rx_{a}(t)z'_{a}(t)}\right|<1\quad\forall t\in\R~\!.
$$
Hence the map $X_{n,a}+r x_{a}N_{n,a}$ is injective in $[-\pi,\pi)\times[-n\tau_{a},n\tau_{a})$. Now let us show that it is an immersion, namely
\begin{equation}
\label{eq:immersion}
(X_{n,a}+r x_{a}N_{n,a})_{t}\wedge(X_{n,a}+r x_{a}N_{n,a})_{\theta}\ne 0\quad\text{on $\R^{2}$.}
\end{equation}
For the sake of simplicity, let us write $X=X_{n,a}$, $N=N_{n,a}$ and $Y=X_{n,a}+r x_{a}N_{n,a}$, and let us compute
$$
N\cdot Y_{t}\wedge Y_{\theta}=N\cdot X_{t}\wedge X_{\theta}+rx_{a}N\cdot(X_{t}\wedge N_{\theta}+N_{t}\wedge X_{\theta})+r^{2}x_{a}^{2}N\cdot N_{t}\wedge N_{\theta}~\!.
$$
By (\ref{316B})--(\ref{316C}) and (\ref{eq:M-na}) we have
$$
N\cdot(X_{t}\wedge N_{\theta}+N_{t}\wedge X_{\theta})=-2|X_{t}|~\!|X_{\theta}|~\!\mathfrak{M}(X_{n,a})~\!.
$$
Moreover we can find a constant $C>0$ independent of $a\in(0,a_{0})$ and $n\ge n_{0}$ such that
\begin{equation}
\label{NNN}
|N\cdot N_{t}\wedge N_{\theta}|\le C\quad\text{on $\R^{2}$.}
\end{equation}
Hence, taking a smaller $r_{0}$, if necessary, but always independent of $a$ and $n$, we obtain that
$$
N\cdot Y_{t}\wedge Y_{\theta}\ge|X_{t}|~\!|X_{\theta}|~\!(1-2rx_{a}\mathfrak{M}(X_{n,a}))-Cr^{2}x_{a}^{2}\ge\frac{1}{2}x_{a}^{2}\left[1-2rx_{a}\mathfrak{M}(X_{n,a})-Cr^{2}\right]>0\quad\text{on $\R^{2}$}
$$
when $|r|\le r_{0}$. It remains to check (\ref{NNN}). One has that
\begin{equation*}
N\cdot N_{t}\wedge N_{\theta}=\frac{(N_{t}\cdot X_{t})(N_{\theta}\cdot X_{\theta})}{|X_{t}\wedge X_{\theta}|}-\frac{(N_{t}\cdot X_{\theta})(N_{\theta}\cdot X_{t})}{|X_{t}\wedge X_{\theta}|}=\frac{(N\cdot X_{tt})(N\cdot X_{\theta\theta})}{|X_{t}\wedge X_{\theta}|}-\frac{(N\cdot X_{t\theta})^{2}}{|X_{t}\wedge X_{\theta}|}
\end{equation*}
%\begin{equation*}
%\begin{split}
%N_{n}\cdot(N_{n})_{t}\wedge (N_{n})_{\theta}&=\frac{((N_{n})_{t}\cdot (X_{n})_{t})((N_{n})_{\theta}\cdot(X_{n})_{\theta})}{|(X_{n})_{t}\wedge (X_{n})_{\theta}|}-\frac{((N_{n})_{t}\cdot(X_{n})_{\theta})((N_{n})_{\theta}\cdot(X_{n})_{t})}{|(X_{n})_{t}\wedge(X_{n})_{\theta}|}\\&=\frac{(N_{n}\cdot (X_{n})_{tt})(N_{n}\cdot(X_{n})_{\theta\theta})}{|(X_{n})_{t}\wedge(X_{n})_{\theta}|}-\frac{(N_{n}\cdot(X_{n})_{t\theta})^{2}}{|(X_{n})_{t}\wedge(X_{n})_{\theta}|}
%\end{split}
%\end{equation*}
By (\ref{eq:XzzN})--(\ref{eq:XttN}) we have that
$$
|N\cdot X_{\theta\theta}|\le|W_{\theta}|~\!,\quad
|N\cdot X_{t\theta}|\le|V_{\theta}|~\!,\quad
|N\cdot X_{tt}|\le|V_{t}|+\eps z'_{a}|V|
$$
where $V=V_{\eps ,a}$ and $W=W_{\eps ,a}$ defined according to (\ref{eq:UVW}). Making computations, one obtains that $|W_{\theta}|=x_{a}$, $|V_{t}|\le Cx_{a}$, $|V_{\theta}|\le Cx_{a}$. Moreover $|V|\le Cx_{a}\le C$, $0<z'_{a}\le x_{a}\le 1$ and $|X_{t}\wedge X_{\theta}|\ge \frac{1}{2}x_{a}^{2}$. Hence (\ref{NNN}) is proved and consequently (\ref{eq:immersion}) holds true.
\qed

\noindent
\emph{Proof of Theorem \ref{T:embedded}.} Set $Y_{n}=X_{n,a_{n}}+\varphi_{n,a_{n}}N_{n,a_{n}}$ and
$$
\mathcal{N}_{n}=\{X_{n,a_{n}}(t,\theta)+rx_{a_{n}}(t)N_{n,a_{n}}(t,\theta)~|~(t,\theta)\in\R^{2}~\!,~|r|<r_{0}\}~\!.
$$
We show that for $n$ large enough
\begin{equation}
\label{eq:inside}
Y_{n}(t,\theta)\in \mathcal{N}_{n}\quad\forall(t,\theta)\in\R^{2}~\!.
\end{equation}
Indeed, by (\ref{eq:phi-n-weighted-estimate}) and since $\mu>1$, $0<x_{a_{n}}<1$ and $\eps_{n}\to 0$, one has that $|\varphi_{n,a_{n}}|< r_{0}x_{a_{n}}$ point-wise on $\R^{2}$,
%\begin{equation}
%\label{eq:phinan-estimate}
%|\varphi_{n,a_{n}}|< r_{0}x_{a_{n}}\quad\text{on~}\R^{2}~\!,
%\end{equation}
namely (\ref{eq:inside}). Because of the shape of $Y_{n}$, its image intersect each segment
$$
S_{n}(t,\theta)=\{X_{n,a_{n}}(t,\theta)+rx_{a_{n}}(t)N_{n,a_{n}}(t,\theta)~|~|r|<r_{0}\}
$$
at least once. In order to complete the proof, we show that $Y_{n}$ is an immersion, namely $(Y_{n})_{t}\wedge(Y_{n})_{\theta}\ne 0$ everywhere and the image of $Y_{n}$ intersects each segment $S_{n}(t,\theta)$ at most once. These facts can be simultaneously obtained by checking that
\begin{equation}
\label{positive-normal-projection}
N_{n,a_{n}}\cdot (Y_{n})_{t}\wedge(Y_{n})_{\theta}>0\quad\text{on~~}\R^{2}~\!.
\end{equation}
Exploiting the computations made in the proof of Lemma \ref{L:L-na}, we obtain
$$
N_{n,a_{n}}\cdot (Y_{n})_{t}\wedge(Y_{n})_{\theta}=|X_{t}\wedge X_{\theta}|\left[1+\frac{\varphi_{n,a_{n}}^{2}}{|X_{t}\wedge X_{\theta}|}N\cdot N_{t}\wedge N_{\theta}\right]
$$
where $X=X_{n,a_{n}}$ and $N=N_{n,a_{n}}$. Using (\ref{eq:phi-n-weighted-estimate}) and (\ref{NNN}), we infer that
$$
\varphi_{n,a_{n}}^{2}\frac{|N\cdot N_{t}\wedge N_{\theta}|}{|X_{t}\wedge X_{\theta}|}\le C\eps_{n}^{2\overline\gamma}\quad\text{for $n$ large.}
$$
Hence (\ref{positive-normal-projection}) holds true.
\qed

\noindent
\textbf{Acknowledgements.}
The research of the first author has been partly supported by the PRIN Project 2015KB9WPT ``Variational methods, with applications to problems in mathematical physics and geometry'', by the ERC Advanced Grant 2013 n. 339958 COMPAT - Complex Patterns for Strongly Interacting Dynamical Systems, and by the Gruppo Nazionale per l'Analisi Matematica, la Probabilit\`a e le loro Applicazioni (GNAMPA) of the Istituto Nazionale di Alta Matematica (INdAM). The research of the second author has been partly supported by Fondecyt Grant 1160135.

\end{document}